\documentclass[12pt, leqno]{amsart}
\usepackage{array}
\usepackage{longtable}
\usepackage{url}
\usepackage{enumerate}
\usepackage{amsmath}
\usepackage{amssymb}
\usepackage{amsthm}
\usepackage{amsfonts}
\usepackage{tikz,pgf}
\usepackage{graphicx}
\usepackage{tikz-cd}
\usepackage{mathtools}
\usepackage{tikz}
\usetikzlibrary{arrows}
  \usepackage{ifthen}
 \newcommand{\mymarginpar}[1]{%
    \marginpar{\ifthenelse{\isodd{\arabic{page}}}{\flushleft 
#1}{\flushright #1}}}

\numberwithin{equation}{section}

 \newcommand{\IC}{\mathbb{C}}

 \newcommand{\IN}{\mathbb{N}}                  
                   
 \newcommand{\IQ}{\mathbb{Q}}                  
 \newcommand{\IR}{\mathbb{R}}                  
                   
 \newcommand{\IT}{\mathbb{T}}                  
 \newcommand{\IZ}{\mathbb{Z}}         
                  
\newcommand{\CA}{\mathcal{A}}
\newcommand{\CB}{\mathcal{B}}

\newcommand{\CE}{\mathcal{E}}

\newcommand{\CT}{\mathcal{T}}

\newcommand{\CK}{\mathcal{K}}

 \theoremstyle{plain} %%%%%%%%%%%%%%%%%%%%%%%%%%%%%%%%%
 \newtheorem{Theorem}{Theorem}[section]
 \newtheorem{Lemma}[Theorem]{Lemma}
 \newtheorem{Proposition}[Theorem]{Proposition}
 
 \newtheorem{Corollary}[Theorem]{Corollary}

 \theoremstyle{definition} %%%%%%%%%%%%%%%%%%%%%%%%%%%%
 \newtheorem{Definition}[Theorem]{Definition}
 \newtheorem{Remark}[Theorem]{Remark}
 \newtheorem{Example}[Theorem]{Example}
 
 \newtheorem{Notation}[Theorem]{Notation}

\allowdisplaybreaks[1]

\begin{document}

\title[AF $C^*$-algebras from non AF groupoids]{AF $C^*$-algebras from non AF groupoids}


\author[Mitscher]{Ian Mitscher}
\author[Spielberg]{Jack Spielberg}

\address{School of Mathematical and Statistical Sciences\\Arizona State University\\USA}
\email{ian.mitscher@gmail.com}
\email{jack.spielberg@asu.edu}

\begin{abstract}

We construct ample groupoids from certain categories of paths, and prove that their $C^*$-algebras coincide with the continued fraction AF algebras of Effros and Shen. The proof relies on recent classification results for simple nuclear $C^*$-algebras. The groupoids are not principal. This provides examples of Cartan subalgebras in the continued fraction AF algebras that are isomorphic, but not conjugate, to the standard diagonal subalgebras.

\end{abstract}

\maketitle

\section{Introduction}
\label{section:introduction}

Approximately finite dimensional (AF) $C^*$-algebras were introduced by Bratteli, generalizing the uniformly hyperfinite (UHF) $C^*$-algebras of Glimm. They have become one of the most intensely studied classes of $C^*$-algebras. They are an integral part, and in fact were the beginning, of Elliott's classification program for (simple) separable nuclear $C^*$-algebras. Despite the many remarkable successes in the study of AF algebras, it can be difficult to decide if a given $C^*$-algebra is AF. Namely, while the definition requires that there be an increasing family of finite dimensional $C^*$-subalgebras whose union is dense, finding such a family is not always practical. The known examples of this problem are crossed products (or fixed point algebras) of an action of a finite group on a non AF algebra. In some instances (\cite{kumjian}) the approximating family of finite dimensional subalgebras is explicitly constructed, while in others (\cite{bratkish}, \cite{elpw}) the proof consists of showing that the algebra satisfies the hypotheses of a classification theorem, and then checking that the invariants match those of an AF algebra.

One context where the problem has been completely solved is that of the $C^*$-algebras of directed graphs. It is known that a graph algebra is AF if and only if the graph does not contain a cycle (\cite[Section 5.4]{raeszy}). A bit more detail is known as well. Since AF algebras are always stably finite, the presence of an infinite projection (equivalently, a corner containing a nonunitary isometry) precludes the AF property. It is known that a graph algebra contains an infinite projection if and only if the graph contains a cycle \textit{with an entrance} (\cite[Proposition 5.4]{raeszy}). Moreover, it is also known that all separable AF algebras are Morita equivalent to graph algebras (\cite{drin,tyler}).

A significant generalization of graph algebras was made by Kumjian and Pask with the definition of \textit{higher rank graphs} (\cite{kumpas}). This vastly expanded the family of $C^*$-algebras represented beyond the graph algebras (which are the rank one case). However, the identification of AF algebras among the higher rank graph algebras is much harder. This problem was studied in detail by Evans and Sims (\cite{evanssims}). They gave many positive results. In particular they define \textit{generalized cycles} that play much of the role that cycles play in the case of directed graphs. They proved that a higher rank graph containing a generalized cycle cannot have an AF $C^*$-algebra, and that if a higher rank graph contains a generalized cycle with an entrance then its $C^*$-algebra is infinite. However it is not true that absence of generalized cycles implies that the $C^*$-algebra is AF (\cite{mitscherthesis}, and Example \ref{Ian's first example} below). Moreover, they present an example of a rank two graph whose $C^*$-algebra seems in every way to be AF (even UHF), but were unable to prove this. Similar difficulties occur in other examples (\cite[Subsection 6.1]{evanssims}). Evans and Sims point out that if their example is AF it would give an example of a non principal groupoid whose $C^*$-algebra is AF. (See \cite[Definition III.1.1]{ren} for the definition of AF groupoid, and the proof that AF groupoids are principal.) Moreover this would present an explicit Cartan subalgebra that is not conjugate to the canonical diagonal subalgebra.

In this paper we give examples that are similar, but in a more general setting, and we are able to prove that they are actually AF algebas. One of the factors complicating the study of higher rank graph algebras is that the requirements for a higher rank graph are quite rigid. A directed graph can be thought of as an arbitrary collection of dots and arrows. Thus it is easy to give examples and to tailor them for certain purposes. By contrast, it is difficult to construct explicit higher rank graphs. The second author introduced \textit{categories of paths} (\cite{spiel1}, see also \cite{spiel2}) as a simultaneous generalization of higher rank graphs and of the quasi-lattice ordered groups of Nica. The main motivation, however, was to try to bring some of the flexibility of construction back to the setting of higher rank graphs. A category of paths can be thought of as a directed graph with identifications, but the requirements are much looser than those for higher rank graphs. The notions of generalized cycle and entrance are the same as for higher rank graphs. It is shown in \cite{spiel1} that the presence of a generalized cycle with an entrance implies that the $C^*$-algebra is infinite.  However the presence of a generalized cycle without an entrance is not sufficient to preclude that the $C^*$-algebra is AF (\cite{mitscherthesis}, and Example \ref{Ian's second example} below). Categories of paths determine $C^*$-algebras very much in the same spirit as higher rank graphs, and these can be defined either via an \'etale groupoid or by generators and relations.

Our examples are constructed from categories of paths that are amalgamations of a 2-graph and a 1-graph. We identify the $C^*$-algebras as being Morita equivalent to the continued fraction AF algebras of \cite{effshen}. Moreover there is a natural full corner in our example that is isomorphic to the Effros Shen algebra $A_\theta$, where $\theta \in (0,1) \setminus \IQ$ is arbitrary. The proof relies on the remarkable recent classification results for separable simple nuclear $C^*$-algebras, for which we cite \cite[Theorem D]{tww}. The classification is in terms of the \textit{Elliott invariant}, which consists of the ordered $K$-theory of the algebra, the position of the unit in $K_0^+$, the space of tracial states, and the pairing between the traces and $K_0$. Thus we must compute the Elliott invariant of our examples, and also verify that they satisfy the hypotheses of the classification theorem. Sections \ref{section K-theory of G_i} - \ref{section invariant measure} are devoted to calculating the Elliott invariant, as well as working out the fine structure of our examples. In section \ref{section identifying C*(G)} we verify the hypotheses of the classification theorem and deduce our main theorem.

This provides examples of Cartan subalgebras of AF algebras that are isomorphic, but not conjugate, to the diagonal subalgebra defined by a dense increasing union of finite dimensional subalgebras. We also identify a subalgebra that is the closure of a (nondense) increasing union of finite dimensional subalgebras, which is isomorphic to the whole algebra, and such that the inclusion induces an isomorphism of Elliott invariants.

We briefly describe the rest of the paper. Section \ref{section background} gives a short account of the basic facts about the construction of groupoids and $C^*$-algebras from categories of paths, and the two examples mentioned above. In section \ref{section:one} we describe the categories of paths and groupoids that are the main subject of the paper. The groupoids are inductive limits, and in section \ref{section K-theory of G_i} we compute the ordered $K$-theory of the $C^*$-algebras of the terms in this limit. In section \ref{section K-theory of G} we compute the ordered $K$-theory of the $C^*$-algebra of the whole groupoid and identify it as that of a continued fraction AF algebra. This requires the use of a somewhat technical sequence of partitions refining the compact open subsets of the unit space. (We say that a partition \textit{refines} a set $S$ if $S$ equals the union of a subcollection of the partition.) In section \ref{section invariant measure} we prove that there exists a unique invariant measure on the unit space of the groupoid. In section \ref{section identifying C*(G)} we verify that the groupoid $C^*$-algebra is classifiable and complete the calculation of its Elliott invariant, proving that we have presented the Effros Shen algebras as $C^*$-algebras of non principal groupoids. In section \ref{section subalgebra} we show that our algebra contains a proper copy of itself in the usual form of the Effros Shen algebra. In section \ref{section stability} we investigate the scale of the algebras and give conditions under which it equals the stabilization of the Effros Shen algebra.

We thank the referee for a suggestion that led to section \ref{section stability}.

\section{Background on categories of paths}
\label{section background}

In this section we briefly recall the basic facts about categories of paths and their $C^*$-algebras. This generalizes the case of higher rank graphs. We refer to \cite{spiel1}, \cite{spiel2}.

\begin{Definition}

A \textit{category of paths} is a cancellative small category in which no nonunit has an inverse.

\end{Definition}

Let $\Lambda$ be a category of paths, with unit space $\Lambda^0$. We write $s,r : \Lambda \to \Lambda^0$ for the \textit{source} and \textit{range} maps. For $\mu \in \Lambda$ we write $\mu\Lambda = \{\mu\nu : \nu \in \Lambda, r(\nu) = s(\mu) \}$, and similarly for $\Lambda \mu$. We say that $\nu$ \textit{extends} $\mu$ if $\nu \in \mu \Lambda$, and in this case we call $\mu$ a \textit{prefix} of $\nu$. The set of prefixes of an element $\nu \in \Lambda$ is denoted $[\nu]$. For two elements $\mu,\nu \in \Lambda$ the \textit{common extensions} of $\mu$ and $\nu$ are the elements of $\mu\Lambda \cap \nu\Lambda$. We write $\mu \Cap \nu$ if $\mu\Lambda \cap \nu\Lambda \not= \varnothing$, and we say that $\mu$ \textit{meets} $\nu$; otherwise we write $\mu \perp \nu$, and we say that $\mu$ and $\nu$ are \textit{disjoint}. There is an important simplifying assumption that was given for higher rank graphs in \cite{raesimyee}.

\begin{Definition}

The category of paths $\Lambda$ is \textit{finitely aligned} if for all $\mu,\nu \in \Lambda$ there exists a (possibly empty) finite set $F \subseteq \Lambda$ such that $\mu\Lambda \cap \nu\Lambda = \bigcup_{\eta \in F} \eta \Lambda$.

\end{Definition}

In the finitely aligned case there is a unique minimal such set $F$, denoted $\mu \vee \nu$, the set of \textit{minimal common extensions} of $\mu$ and $\nu$.

\begin{Definition}

Let $\Lambda$ be a finitely aligned category of paths. A subset $x \subseteq \Lambda$ is

\begin{itemize}

\item[] \textit{directed} if for all $\mu,\nu \in x$, $x \cap \mu\Lambda \cap \nu\Lambda \not= \varnothing$,

\item[] \textit{hereditary} if for all $\mu \in x$, $[\mu] \subseteq x$.

\end{itemize}

\end{Definition}

All elements of a directed hereditary set have the same range, denoted $r(x)$. We let $\Lambda^*$ denote the set of all directed hereditary subsets of $\Lambda$, and for $u \in \Lambda^0$ we write $u \Lambda^* = \{ x \in \Lambda^* : r(x) = u \}$. For $\mu \in \Lambda$ we write $Z(\mu) = \{ x \in \Lambda^* : \mu \in x \}$. We think of $Z(\mu)$ as a \textit{generalized cylinder set}.

\begin{Theorem}

Let $\Lambda$ be a finitely aligned category of paths. The collection $\{ Z(\mu) : \mu \in \Lambda \}$ generates a totally disconnected locally compact Hausdorff topology on $\Lambda^*$, and $Z(\mu)$ is compact open in this topology. The collection $\CB^* = \{ Z(\lambda) \setminus \cup_{i=1}^k Z(\mu_i) : \lambda, \mu_1, \ldots, \mu_k \in \Lambda, \}$ is a base (of compact open sets) for this topology.

\end{Theorem}

We let $\Lambda^{**}$ denote the maximal elements of $\Lambda^*$ with respect to inclusion, and put $\partial \Lambda = \overline{\Lambda^{**}}$, the \textit{boundary} of $\Lambda$. Left cancellation in $\Lambda$ implies that for $\mu \in \Lambda$, the map $\nu \in s(\mu)\Lambda \mapsto \mu\nu \in \mu\Lambda$ is bijective. It extends to a partial homeomorphism of $\Lambda^* : x \in s(\mu) \Lambda^* \mapsto \mu x = \cup \{[\mu\nu] : \nu \in x \}$. Moreover, $\mu x \in \partial \Lambda$ if and only if $x \in \partial \Lambda$.

We define a groupoid with unit space $\Lambda^*$ as follows. Let $\Lambda * \Lambda * \Lambda^* = \{ (\mu,\nu,x) \in \Lambda \times \Lambda \times \Lambda^* : s(\mu) = s(\nu) = r(x) \}$. Let $\sim$ be the equivalence relation on $\Lambda * \Lambda * \Lambda^*$ generated by $(\mu,\nu, \eta x) \sim (\mu\eta,\nu\eta,x)$, and define $G(\Lambda) = \Lambda * \Lambda * \Lambda^* / \sim$. The groupoid operations are given by

\begin{enumerate}[(i)]

\item $s([\mu,\nu,x]) = \nu x$, $r([\mu,\nu,x]) = \mu x$

\item $[\mu,\nu,x]^{-1} = [\nu,\mu,x]$

\item \label{multiplication in G(Lambda)} $[\mu,\nu,x] \cdot [\nu,\eta,x] = [\mu,\eta,x]$.

\end{enumerate}
More precisely, we use \cite[Lemma 4.12]{spiel1}: if $\mu,\mu' \in \Lambda$ and $x,x' \in \Lambda^*$ are such that $\mu x = \mu' x'$ then there exist $\eta,\eta' \in \Lambda$ and $y \in \Lambda^*$ such that $x = \eta y$, $x' = \eta' y$, and $\mu \eta = \mu' \eta'$. Then the composable pairs are given by $G^{(2)} = \{([\mu,\nu,x],[\xi,\eta,y]) : \nu x = \xi y \}$, and multiplication is given as follows. Let $([\mu,\nu,x],[\xi,\eta,y]) \in G^{(2)}$. Then $\nu x = \xi y$, so the lemma cited above yields $z \in \Lambda^*$ and $\delta,\epsilon \in \Lambda$ such that $x = \delta z$, $y = \epsilon z$, and $\nu \delta = \xi \epsilon$. Then
\[
[\mu,\nu,x] \cdot [\xi,\eta,y]
= [\mu\delta,\nu\delta,z] \cdot [\xi\epsilon,\eta\epsilon,z] = [\mu\delta,\eta\epsilon,z].
\]
The unit space $G(\Lambda)^{(0)} = \{[r(x),r(x),x] : x \in \Lambda^*\}$ can be identified with $\Lambda^*$. For a subset $E \subseteq s(\nu)\Lambda^*$ we write $[\mu,\nu,E] := \{[\mu,\nu,x] : x \in E \}$.

\begin{Theorem}

The collection of sets $[\mu,\nu,E]$ for $E \subseteq \Lambda^*$ compact open is a base for a topology making $G(\Lambda)$ into an ample Hausdorff (\'etale) groupoid.

\end{Theorem}

The boundary $\partial \Lambda$ is a closed invariant subset of $G(\Lambda)^{(0)} \equiv \Lambda^*$.

\begin{Definition}

The \textit{Toeplitz algebra} of $\Lambda$ is $\CT C^*(\Lambda) := C^*(G(\Lambda))$. The \textit{Cuntz-Krieger algebra} of $\Lambda$ is $C^*(\Lambda) := C^*(G(\Lambda)|_{\partial \Lambda})$.

\end{Definition}

Generalized cycles and entrances are defined exactly as for higher rank graphs (\cite{evanssims}).

\begin{Definition} \label{def generalized cycle}

Let $\Lambda$ be a category of paths. A \textit{generalized cycle} in $\Lambda$ is an ordered pair $(\mu,\nu) \in \Lambda \times \Lambda$ such that $\mu \not= \nu$, and

\begin{enumerate}[(i)]

\item $s(\mu) = s(\nu)$

\item $r(\mu) = r(\nu)$

\item for all $\eta \in \mu\Lambda$, $\nu \Cap \eta$.

\end{enumerate}

The generalized cycle $(\mu,\nu)$ \textit{has an entrance} if $(\nu,\mu)$ is not a generalized cycle.

\end{Definition}

It is proved in \cite{spiel1} that if a finitely aligned category of paths contains a generalized cycle with an entrance then its Cuntz-Krieger algebra is infinite. In particular,

\begin{Theorem}

If the finitely aligned category of paths $\Lambda$ contains a generalized cycle with an entrance, then $C^*(\Lambda)$ is not AF.

\end{Theorem}

It is proved in \cite{evanssims} that the $C^*$-algebra of a higher rank graph that contains a generalized cycle without an entrance is not AF. However, this is not generally true for categories of paths. The following example is due to the first author (\cite{mitscherthesis}).

\begin{Example} \label{Ian's second example}

Let $\Lambda$ be the category of paths given by the directed graph
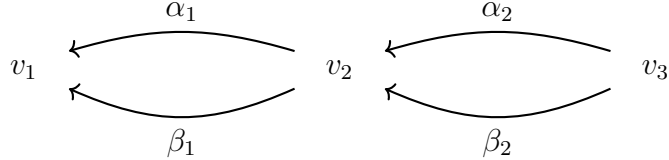
\[
\begin{tikzpicture}[yscale = 5, xscale = 3]
	\draw[thick, <-] (0,-.05) to[bend right=15] node[below]{$\beta_1$} (1,-.05);
	\draw[thick, <-] (1.4,-.05) to[bend right=15] node[below]{$\beta_2$} (2.4,-.05);
	\draw[thick, <-] (0,.05) to[bend left=10] node[above]{$\alpha_1$} (1,.05);
	\draw[thick, <-] (1.4,.05) to[bend left=10] node[above]{$\alpha_2$} (2.4,.05);
		
	\node (bleft) at (-.2,0) [rectangle] {$v_1$};
	\node (bleft) at (1.2,0) [rectangle] {$v_2$};
	\node (bleft) at (2.6,0) [rectangle] {$v_3$};
	%\node [below=1cm, align=flush center] at (0,-.1)
	%{Figure 5.3};
\end{tikzpicture}
\]
with commuting relations $\alpha_1\beta_2 = \beta_1\alpha_2$ and $\alpha_1\alpha_2 = \beta_1\beta_2$. Explicitly, 
\begin{align*}
\Lambda &= \{v_1, v_2, v_3, \alpha_1, \alpha_2, \beta_1, \beta_2, \alpha_1\alpha_2, \alpha_1\beta_2\} \\
\noalign{where} \Lambda^0 &= \{ v_1, v_2, v_3\}\\
\noalign{and} \Lambda^2 &= \{ (\alpha_1, \alpha_2), (\alpha_1, \beta_2), (\beta_1, \alpha_2), (\beta_1, \beta_2)\}
\end{align*}
with the obvious range and source maps, and ranges for multiplication. It is clear that this is a small category, without inverses since the commuting relations preserve length, and there is no (non-trivial) cancellation to check. Note that $(\alpha_1, \beta_1)$ is a generalized cycle without an entrance. Since the category is finite, $G(\Lambda)$ is finite, and hence $C^*(\Lambda)$ is finite dimensional.

\end{Example}

We also note that even for higher rank graphs, the absence of generalized cycles is not sufficient for the $C^*$-algebra to be AF. We give the following example, due to the first author (\cite{mitscherthesis}).

\begin{Example} \label{Ian's first example}

Let $\Lambda$ be the 2-graph which has skeleton
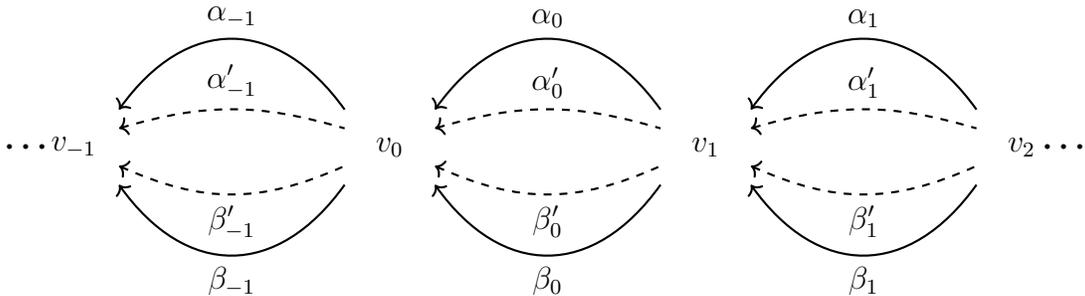
\[
\begin{tikzpicture}[yscale = 5, xscale = 3]
		\draw[thick, dashed, <-] (0,-.05) to[bend right=15] node[below]{$\beta_{-1}'$} (1,-.05);
		\draw[thick, dashed, <-] (1.4,-.05) to[bend right=15] node[below]{$\beta_0'$} (2.4,-.05);
		\draw[thick, dashed, <-] (2.8,-.05) to[bend right=15] node[below]{$\beta_1'$} (3.8,-.05);
		\draw[thick, <-] (0,-.1) to[bend right=40] node[below]{$\beta_{-1}$} (1,-.1);
		\draw[thick, <-] (1.4,-.1) to[bend right=40] node[below]{$\beta_0$} (2.4,-.1);
		\draw[thick, <-] (2.8,-.1) to[bend right=40] node[below]{$\beta_1$} (3.8,-.1);
		\draw[thick, dashed, <-] (0,.05) to[bend left=10] node[above]{$\alpha_{-1}'$} (1,.05);
		\draw[thick, dashed, <-] (1.4,.05) to[bend left=10] node[above]{$\alpha_0'$} (2.4,.05);
		\draw[thick, dashed, <-] (2.8,.05) to[bend left=10] node[above]{$\alpha_1'$} (3.8,.05);
		\draw[thick, <-] (0,.1) to[bend left=40] node[above]{$\alpha_{-1}$} (1,.1);
		\draw[thick, <-] (1.4,.1) to[bend left=40] node[above]{$\alpha_0$} (2.4,.1);
		\draw[thick, <-] (2.8,.1) to[bend left=40] node[above]{$\alpha_1$} (3.8,.1);
		
		\node (b40) at (4.2,0) [rectangle] {$\boldsymbol{\cdots}$};
		\node (b50) at (-.4,0) [rectangle] {$\boldsymbol{\cdots}$};
		%\node [below=1cm, align=flush center] at (0,-.2)
		%{Figure 5.2};
		
		\node (bleft) at (-.2,0) [rectangle] {$v_{-1}$};
		\node (bleft) at (1.2,0) [rectangle] {$v_0$};
		\node (bleft) at (2.6,0) [rectangle] {$v_1$};
		\node (bleft) at (4,0) [rectangle] {$v_2$};
\end{tikzpicture}
\]
and factorization rules $\sigma_i\tau_{i+1}' = \tau_i'\sigma_{i+1}$ where $\sigma, \tau\in \{\alpha, \beta\}$. As in \cite[Section 2.4]{evanssims}, this defines a 2-graph. We will use the convention that $d(\alpha_i) = d(\beta_i) = (1,0)$ and $d(\alpha_i') = d(\beta_i') = (0,1)$. By \cite[Proposition 3.16]{evans2} (see also \cite[Corollary 5.1]{evans}) and the remarks following \cite[Figure 4.3]{evans}, $K_1(C^*(\Lambda))\cong \mathbb{Z}[\frac{1}{2}]$ so that $C^*(\Lambda)$ is not an AF-algebra. We will show that $\Lambda$ has no generalized cycles.

To facilitate the argument, we will drop subscripts when the range is understood. Thus we may write $\alpha_i \beta_{i+1}' = \alpha_i \beta' = v_i \alpha \beta'$. We note the following notational device. Let $\lambda \in v_i \Lambda$ with $d(\lambda) = (m,0)$. Write $\lambda = \lambda(i)_i \lambda(i+1)_{i+1} \cdots \lambda(i+m-1)_{i+m-1}$, where $\lambda(j) \in \{\alpha,\beta\}$ for $i \le j < i+m$. Then the factorization rules give
\begin{align*}
\lambda \alpha_{i+m}'
&= \lambda(i)_i \lambda(i+1)_{i+1} \cdots \lambda(i+m-1)_{i+m-1} \alpha_{i+m}' \\
&= \lambda(i)_i \lambda(i+1)_{i+1} \cdots \lambda(i+m-2)_{i+m-2} \alpha_{i+m-1}' \lambda(i+m-1)_{i+m} \\
&= \cdots \\
&= \alpha_i' \lambda(i)_{i+1} \cdots \lambda(i+m-1)_{i+m} \\
&= \alpha' \lambda.
\end{align*}
Similarly, $\lambda \beta' = \beta' \lambda$, and in general, if $\lambda' \in v_{i+m} \Lambda$ with $d(\lambda') = (0,n)$ we have $\lambda\lambda' = \lambda' \lambda$.

Now fix $\mu\ne\nu\in\Lambda$ with $r(\mu)=r(\nu)$ and $s(\mu)=s(\nu)$. Write $\mu = \sigma\sigma'$, $\nu = \tau\tau'$ where $d(\sigma),d(\tau)\in\mathbb{N}\times\{0\}$ and $d(\sigma'),d(\tau')\in\{0\}\times\mathbb{N}$. First assume $d(\sigma)\ge d(\tau)$ and suppose that $\sigma$ does not extend $\tau$. Then there is some $j$ such that $\sigma_j\ne \tau_j$. This implies that $\mu$ and $\nu$ have no common extension, since if $\lambda$ were such an extension, by factoring $\lambda = \sigma\lambda'$ we would have $\sigma_j = \lambda_j = \tau_j$.

Now suppose $\sigma$ extends $\tau$ and consider the case where $d(\sigma) > d(\tau)$. Let $v_j = s(\tau)$ and suppose $\sigma_j = \alpha_j$. Then we see that $\mu, \nu\beta$ have no common extension by the same reasoning as above, since $\nu\beta = \tau\beta_j\tau'$ while $\mu_j = \alpha_j$. Similarly, if $\sigma_j = \beta_j$, then $\mu, \nu\alpha$ have no common extension. If $d(\sigma) = d(\tau)$, then $\sigma = \tau$ and hence $\sigma'\ne \tau'$ so there is some $j$ with $\sigma'_j\ne \tau_j'$, and reasoning as above shows this precludes the existence of a common extension of $\mu$ and $\nu$.

If $d(\sigma) < d(\tau)$, we may write $\mu = \sigma'\sigma$ and $\nu = \tau'\tau$. Then the previous argument applies to $\sigma'$ and $\tau'$. Thus in all cases we can find $\eta\in\Lambda$ with $\mu, \nu\eta$ having no common extension, so $(\mu,\nu)$ is not a generalized cycle.
\end{Example}

\section{The examples}
\label{section:one}

We will build a family of categories of paths as amalgamations of a 2-graph and a 1-graph (\cite[Section 11]{spiel1}). We begin with the following directed graph $E$:
\[
\begin{tikzpicture}[scale = 1.5]
\draw[<-,thick] (0,0) to node[above]{$a_1$} (1,0);
\draw[<-,thick] (1.4,0) to node[above]{$a_2$} (2.4,0);
\draw[<-,thick] (2.8,0) to node[above]{$a_3$} (3.8,0);

\node (bleft) at (-.2,0) [rectangle] {$w_1$};
\node (bleft) at (1.2,0) [rectangle] {$w_2$};
\node (bleft) at (2.6,0) [rectangle] {$w_3$};
\node (bleft) at (4,0) [rectangle] {$w_4$};
\node (b40) at (4.5,0) [rectangle] {$\boldsymbol{\cdots}$};
\end{tikzpicture}
\]
Thus $E^0 = \{w_1,w_2, \ldots\}$, $E^1 = \{a_1,a_2,\ldots\}$, and $s(a_i) = w_{i+1}$, $r(a_i) = w_i$. Let $E^*$ be the set of finite paths in $E$, that is, the 1-graph associated with $E$. Let $f : \IN^2 \to \IN$ be given by $f(n) = n_1 + n_2$. Using \cite[Definition 1.9]{kumpas}, set $\Lambda_1 = f^*(E^*)$; $\Lambda_1$ is a 2-graph. Then $\Lambda_1 = \{ (\mu,n) \in E^* \times \IN^2 : |\mu| = f(n) \}$. For $i \ge 1$ we let $\alpha_i = (a_i,e_1)$ and $\beta_i = (a_i,e_2)$; these are the edges in $\Lambda_1$:
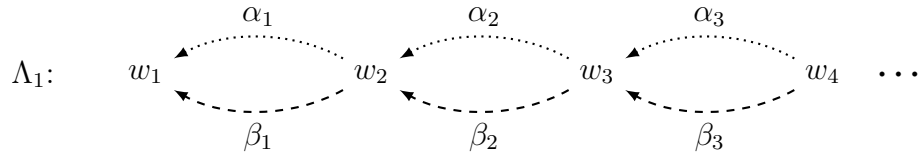
\[
\begin{tikzpicture}[xscale=3,yscale=2]

\node (00) at (.5,0) [rectangle] {$\Lambda_1$:};

\node (10) at (1,0) [rectangle] {$w_1$};
\node (20) at (2,0) [rectangle] {$w_2$};
\node (30) at (3,0) [rectangle] {$w_3$};
\node (40) at (4,0) [rectangle] {$w_4$};
\node (4half0) at (4.35,0) [rectangle] {$\boldsymbol{\cdots}$};

\draw[-latex,thick,dotted] (20) to[bend right=40] node[above] {$\alpha_1$} (10);
\draw[-latex,thick,dotted] (30) to[bend right=40] node[above] {$\alpha_2$} (20);
\draw[-latex,thick,dotted] (40) to[bend right=40] node[above] {$\alpha_3$} (30);

\draw[-latex,thick,dashed] (20) to[bend left=40] node[below] {$\beta_1$} (10);
\draw[-latex,thick,dashed] (30) to[bend left=40] node[below] {$\beta_2$} (20);
\draw[-latex,thick,dashed] (40) to[bend left=40] node[below] {$\beta_3$} (30);

\end{tikzpicture}
\]
The commuting squares in $\Lambda_1$ are $\alpha_i \beta_{i+1} = \beta_i \alpha_{i+1}$. We may write a typical element of $\Lambda_1$ as a product of edges: $(a_i a_{i+1} \cdots a_j,n) = \alpha_i \alpha_{i+1} \cdots \alpha_{i + n_1 - 1} \beta_{i + n_1} \beta_{i + n_1 + 1} \cdots \beta_{i + n_1 + n_2 - 1}$. We will often write this in the form $\alpha^{n_1} \beta^{n_2}$ if the source (and range) are understood (otherwise we might write, for example, $w_i \alpha^{n_1} \beta^{n_2}$). By the factorization rule in $\Lambda_1$ the edges may be ``permuted'' as desired; thus, for example, $\alpha_i \alpha_{i+1} \beta_{i+2} \beta_{i+3} = \beta_i \alpha_{i+1} \beta_{i+2} \alpha_{i+3}$, etc.

Now choose nonnegative integers $k_1$, $k_2$, $\ldots$, with $k_i > 0$ for infinitely many $i$, and let $\Lambda_2$ be the 1-graph with vertex set $\Lambda_2^0 = \{u_1,u_2,\ldots\}$ and edge set $\Lambda_2^1 = \{\gamma_i^{(j)} : i \ge 1, 1 \le j \le k_i \}$, with $s(\gamma_i^{(j)}) = u_{i+1}$, $r(\gamma_i^{(j)})= u_i$. (Note that if $k_i = 0$ then $u_i \Lambda_2 u_{i+1} = \varnothing$.)
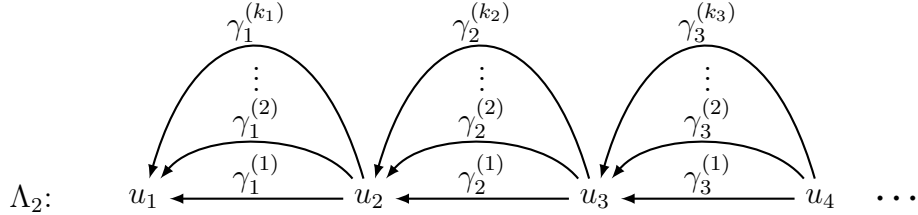
\[
\begin{tikzpicture}[xscale=3,yscale=2]

\node (00) at (.5,0) [rectangle] {$\Lambda_2$:};

\node (10) at (1,0) [rectangle] {$u_1$};
\node (20) at (2,0) [rectangle] {$u_2$};
\node (30) at (3,0) [rectangle] {$u_3$};
\node (40) at (4,0) [rectangle] {$u_4$};

\node (4half0) at (4.35,0) [rectangle] {$\boldsymbol{\cdots}$};
\node (1half1) at (1.5,.85) [rectangle] {$\boldsymbol{\vdots}$};
\node (2half1) at (2.5,.85) [rectangle] {$\boldsymbol{\vdots}$};
\node (23half1) at (3.5,.85) [rectangle] {$\boldsymbol{\vdots}$};

\draw[-latex,thick] (20) -- (10) node[pos=0.5,inner sep=0.5pt,above=1pt] {$\gamma^{(1)}_1$};
\draw[-latex,thick] (20) .. controls (1.75,.45) and (1.25,.45) .. (10) node[pos=0.5,inner sep=0.5pt,above=1pt] {$\gamma^{(2)}_1$};
\draw[-latex,thick] (20) .. controls (1.75,1.3) and (1.25,1.3) .. (10) node[pos=0.5,inner sep=0.5pt,above=1pt] {$\gamma^{(k_1)}_1$};

\draw[-latex,thick] (30) -- (20) node[pos=0.5,inner sep=0.5pt,above=1pt] {$\gamma^{(1)}_2$};
\draw[-latex,thick] (30) .. controls (2.75,.45) and (2.25,.45) .. (20) node[pos=0.5,inner sep=0.5pt,above=1pt] {$\gamma^{(2)}_2$};
\draw[-latex,thick] (30) .. controls (2.75,1.3) and (2.25,1.3) .. (20) node[pos=0.5,inner sep=0.5pt,above=1pt] {$\gamma^{(k_2)}_2$};

\draw[-latex,thick] (40) -- (30) node[pos=0.5,inner sep=0.5pt,above=1pt] {$\gamma^{(1)}_3$};
\draw[-latex,thick] (40) .. controls (3.75,.45) and (3.25,.45) .. (30) node[pos=0.5,inner sep=0.5pt,above=1pt] {$\gamma^{(2)}_3$};
\draw[-latex,thick] (40) .. controls (3.75,1.3) and (3.25,1.3) .. (30) node[pos=0.5,inner sep=0.5pt,above=1pt] {$\gamma^{(k_3)}_3$};

\end{tikzpicture}
\]

\begin{Definition} \label{def Lambda}

Let $\sim$ be the equivalence relation on $\Lambda_1^0 \sqcup \Lambda_2^0$ generated by $w_i \sim u_i$ for $i \ge 1$. Let $\Lambda$ be the category of paths given by the amalgamation of $\Lambda_1$ and $\Lambda_2$ over $\sim$, as in \cite[Definition 11.3]{spiel1}. Let $v_i = [u_i] = [w_i]$ for $i \ge 1$. We will make a small notational abuse to write $v_i \Lambda_1$ for the set $\{v_i \alpha^m \beta^n : m,n \ge 0\} \subseteq \Lambda$, and identify it with $w_i \Lambda_1$ (and similarly for $v_i \Lambda_2$).

\end{Definition}

We may picture $\Lambda$ as follows:
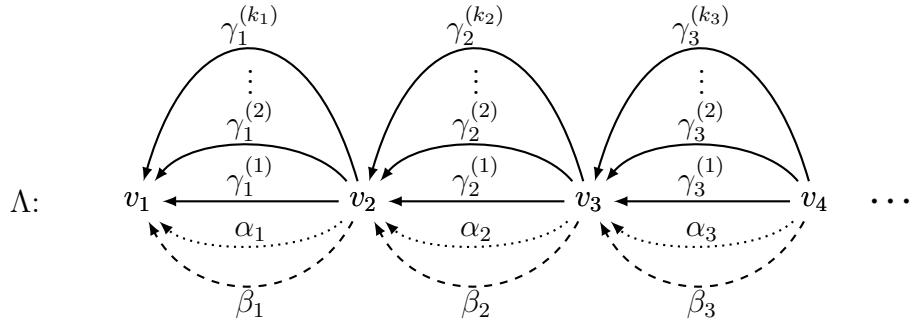
\[
\begin{tikzpicture}[xscale=3,yscale=2]

\node (00) at (.5,0) [rectangle] {$\Lambda$:};

\node (10) at (1,0) [rectangle] {$v_1$};
\node (20) at (2,0) [rectangle] {$v_2$};
\node (30) at (3,0) [rectangle] {$v_3$};
\node (40) at (4,0) [rectangle] {$v_4$};

\node (4half0) at (4.35,0) [rectangle] {$\boldsymbol{\cdots}$};
\node (1half1) at (1.5,.85) [rectangle] {$\boldsymbol{\vdots}$};
\node (2half1) at (2.5,.85) [rectangle] {$\boldsymbol{\vdots}$};
\node (23half1) at (3.5,.85) [rectangle] {$\boldsymbol{\vdots}$};

\draw[-latex,thick] (20) -- (10) node[pos=0.5,inner sep=0.5pt,above=1pt] {$\gamma^{(1)}_1$};
\draw[-latex,thick] (20) .. controls (1.75,.45) and (1.25,.45) .. (10) node[pos=0.5,inner sep=0.5pt,above=1pt] {$\gamma^{(2)}_1$};
\draw[-latex,thick] (20) .. controls (1.75,1.3) and (1.25,1.3) .. (10) node[pos=0.5,inner sep=0.5pt,above=1pt] {$\gamma^{(k_1)}_1$};

\draw[-latex,thick] (30) -- (20) node[pos=0.5,inner sep=0.5pt,above=1pt] {$\gamma^{(1)}_2$};
\draw[-latex,thick] (30) .. controls (2.75,.45) and (2.25,.45) .. (20) node[pos=0.5,inner sep=0.5pt,above=1pt] {$\gamma^{(2)}_2$};
\draw[-latex,thick] (30) .. controls (2.75,1.3) and (2.25,1.3) .. (20) node[pos=0.5,inner sep=0.5pt,above=1pt] {$\gamma^{(k_2)}_2$};

\draw[-latex,thick] (40) -- (30) node[pos=0.5,inner sep=0.5pt,above=1pt] {$\gamma^{(1)}_3$};
\draw[-latex,thick] (40) .. controls (3.75,.45) and (3.25,.45) .. (30) node[pos=0.5,inner sep=0.5pt,above=1pt] {$\gamma^{(2)}_3$};
\draw[-latex,thick] (40) .. controls (3.75,1.3) and (3.25,1.3) .. (30) node[pos=0.5,inner sep=0.5pt,above=1pt] {$\gamma^{(k_3)}_3$};

%%%%%%%%%%%%%%
%%%%%%%%%%%%%%

\node (10) at (1,0) [rectangle] {$v_1$};
\node (20) at (2,0) [rectangle] {$v_2$};
\node (30) at (3,0) [rectangle] {$v_3$};
\node (40) at (4,0) [rectangle] {$v_4$};

\draw[-latex,thick,dotted] (20) .. controls (1.75,-.35) and (1.25,-.35) .. (10) node[pos=0.5,inner sep=0.5pt,above=1pt] {$\alpha_1$};
\draw[-latex,thick,dotted] (30) .. controls (2.75,-.35) and (2.25,-.35) .. (20) node[pos=0.5,inner sep=0.5pt,above=1pt] {$\alpha_2$};
\draw[-latex,thick,dotted] (40) .. controls (3.75,-.35) and (3.25,-.35) .. (30) node[pos=0.5,inner sep=0.5pt,above=1pt] {$\alpha_3$};

\draw[-latex,thick,dashed] (20) .. controls (1.75,-.7) and (1.25,-.7) .. (10) node[pos=0.5,inner sep=0.5pt,below=1pt] {$\beta_1$};
\draw[-latex,thick,dashed] (30) .. controls (2.75,-.7) and (2.25,-.7) .. (20) node[pos=0.5,inner sep=0.5pt,below=1pt] {$\beta_2$};
\draw[-latex,thick,dashed] (40) .. controls (3.75,-.7) and (3.25,-.7) .. (30) node[pos=0.5,inner sep=0.5pt,below=1pt] {$\beta_3$};

\end{tikzpicture}
\]
We remark that the sequence $(k_i)$ is a necessary part of the definition of $\Lambda$. We will use it when working with $\Lambda$ and other objects constructed from $\Lambda$ (such as $G(\Lambda)$), understanding that a choice of this sequence has been made, even if this was not explicit.

\begin{Remark}

It is fruitful to think of $\Lambda_1$ as the base of our construction, having a regular structure, while $\Lambda_2$ represents a sporadic introduction of ``impurities'' into that regular structure. (For example, in Theorem \ref{thm stability} we will show that $C^*(\Lambda)$ is stable if and only if the terms of the sequence $(k_i)$ are rarely nonzero.)  A similar construction could be made with $\Lambda_1$ changed to an analogous higher rank graph with rank larger than two. Thus we trust that the reader will not be misled by the apparent mixup of rank and subscript.

\end{Remark}

\begin{Lemma} \label{lemma:normal form}

Let $\lambda \in \Lambda \setminus \Lambda^0$. There exist unique $m \ge 1$ and elements $\theta_1$, $\ldots$, $\theta_m \in \Lambda_1 \cup \Lambda_2$ such that
\begin{itemize}

\item[$\boldsymbol{\cdot}$] $\theta_i \not\in \Lambda_1^0 \cup \Lambda_2^0$

\item[$\boldsymbol{\cdot}$] $s(\theta_i) \sim r(\theta_{i+1})$

\item[$\boldsymbol{\cdot}$] $s(\theta_i) \not= r(\theta_{i+1})$

\item[$\boldsymbol{\cdot}$] $\lambda = \theta_1 \cdots \theta_m$.

\end{itemize}

\noindent
The decomposition $\lambda = \theta_1 \cdots \theta_m$ is called the {\em normal form} of $\lambda$.

\end{Lemma}

\begin{proof}
This follows from \cite[Lemma 11.2]{spiel1}.
\end{proof}

Thus $\theta_i \in \Lambda_1$ if and only if $\theta_{i+1} \in \Lambda_2$. We write $|\lambda| = \sum_{i=1}^m |\theta_i|$, where for $\theta \in \Lambda_1 \cup \Lambda_2$ we let $|\theta|$ denote the number of edges in $\theta$. If $r(\lambda) = v_j$ and $|\lambda| = n$, we may write $\lambda = \lambda_j \lambda_{j+1} \ldots \lambda_{j + n - 1}$, where $\lambda_i \in \{\alpha_i, \beta_i\} \cup \{\gamma_i^{(t)} : 1 \le t \le k_i\}$. Of course this representation of $\lambda$ might not be unique, as a factor $\theta_i \in \Lambda_1$ may be written in several ways as a path.

\begin{Lemma} \label{lemma:common extension}

Let $\lambda,\mu \in \Lambda$. Write $\lambda = \theta_1 \cdots \theta_m$, $\mu = \phi_1 \cdots \phi_n$ in normal form. Then $\lambda$ and $\mu$ have a common extension if and only if $\theta_i = \phi_i$ for $i < \min\{m,n\}$, and one of the following conditions holds:

\begin{enumerate}[(1)]

\item \label{lemma:common extension 1} $m \not= n$, and if, say, $m < n$, then $\phi_m$ extends $\theta_m$. (Note that in this case, $\mu$ extends $\lambda$.)

\item \label{lemma:common extension 2} $m = n$, and $\theta_m$ and $\phi_m$ have a common extension in $\Lambda_j$ (where $j$ is such that $\theta_m,\phi_m \in \Lambda_j$). (Note that in this case, if $j = 2$ then one of $\lambda,\mu$ extends the other.)

\end{enumerate}

\end{Lemma}

\begin{proof}
This follows from \cite[Lemma 11.4]{spiel1}.
\end{proof}

\begin{Corollary} \label{cor:eqinlam1c}

Let $\lambda,\mu \in \Lambda$ have a common extension, and suppose that neither extends the other. Then $\lambda$ and $\mu$ must be in situation (2) of Lemma \ref{lemma:common extension}, and with notation as in Lemma \ref{lemma:common extension} we have that $j = 1$ and that $d(\theta_m)$, $d(\phi_m)$ are not comparable (in $\IN^2$).

\end{Corollary}

It follows from Corollary \ref{cor:eqinlam1c} that $\Lambda$ is finitely aligned. In fact, since both $\Lambda_1$ and $\Lambda_2$ are \textit{singly aligned} (or \textit{right LCM}) then so is $\Lambda$; this means that if $\lambda$ and $\mu$ have a common extension then they have a unique minimal common extension.

\begin{Proposition} \label{prop no generalized cycle}

$\Lambda$ contains no generalized cycles.

\end{Proposition}

\begin{proof}
Let $\mu,\nu \in \Lambda$ with $r(\mu) = r(\nu)$, $s(\mu) = s(\nu)$ and $\mu \not= \nu$. Then neither of $\mu,\nu$ can extend the other. If $\mu \perp \nu$ then $(\mu,\nu)$ is not a generalized cycle. So suppose that $\mu \Cap \nu$. By Corollary \ref{cor:eqinlam1c}, $\mu$ and $\nu$ have normal forms $\mu = \theta_1 \cdots \theta_{m-1} \phi$ and $\nu = \theta_1 \cdots \theta_{m-1} \phi'$, where $\phi,\phi' \in \Lambda_1$ and $d(\phi),d(\phi')$ are not comparable in $\IN^2$. Letting $\phi = \alpha^p \beta^q$ and $\phi' = \alpha^{p'} \beta^{q'}$, we may assume without loss of generality that $p < p'$ and $q > q'$. Choose $i > j$, where $v_j = s(\mu)$, such that $k_i \not= 0$. Let $\eta = \beta^{i - |\mu| - 1} \gamma_i^{(1)}$. Then Lemma \ref{lemma:common extension} implies that $\mu\eta \perp \nu$, and hence $(\mu,\nu)$ is not a generalized cycle.
\end{proof}

We will describe the elements of $\Lambda^*$ explicitly. First note that the finite directed hereditary subsets of $\Lambda$ are precisely those that contain a maximal element, and these are in one-to-one correspondence with the elements of $\Lambda$, via $\lambda \leftrightarrow [\lambda]$.

\begin{Lemma} \label{lemma:infinite directed hereditary sets in Lambda one}

The infinite directed hereditary subsets of $\Lambda_1$ are the following: for $\ell \ge 1$,

\begin{itemize}

\item[$\boldsymbol{\cdot}$] $w_\ell \alpha^m \beta^\infty := \{ w_\ell \alpha^i \beta^j : 0 \le i \le m, 0 \le j \}$, $m = 0$, 1, 2, $\ldots$

\item[$\boldsymbol{\cdot}$] $w_\ell \alpha^\infty \beta^n := \{ w_\ell \alpha^i \beta^j : 0 \le i, 0 \le j  \le n \}$, $n = 0$, 1, 2, $\ldots$

\item[$\boldsymbol{\cdot}$] $w_\ell \alpha^\infty \beta^\infty := w_\ell \Lambda_1$.

\end{itemize}

\end{Lemma}

\begin{proof}
This follows from the commutation relations in $\Lambda_1$.
\end{proof}

\begin{Notation} \label{notation Lambda_1^*}

We write $v_\ell \Lambda_1^\infty = \{v_\ell \alpha^m \beta^n : m + n = \infty \}$ for the infinite elements of $\Lambda_1^*$ with range $v_\ell$. The infinite elements of $v_\ell \Lambda_2^*$ can be identified with the infinite path space $v_\ell \Lambda_2^\infty$ of the directed graph $\Lambda_2$. (We note that $\Lambda_2^\infty = \varnothing$ unless $k_i > 0$ eventually.)

\end{Notation}

\begin{Theorem} \label{theorem:infinite directed hereditary sets in Lambda}

Let $x$ be an infinite element of $\Lambda^*$. Then $x$ has one of the following forms:

\begin{enumerate}

\item \label{idhsil one} there is an infinite sequence $\theta_1$, $\theta_2$, $\ldots$ such that $\theta_1 \theta_2 \cdots \theta_m$ is in normal form for all $m$, and $x = \cup_{m=1}^\infty [\theta_1 \cdots \theta_m]$;

\item \label{idhsil two} there is a finite normal form $\theta_1 \cdots \theta_{M-1}$, with $\theta_{M-1} \in \Lambda_2$ if $M > 1$, and an infinite element $y \in \Lambda_1^\infty$ such that
\[
x = \theta_1 \cdots \theta_{M-1} y := \bigcup \{ [\theta_1 \cdots \theta_{M-1} \eta] : \eta \in y \}.
\]

\item \label{idhsil three} there is a finite normal form $\theta_1 \cdots \theta_{M-1}$, with $\theta_{M-1} \in \Lambda_1$ if $M > 1$, and an element $y \in \Lambda_2^\infty$ such that
\[
x = \theta_1 \cdots \theta_{M-1} y := \bigcup \{ [\theta_1 \cdots \theta_{M-1} \eta] : \eta \in y \}.
\]

\end{enumerate}

\end{Theorem}

\begin{proof}
Let $x \in \Lambda^*$ be infinite. Let $M = \sup \{m : \text{ there is } \lambda \in x \text{ such that } \lambda = \theta_1 \cdots \theta_m \text{ in normal form}\}$. We first consider the case that $M = \infty$. For $m \ge 1$ let $\lambda_m \in x$ have normal form $\lambda_m = \theta_{m,1} \cdots \theta_{m,m}$. By Lemma \ref{lemma:common extension}, $\theta_{m,i} = \theta_{n,i}$ for all $i < \min\{m,n\}$. Therefore $\theta_i := \theta_{m,i}$ for any $m > i$ is well defined, and $\theta_1 \cdots \theta_m \in x$ for all $m$. It follows from Lemma \ref{lemma:common extension} that $x = \cup_{m=1}^\infty [\theta_1 \cdots \theta_m]$.

Next suppose that $M < \infty$ (still assuming that $x$ is infinite). Let $\theta_1 \cdots \theta_M \in x$ be in normal form. If $\theta_M \in \Lambda_1$ then there is a unique infinite element $y \in \Lambda_1^*$ such that $x = \bigcup \{ [\theta_1 \cdots \theta_{M-1} \eta] : \eta \in y \}$. Similarly, if $\theta_M \in \Lambda_2$ then there is a unique infinite element $y \in \Lambda_2^*$ such that $x = \bigcup \{ [\theta_1 \cdots \theta_{M-1} \eta] : \eta \in y \}$.
\end{proof}

\begin{Notation} \label{notation infinite words}

In case \eqref{idhsil one} of Theorem \ref{theorem:infinite directed hereditary sets in Lambda} we say that $x = \theta_1 \theta_2 \cdots$ is in \textit{normal form}. In cases \eqref{idhsil two} and \eqref{idhsil three} we say that $x = \theta_1 \cdots \theta_{M-1} y$ is in \textit{normal form}.

We may represent infinite directed hereditary subsets of $\Lambda$ as \textit{infinite words} as follows. If $x$ is as in Theorem \ref{theorem:infinite directed hereditary sets in Lambda}\eqref{idhsil one} for the sequence $\theta_1$, $\theta_2$, $\ldots$, write $\theta_i = \mu_{i,1} \ldots \mu_{i,n_i}$ as a sequence of edges in $\Lambda_1$ or $\Lambda_2$. Then we associate $x$ with the infinite word
\[
\mu_{1,1} \ldots \mu_{1,n_1} \,
\mu_{2,1} \ldots \mu_{2,n_2} \,
\mu_{3,1} \ldots \mu_{3,n_3} \,
\cdots
\]
This representation need not be unique, since if $\theta_i \in \Lambda_1$ there may be various expressions for $\theta_i$ as a sequence of edges. However, if the infinite directed hereditary set $x$ is represented by two infinite words, $x_1 x_2 \cdots$ and $x_1' x_2' \cdots$, then $x_i \in \Lambda_1$ if and only if $x_i' \in \Lambda_1$, and if, say, $x_i$, $x_i' \in \Lambda_2$ then $x_i = x_i'$. It is easily seen that the same result holds if $x$ is as in Theorem \ref{theorem:infinite directed hereditary sets in Lambda}\eqref{idhsil two} or \eqref{idhsil three}.

\end{Notation}

\begin{Theorem} \label{thm:maximal elements}

Let $x \in \Lambda^*$. Then $x \in \Lambda^{**}$ if and only if $x$ is infinite, and $x$ has the form either of Theorem \ref{theorem:infinite directed hereditary sets in Lambda}\eqref{idhsil one} or \eqref{idhsil three}, or of Theorem \ref{theorem:infinite directed hereditary sets in Lambda}\eqref{idhsil two} with $y = \alpha^\infty \beta^\infty$.

\end{Theorem}

\begin{proof}
Let $x \in \Lambda^{**}$. It is clear that $x$ must be infinite. Suppose that $x$ has the form of Theorem \ref{theorem:infinite directed hereditary sets in Lambda}\eqref{idhsil two}. Let $x = \theta_1 \cdots \theta_{m-1} y$, $y \in \Lambda_1^\infty$, as in Theorem \ref{theorem:infinite directed hereditary sets in Lambda}\eqref{idhsil two}. Then $x \subseteq \theta_1 \cdots \theta_{m-1} \alpha^\infty \beta^\infty$. Since $x$ is maximal, $x = \theta_1 \cdots \theta_{m-1} \alpha^\infty \beta^\infty$.

Conversely, let $x \in \Lambda^*$ be infinite. Suppose first that $x$ is as in Theorem \ref{theorem:infinite directed hereditary sets in Lambda}\eqref{idhsil one}. Let $\theta_1$, $\theta_2$, $\ldots$ be as in Theorem \ref{theorem:infinite directed hereditary sets in Lambda}\eqref{idhsil one}. We will show that $x$ is maximal. Let $x' \in \Lambda^*$ with $x \subseteq x'$. Let $\mu \in x'$, and write $\mu = \phi_1 \cdots \phi_n$ in normal form. Let $\lambda = \theta_1 \cdots \theta_{n+1} \in x \subseteq x'$. Since $\lambda$, $\mu \in x'$ and $x'$ is directed, $\lambda$ and $\mu$ have a common extension. By Lemma \ref{lemma:common extension} \eqref{lemma:common extension 1} we have that $\mu \in [\lambda] \subseteq x$, since $x$ is hereditary. Thus $x' \subseteq x$. Now suppose that $x$ is as in Theorem \ref{theorem:infinite directed hereditary sets in Lambda}\eqref{idhsil two} with $y = \alpha^\infty \beta^\infty$. Write $x = \theta_1 \cdots \theta_{m-1} \alpha^\infty \beta^\infty$ as in Theorem \ref{theorem:infinite directed hereditary sets in Lambda}\eqref{idhsil two}. Again let $x' \in \Lambda^*$ with $x \subseteq x'$. Let $\mu \in x'$. For each $k \ge 0$, let $\lambda_k = \theta_1 \cdots \theta_{m-1} \alpha^k \beta^k \in x \subseteq x'$. Then $\lambda_k$ and $\mu$ have a common extension. Write $\mu = \phi_1 \cdots \phi_n$ in normal form. If $n < m$ then Lemma \ref{lemma:common extension} \eqref{lemma:common extension 1} implies that $\mu \in [\lambda_k] \subseteq x$. If $n = m$ then $\phi_m \in \Lambda_1$ and hence $\phi_m \in [\alpha^k \beta^k]$ for some $k$. Then again $\mu \in [\lambda_k] \subseteq x$. If $n > m$ then by Lemma \ref{lemma:common extension} \eqref{lemma:common extension 1}, $\phi_n$ extends $\alpha^k \beta^k$ for all $k$, which is impossible. Therefore $x' \subseteq x$. Finally, suppose that $x$ is as in Theorem \ref{theorem:infinite directed hereditary sets in Lambda}\eqref{idhsil three}. Write $x = \theta_1 \cdots \theta_{m-1} y$, where $y = \gamma_\ell^{(i_\ell)} \gamma_{\ell+1}^{(i_{\ell+1})} \cdots \in v_\ell \Lambda_2^\infty$. Let $x' \in \Lambda^*$ with $x \subseteq x'$. Let $\mu \in x'$. Write $\mu = \phi_1 \cdots \phi_n$ in normal form. Choose $p \ge \ell$ such that $| \mu| < |\theta_1 \cdots \theta_{m-1}| + p - \ell + 1$. Then $n \le m$, and by Lemma \ref{lemma:common extension}, $\mu \in [\theta_1 \cdots \theta_{m-1} \gamma_\ell^{(i_\ell)} \cdots \gamma_p^{(i_p)}] \subseteq x$. Therefore $x' \subseteq x$.
\end{proof}

\begin{Proposition}

The closure of $\Lambda^{**}$ in $\Lambda^*$ equals the subset of infinite elements of $\Lambda^*$.

\end{Proposition}

\begin{proof}
Let $x \in \Lambda^*$ be infinite. We may as well assume that $x \not\in \Lambda^{**}$. By Theorem \ref{thm:maximal elements}, $x = \theta_1 \cdots \theta_{m-1} y$, where $y \in \Lambda_1^\infty$ and $y \not= \alpha^\infty \beta^\infty$. By Lemma \ref{lemma:infinite directed hereditary sets in Lambda one} we let, say, $y = \alpha^\ell \beta^\infty$ for some $0 \le \ell < \infty$. Let $Z(\lambda) \setminus \cup_{i=1}^t Z(\mu_i) \in \CB^*$ be a neighborhood of $x$. Then $\lambda \in x$ and for all $i$, $\mu_i \not\in x$. Let $\lambda' = \theta_1 \cdots \theta_{m-1} \alpha^\ell \beta^j$, where $j$ is so large that $|\lambda'|  > |\mu_i|$ for all $1 \le i \le t$. Further, by increasing $j$ if necessary, we choose $j$ so that letting $v_h = s(\lambda')$ gives $k_{h+1}>0$. Then $x \in Z(\lambda') \setminus \cup_{i=1}^t Z(\mu_i) \subseteq Z(\lambda) \setminus \cup_{i=1}^t Z(\mu_i)$. Fix $i$, $1 \le i \le t$. Since $\mu_i \not\in x$ there are two possibilities: either the normal form of $\mu_i$ differs from that of $\lambda'$ in a term before the $m$th term, or $\mu_i$ agrees with $\lambda'$ in the first $m-1$ terms, and has $m$th term equal to $\alpha^{\ell+1}\delta$ for some $\delta \in \Lambda_1$.  Put $\lambda'' = \lambda' \gamma_{h+1}^{(1)}$. Then in all cases, $\mu_i$ and $\lambda''$ have no common extension. Let $x' = \lambda'' \alpha^\infty \beta^\infty \in \Lambda^{**}$. Then $x' \in \Lambda^{**} \cap (Z(\lambda') \setminus \cup_{i=1}^t Z(\mu_i))$. Therefore $x \in \overline{\Lambda^{**}}$.

Now let $x \in \Lambda^*$ be finite, so there is $\lambda \in \Lambda$ with $x = [\lambda]$. Let $n = |\lambda|$. Then $Z(\lambda) \setminus (Z(\lambda \alpha_{n+1}) \cup Z(\lambda \beta_{n+1}) \cup \cup_{j=1}^{k_{n+1}} Z(\lambda \gamma_{n+1}^{(j)})) = \{x\}$. Therefore $x \not\in \overline{\Lambda^{**}}$.
\end{proof}

 We will simplify the situation further by restricting to a transversal in $\partial \Lambda$ (in the sense of \cite[Example 2.7]{muhrenwil}).

\begin{Lemma} \label{lemma transversal}

Let $X = v_1 \partial \Lambda$. Then $X$ is a compact open transversal in $G(\Lambda)|_{\partial \Lambda}$.

\end{Lemma}

\begin{proof}
$X$ is the intersection of $\partial \Lambda $ with the compact open subset $v_1 \Lambda^*$, hence $X$ is compact open. To see that $X$ is a transversal, let $x \in \Lambda^*$, and let $\mu \in v_1 \Lambda r(x)$ be arbitrary. Then $\mu x \in X$. Since $G(\Lambda)|_{\partial \Lambda}$ is \'etale, $r$ and $s$ are open maps. Since $\partial\Lambda \cdot G(\Lambda) \cdot X$ is open in $G(\Lambda)|_{\partial \Lambda}$, the restrictions of $r$ and $s$ to it are open as well. This verifies the hypotheses of \cite[Example 2.7]{muhrenwil}.
\end{proof}

It now follows from \cite[Theorem 2.8]{muhrenwil} that $C^*(\Lambda)$ is Morita equivalent to $C^*(G(\Lambda)|_X)$.

\begin{Definition} \label{def the groupoid G}

We will write $G := G(\Lambda)|_X$, the restriction of $G(\Lambda)$ to the transversal $X$.

\end{Definition}

For the rest of this paper (except for section \ref{section stability}) we will study $G$ and $C^*(G)$. Thus $G = \{ [\mu,\nu,x] \in G(\Lambda) : r(\mu) = r(\nu) = v_1,\ x \in s(\mu) \partial \Lambda \}$. Note that for $[\mu,\nu,x] \in G$ it is necessarily the case that $|\mu| = |\nu|$.

\begin{Definition} \label{def subgroupoids}

For $i \ge 1$ let $G_i$ be the subgroupoid of $G$ generated by elements of the form $[\mu,\nu,x] \in G$ such that $|\mu| = |\nu| \le i$.

\end{Definition}

\begin{Definition} \label{def X_p}

For $p \ge 1$, let $X_p = v_p \partial \Lambda$ (thus $X_1 = X$).

\end{Definition}

The following simple lemma will be useful several times.

\begin{Lemma} \label{simple lemma}

For $p \ge 1$ and $x \in X_{p+1}$,
\[
[\beta^p,\alpha^p,x] = \prod_{j=0}^{p-1} [\beta,\alpha,\alpha^j \beta^{p - j - 1} x].
\]

\end{Lemma}

\begin{proof}
The proof is by induction on $p$. This is clearly true when $p = 1$. Let $p \ge 1$ and suppose the formula is true for $p$. Then
\begin{align*}
\prod_{j=0}^p [\beta,\alpha,\alpha^j \beta^{p-j} x]
&= \prod_{j=0}^{p-1} [\beta,\alpha,\alpha^j \beta^{p - j - 1} \beta x] [\beta,\alpha,\alpha^p x] \\
&= [\beta^p,\alpha^p,\beta x][\beta,\alpha,\alpha^p x], \text{ by the inductive hypothesis,} \\
&= [\beta^{p+1},\alpha^p\beta,x][\alpha^p\beta,\alpha^{p+1},x] \\
&= [\beta^{p+1},\alpha^{p+1},x],
\end{align*}
which is the formula for the case $p+1$.
\end{proof}

\begin{Theorem} \label{thm elements of G_i}

$G_i = \{ [\mu\theta,\mu'\theta',x] : |\mu| = |\mu'| \le i \text{ and } \theta, \theta' \in \Lambda_1 \}$.

\end{Theorem}

Before beginning the proof we note the following. Let $[\mu\theta,\mu'\theta',x]$ be an element of the righthand side in Theorem \ref{thm elements of G_i}. We may write $\theta = \phi \eta$, $\theta' = \phi' \eta$ where $d(\phi) \perp d(\phi')$ in $\IN^2$. Then $[\mu\theta,\mu'\theta',x] = [\mu\phi,\mu'\phi',\eta x]$. Thus when describing the elements of $G_i$ we may assume that $d(\theta) \perp d(\theta')$. (This means that one of $\theta$, $\theta'$ equals $\alpha^\ell$ and the other equals $\beta^\ell$.)

\begin{proof}
We begin with the proof of the containment ``$\supseteq$''. Let $|\mu| = |\mu'| = m \le i$ and let $\ell \ge 0$. Then for $x \in X_{|\mu| + \ell + 1}$,
\begin{align*}
[\mu \beta^\ell, \mu' \alpha^\ell, x]
&= [\mu,\beta^m,\beta^\ell x] [\beta^{m+\ell},\alpha^{m+\ell},x] [\alpha^m,\mu',\alpha^\ell x] \\
&= [\mu,\beta^m,\beta^\ell x] \bigl( \prod_{j=0}^{m+\ell-1}[\beta,\alpha,\alpha^j\beta^{m + \ell - j - 1} x] \bigr) [\alpha^m,\mu',\alpha^\ell x]
\end{align*}
(by Lemma \ref{simple lemma}), which is a product of generators of $G_i$. The case where the roles of $\alpha$ and $\beta$ are switched is obtained by taking inverses.

For the reverse containment, let us denote by $S$ the set on the righthand side of the statement. It is clear that all generators of $G_i$ belong to $S$. The proof will be finished if we show that the product of an element of $S$ and a generator is again an element of $S$. Note that if $[\mu,\mu',x]$ is a generator, and $m = |\mu| = |\mu'| < i$, then we may write $x = \epsilon y$ where $|\epsilon| = i - m$. Then $[\mu,\mu',x] = [\mu\epsilon,\mu'\epsilon,y]$, and $|\mu \epsilon| = |\mu' \epsilon| = i$. Thus we may assume that the first two coordinates of a generator of $G_i$ have length $i$. Now let $|\mu| = |\mu'| = i$, let $\ell \ge 0$, and let $|\nu| = |\nu'| = i$. Then assuming that $x$, $y$ are such that $\mu' \beta^\ell x = \nu y$, the product $[\mu \alpha^\ell, \mu' \beta^\ell,x] [\nu,\nu',y]$ is defined. Then there are $\xi$, $\eta$, and $z$ such that $x = \xi z$, $y = \eta z$, and $\mu' \beta^\ell \xi = \nu \eta$, and then
\begin{align*}
[\mu \alpha^\ell, \mu' \beta^\ell,x] [\nu,\nu',y]
&= [\mu \alpha^\ell \xi, \mu' \beta^\ell \xi, z] [\nu \eta, \nu' \eta, z] \\
&= [\mu \alpha^\ell \xi, \nu' \eta, z].
\end{align*}
Write $\mu' = \theta_1 \cdots \theta_p$ and $\nu = \phi_1 \cdots \phi_q$ in normal form. We claim that $p = q$. For suppose, e.g., that $p < q$. We know that $\mu'\beta^\ell\xi = \nu \eta$. Then by Lemma \ref{lemma:common extension} we have that $\mu^\prime = \phi_1 \cdots \phi_{p-1} \theta_p$ and that $\phi_p$ extends $\theta_p$. Since $|\mu'| = |\nu|$ we have that $|\theta_p| = |\phi_p \cdots \phi_q|$. But $|\theta_p| \le |\phi_p| < |\phi_p \cdots \phi_q|$, a contradiction. Thus $p = q$. Then $\theta_p$ and $\phi_p$ have a common extension. If $\theta_p$, $\phi_p \in \Lambda_2$ then $\theta_p = \phi_p$, hence $\mu' = \nu$. In this case it follows that  $\eta = \beta^\ell \xi$, and hence that $[\mu \alpha^\ell \xi,\nu' \eta,z] = [\mu \alpha^\ell,\nu' \beta^\ell,\xi z] \in S$. On the other hand, if $\theta_p$, $\phi_p \in \Lambda_1$, we know only that $|\theta_p| = |\phi_p|$ and $\theta_p \beta^\ell \xi = \phi_p \eta$. Write $\xi = \xi' \xi''$, $\eta = \eta' \eta''$, where $\xi'$, $\eta' \in \Lambda_1$, and $\xi''$, $\eta''$ have normal forms with first factor in $\Lambda_2$ (if nontrivial). Then $\theta_p \beta^\ell \xi'$ and $\phi_p \eta'$ are the first factors in the normal forms of $\theta_p \beta^\ell \xi$ and of $\phi_p \eta$. By uniqueness of normal form, these first factors are equal, and hence also $\xi'' = \eta''$. Thus again we find that $[\mu \alpha^\ell \xi,\nu' \eta,z] = [\mu \alpha^\ell \xi',\nu' \eta', \xi'' z] \in S$. A similar argument treats the case where an element of $S$ is multiplied on the left by a generator of $G_i$.
\end{proof}

\section{$K$-theory of $C^*(G_i)$}
\label{section K-theory of G_i}

Let $i$ be fixed throughout this section. We now begin the analysis of $C^*(G_i)$. First we give a decomposition of $X = G_i^{(0)}$.

\begin{Definition} \label{def subsets of X}

Let
\begin{align*}
U_i &= \{x = x_1 x_2 \cdots \in X : x_\ell \in \Lambda_2 \text{ for some } \ell > i \}, \\
F_i &= X \setminus U_i, \\
\noalign{ and for $\ell > i$, let}
E_\ell &= \{x \in X : x_\ell \in \Lambda_2, \text{ and } x_j \in \Lambda_1 \text{ for } i < j < \ell \}, \\
\Omega_\ell &= \{ \lambda \in v_1 \Lambda v_\ell : \lambda_j \in \Lambda_1 \text{ for } i < j < \ell \}. \\
\noalign{(Note that $E_\ell = \varnothing$ if $k_\ell = 0$.)}
\end{align*}

\end{Definition}

\begin{Lemma} \label{lem equivalent projections}

Let $\ell > i$ and let $\mu$, $\nu \in \Omega_\ell$. Then $\chi_{Z(\mu)}$ and $\chi_{Z(\nu)}$ are equivalent projections, and $\chi_{Z(\mu) \setminus Z(\mu\beta)}$ and $\chi_{Z(\nu) \setminus Z(\nu\beta)}$ are equivalent projections, in $C^*(G_i)$.

\end{Lemma}

\begin{proof}
Write $\mu = \mu_1 \mu_2$ and $\nu = \nu_1 \nu_2$ where $|\mu_1| = |\nu_1| = i$ and $\mu_2$, $\nu_2 \in \Lambda_1$. Without loss of generality we may suppose that $\mu_2 = \theta \alpha^p$ and $\nu_2 = \theta \beta^p$, with $\theta \in \Lambda_1$. Now
\begin{align*}
[\nu,\mu,Z(v_\ell)]
&= [\nu_1 \theta \beta^p, \alpha^i \theta \beta^p, Z(v_\ell)] \cdot [\beta^p, \alpha^p, Z(\alpha^i \theta)] \cdot [\alpha^i \theta \alpha^p, \mu_1 \theta \alpha^p, Z(v_\ell)] \\
&= [\nu_1 \theta \beta^p, \alpha^i \theta \beta^p, Z(v_\ell)] \cdot \prod_{j=0}^{p-1} [\beta, \alpha, Z(\beta^{p-j-1} \alpha^i \theta)] \cdot [\alpha^i \theta \alpha^p, \mu_1 \theta \alpha^p, Z(v_\ell)],
\end{align*}
by Lemma \ref{simple lemma}. Therefore $u = \chi_{[\nu,\mu,Z(v_\ell)]} \in C^*(G_i)$, and $u^*u = \chi_{[v_1,v_1,Z(\mu)]}$, $uu^* = \chi_{[v_1,v_1,Z(\nu)]}$. Now we may let $v = \chi_{[\nu,\mu,Z(v_\ell) \setminus Z(v_\ell\beta)]}$, and we have that $v^*v = \chi_{Z(\mu) \setminus Z(\mu\beta)}$ and $vv^* = \chi_{Z(\nu) \setminus Z(\nu\beta)}$.
\end{proof}

\begin{Lemma} \label{lemma structure of E ell}

For each $\ell > i$, $E_\ell$ is an open $G_i$-invariant subset of $X$. Each point of $E_\ell$ has trivial isotropy, and has finite orbit equivalent to $\Omega_\ell$.

\end{Lemma}

\begin{proof}
$E_\ell = \sqcup_{\lambda \in \Omega_\ell} \sqcup_{r = 1}^{k_\ell} Z(\lambda \gamma_\ell^{(r)})$. Therefore $E_\ell$ is open. Now we show that $E_\ell$ is $G_i$-invariant. Let $x \in E_\ell$. Then $x = \lambda \phi \gamma_\ell^{(r)} x'$, where $|\lambda| \le i$ and $\phi \in \Lambda_1$. Let $g \in G_i x$. By Theorem \ref{thm elements of G_i}, $g = [\mu\theta,\mu'\theta',y]$, where $|\mu| = |\mu'| \le i$, $\theta$, $\theta' \in \Lambda_1$, and $\mu'\theta' y = x$.  Then $\mu'\theta' y = \lambda \phi \gamma_\ell^{(r)} x'$. It follows that $|\mu'\theta'| < \ell$, and that $y = \theta'' \gamma_\ell^{(r)} x'$, where $\theta'' \in \Lambda_1$. Then $r(g) = \mu\theta y = \mu \theta \theta'' \gamma_\ell^{(r)} x' \in E_\ell$.

With $x$ and $g$ as above we see that $r(g) \in \Omega_\ell \gamma_\ell^{(r)} x'$, so the orbit of $x$ under $G_i$ is $\Omega_\ell \gamma_\ell^{(r)} x'$, which is equivalent to $\Omega_\ell$. If we assume that $r(g) = x$ (so that $g$ is in the isotropy subgroup of $G_i$ at $x$), then $\mu \theta \theta'' \gamma_\ell^{(r)} x' = \mu \theta y = \mu' \theta' y = \mu' \theta' \theta'' \gamma_\ell^{(r)} x'$, and hence $\mu \theta = \mu' \theta'$. Therefore $g = [v_1,v_1, x] \in G^{(0)}$.
\end{proof}

\begin{Proposition}

$U_i$ is an open $G_i$-invariant subset of $X$, and
\[
C^*(G_i|_{U_i}) \cong  \bigoplus_{\ell > i} \bigl(M_{\Omega_\ell \times \Omega_\ell} \otimes C(X_{\ell+1} \bigr)^{k_\ell}.
\]
(Notice that if $k_\ell = 0$ then the $\ell$th summand is not present.)

\end{Proposition}

\begin{proof}
Note that $U_i = \sqcup_{\ell > i} E_\ell$, and hence $U_i$ is open and $G_i$-invariant. Moreover, $C^*(G_i|_{U_i}) = \bigoplus_{\ell > i} C^*(G_i|_{E_\ell})$. But $E_\ell = \sqcup_{r=1}^{k_\ell} \Omega_\ell \gamma_\ell^{(r)} X_{\ell+1}$, so $G_i|_{E_\ell} = \sqcup_{r=1}^{k_\ell} (\Omega_\ell \times \Omega_\ell) \times \gamma_\ell^{(r)} X_{\ell+1}$. It follows that $C^*(G_i|_{E_\ell}) \cong \oplus_{r=1}^{k_\ell} M_{\Omega_\ell \times \Omega_\ell} \otimes C(X_{\ell+1}) = \bigl( M_{\Omega_\ell \times \Omega_\ell} \otimes C(X_{\ell+1}) \bigr)^{k_\ell}$. The isomorphism in the proposition now follows.
\end{proof}

\begin{Remark} \label{remark U i amenable}
It follows from Lemma \ref{lemma structure of E ell} that $G_i|_{U_i}$ is an AF groupoid (see \cite[Definition III.1.1]{ren}), and hence is amenable.
\end{Remark}

\begin{Remark} \label{rem isomorphism in E ell}

The isomorphism $C^*(G_i|_{E_\ell}) \cong \bigl( M_{\Omega_\ell \times \Omega_\ell} \otimes C(X_{\ell+1}) \bigr)^{k_\ell}$ is given by
\[
\chi_{[\mu,\nu,\gamma_\ell^{(r)} F]} \mapsto (0, \ldots, e_{\mu,\nu} \otimes \chi_F, \ldots,0),
\]
(in the $r$th coordinate) for $F \subseteq X_{\ell+1}$ a compact open subset.

\end{Remark}

\begin{Notation} \label{notation C*(Gi,Ui)}

For $f \in C^*(G_i|_{U_i})$ we write $f = (f_{i+1},f_{i+2},\ldots)$ with $f_\ell = (f_{\ell,1}, \ldots, f_{\ell,k_\ell})$, where $f_{\ell,r} \in M_{\Omega_\ell \times \Omega_\ell} \otimes C(X_{\ell+1}) \cong C(X_{\ell+1},M_{\Omega_\ell \times \Omega_\ell})$.

\end{Notation}

\begin{Proposition} \label{prop K-theory of G i U i}

$K_1(C^*(G_i|_{U_i})) = 0$, and $K_0(C^*(G_i|_{U_i})) \cong \bigoplus_{\ell > i} C(X_{\ell+1},\IZ)^{k_\ell}$, with generators $\bigl[ \chi_{[\alpha^{\ell-1},\alpha^{\ell-1},\gamma_\ell^{(r)} F]} \bigr]_0$ for $F \subseteq X_{\ell+1}$ a compact open subset.

\end{Proposition}

\begin{proof}
Since $X_{\ell+1}$ is totally disconnected, $C(X_{\ell+1})$ is AF. It follows that $C^*(G_i|_{U_i})$ is AF, and hence that $K_1(C^*(G_i|_{U_i})) = 0$. Moreover, $K_0(M_{\Omega_\ell \times \Omega_\ell} \otimes C(X_{\ell+1})) \cong C(X_{\ell+1},\IZ)$. The description of $K_0(C^*(G_i|_{U_i}))$ follows. The description of generators of the group follows from Remark \ref{rem isomorphism in E ell}.
\end{proof}

The open invariant set $U_i$ determines an ideal in $C^*(G_i)$, and its complement $F_i$ is a closed invariant set determining the quotient $C^*$-algebra. There is thus a short exact sequence
\begin{equation} \label{equation first exact sequence}
0 \longrightarrow C^*(G_i|_{U_i}) \longrightarrow C^*(G_i) \longrightarrow C^*(G_i|_{F_i}) \longrightarrow 0.
\end{equation}
We now study $C^*(G_i|_{F_i})$. Note that $F_i = \{ \lambda \alpha^m \beta^n : |\lambda| \le i \text{ and } m + n = \infty \}$. Let
\begin{align*}
F_i^\infty &= \{ \lambda \alpha^\infty \beta^\infty : |\lambda| \le i \} \\
F_i^0 &= F_i \setminus F_i^\infty.
\end{align*}
Then $F_i^\infty$ is a finite invariant set, and hence is also closed. Therefore $F_i^0$ is a relatively open invariant subset of $F_i$. We have a further exact sequence
\begin{equation} \label{equation second exact sequence}
0 \longrightarrow C^*(G_i|_{F_i^0}) \overset{\iota}{\longrightarrow} C^*(G_i|_{F_i}) \overset{\pi}{\longrightarrow} C^*(G_i|_{F_i^\infty}) \longrightarrow 0.
\end{equation}
We first analyze $G_i|_{F_i^\infty}$. The following definition will be convenient in several places.

\begin{Definition} \label{def Phi_i}

Let $\Phi_i = \{\lambda \in v_1\Lambda : |\lambda| \le i,\ \lambda_{|\lambda|} \in \Lambda_2\}$.

\end{Definition}

\begin{Lemma} \label{lem F i infty}

$F_i^\infty$ is a single orbit for $G_i$. The isotropy at points of $F_i^\infty$ is infinite cyclic. Then the map $\lambda \mapsto \lambda \alpha^\infty \beta^\infty$ is a bijection of $\Phi_i$ onto $F_i^\infty$.

\end{Lemma}

\begin{proof}
Let $x$, $y \in F_i^\infty$. Then $x = \lambda \alpha^\infty \beta^\infty$, $y = \mu \alpha^\infty \beta^\infty$, where as in the proof of Theorem \ref{thm elements of G_i}, we may assume that $|\lambda| = |\mu| = i$. Then $g = [\lambda,\mu, \alpha^\infty \beta^\infty] \in G_i$, and $s(g) = \mu \alpha^\infty \beta^\infty = y$, $r(g) = \lambda \alpha^\infty \beta^\infty = x$. This proves that $G_i|_{F_i^\infty}$ has a single orbit. Now let $g$ be in the isotropy group at, say, $v_1 \alpha^\infty \beta^\infty $. Then $g = [\lambda,\mu,\alpha^\infty \beta^\infty]$, with $|\lambda| = |\mu| \le i$, and $\mu \alpha^\infty \beta^\infty = v_1 \alpha^\infty \beta^\infty = \lambda \alpha^\infty \beta^\infty $. Then $\lambda$, $\mu \in \Lambda_1$. We may assume without loss of generality that $\lambda = \alpha^p$ and $\mu = \beta^p$ for some $p \ge 0$. Then by Lemma \ref{simple lemma} we have that $g = [\alpha^p,\beta^p,\alpha^\infty \beta^\infty] = [\alpha,\beta,\alpha^\infty \beta^\infty]^p$. The last statement of the lemma is clear.
\end{proof}

\begin{Remark} \label{remark F i infty amenable}
It is easy to verify that the range and source maps on $G_i|_{F_i^\infty}$ are open maps. It then follows from Lemma \ref{lem F i infty} and \cite[Corollary 9.78]{will} that $G_i|_{F_i^\infty}$ is amenable.
\end{Remark}

\begin{Corollary} \label{cor K-theory of F i infty}

$C^*(G_i|_{F_i^\infty}) \cong  M_{\Phi_i \times \Phi_i} \otimes C(\IT)$. The isomorphism is given by
\begin{align*}
[\alpha,\beta,\alpha^\infty \beta^\infty ] &\mapsto e_{v_1,v_1} \otimes z \\
[\lambda,\alpha^{|\lambda|},\alpha^\infty \beta^\infty] &\mapsto e_{\lambda,v_1} \otimes 1.
\end{align*}

\end{Corollary}

\begin{proof}
The isomorphism follows from, e.g., \cite[Theorem 3.1]{muhrenwil}. The explicit version is clear from the proof of Lemma \ref{lem F i infty}.
\end{proof}

Now we study $F_i^0 = \{\lambda \alpha^m \beta^n : |\lambda| \le i, \max\{m,n\} = \infty, \min\{m,n\} < \infty \}$.

\begin{Lemma} \label{lem F i 0}

$G_i|_{F_i^0}$ is principal with two orbits:
\begin{align*}
F_i^{0,1} &= \{\lambda \alpha^\infty \beta^n : |\lambda| \le i, n \ge 0 \} \\
F_i^{0,2} &= \{\lambda \alpha^m \beta^\infty : |\lambda| \le i, m \ge 0 \}.
\end{align*}

\end{Lemma}

\begin{proof}
We first show that $G_i|_{F_i^0}$ is principal. Let $g = [\lambda,\mu,y] \in G_i|_{F_i^0}$ with $\lambda y = r(g) = s(g) = \mu y$. Write $\lambda = \lambda'\lambda''$ and $\mu = \mu'\mu''$, where $\lambda'$, $\mu'$ have source-most edge in $\Lambda_2$ and $\lambda''$, $\mu'' \in \Lambda_1$. We may assume that, say, $\lambda'' = \alpha^p$, $\mu'' = \beta^q$, and $y = \alpha^\infty \beta^n$. Then our assumption implies that $\lambda' \alpha^\infty \beta^n = \mu' \alpha^\infty \beta^{q + n}$. It follows that $q = 0$ and $\lambda' = \mu'$. Then since $|\lambda| = |\mu|$ we have that $\lambda' \alpha^p = \mu'$. Since the source-most edge of $\mu'$ is in $\Lambda_2$, $p=0$. Therefore $\lambda' = \mu'$, and hence $\lambda = \mu$. Therefore $g \in G^{(0)}$.

Now consider $\lambda \in v_1\Lambda$ with $|\lambda| \le i$, and $n \ge 0$. By Lemma \ref{simple lemma} we have
\begin{align*}
[\lambda \beta^n,\alpha^{|\lambda| + n}, \alpha^\infty]
&= [\lambda,\alpha^{|\lambda|},\alpha^\infty\beta^n] [\beta^n,\alpha^n,\alpha^\infty] \\
&= [\lambda,\alpha^{|\lambda|},\alpha^\infty\beta^n] \prod_{j=0}^{n-1} [\beta,\alpha,\alpha^\infty \beta^{n-j-1}] \\
&\in G_i|_{F_i^0},
\end{align*}
has source $\alpha^\infty$, and has range $\lambda \alpha^\infty \beta^n$, a typical element of $F_i^{0,1}$. Thus $F_i^{0,1}$ is contained in the orbit of $\alpha^\infty$. A similar argument shows that $F_i^{0,2}$ is contained in the orbit of $\beta^\infty$. To finish the proof we observe that $\alpha^\infty$ and $\beta^\infty$ are not in the same orbit.
\end{proof}

\begin{Remark} \label{remark F i 0 amenable}

Note that the sets $F_i^{0,1}$ and $F_i^{0,2}$ of Lemma \ref{lem F i 0} are disjoint, open and invariant. Therefore $G_i|_{F_i^0} = G_i|_{F_i^{0,1}} \sqcup G_i|_{F_i^{0,2}}$. Again, it is easy to see that the range and source maps in $G_i|_{F_i^0}$ are open. Then \cite[Corollary 9.78]{will} implies that $G_i|_{F_i^0}$ is amenable.

\end{Remark}

\begin{Proposition} \label{prop G i amenable}

$G_i$ is amenable.

\end{Proposition}

\begin{proof}
From Remarks \ref{remark F i infty amenable} and \ref{remark F i 0 amenable}, and \cite[Proposition 9.83]{will}, it follows that $G_i|_{F_i}$ is amenable. Then it follows from Remark \ref{remark U i amenable} and the same Proposition that $G_i$ is amenable.
\end{proof}

\begin{Corollary} \label{cor K-theory of F i 0}

$C^*(G_i|_{F_i^0}) \cong M_{\Phi_i \times \Phi_i} \otimes (\CK \oplus \CK)$. The isomorphism is given by
\begin{align*}
[\beta^n,\alpha^n,\alpha^\infty] &\mapsto e_{v_1,v_1} \otimes (e_{n,0} \oplus 0) \\
[\lambda,\alpha^{|\lambda|},\alpha^\infty] &\mapsto e_{\lambda,v_1} \otimes (e_{0,0} \oplus 0) \\
[\alpha^n,\beta^n,\beta^\infty] &\mapsto e_{v_1,v_1} \otimes (0 \oplus e_{n,0}) \\
[\lambda,\beta^{|\lambda|},\beta^\infty] &\mapsto e_{\lambda,v_1} \otimes (0 \oplus e_{0,0}).
\end{align*}

\end{Corollary}

\begin{proof}
The isomorphism follows from, e.g., \cite[Theorem 3.1]{muhrenwil}. The explicit version is clear from the proof of Lemma \ref{lem F i 0}.
\end{proof}

The exact sequence \eqref{equation second exact sequence} becomes
\begin{equation*}
0 \to M_{\Phi_i \times \Phi_i} \otimes (\CK \oplus \CK) \overset{\iota}{\longrightarrow} C^*(G_i|_{F_i}) \overset{\pi}{\longrightarrow} M_{\Phi_i \times \Phi_i} \otimes C(\IT) \to 0.
\end{equation*}

The corresponding long exact sequence in $K$-theory is
\[
\begin{tikzpicture}[xscale=3,yscale=1.5]

\node (00) at (0,0) [rectangle] {$\IZ$};
\node (10) at (1,0) [rectangle] {$K_0(C^*(G_i|_{F_i}))$};
\node (20) at (2,0) [rectangle] {$\IZ^2$};
\node (01) at (0,1) [rectangle] {$0$};
\node (11) at (1,1) [rectangle] {$K_1(C^*(G_i|_{F_i}))$};
\node (21) at (2,1) [rectangle] {$\IZ$};

\draw[-latex,thick] (00) -- (01);
\draw[-latex,thick] (01) -- (11) node[pos=0.5,inner sep=0.5pt, above=1pt] {$\iota_{*,1}$};
\draw[-latex,thick] (11) -- (21) node[pos=0.5,inner sep=0.5pt, above=1pt] {$\pi_{*,1}$};
\draw[-latex,thick] (21) -- (20) node[pos=0.5,inner sep=0.5pt, right=1pt] {$\partial_1$};
\draw[-latex,thick] (20) -- (10) node[pos=0.5,inner sep=0.5pt, below=1pt] {$\iota_{*,0}$};
\draw[-latex,thick] (10) -- (00) node[pos=0.5,inner sep=0.5pt, below=1pt] {$\pi_{*,0}$};

\end{tikzpicture}
\]
In Corollaries \ref{cor K-theory of F i infty} and \ref{cor K-theory of F i 0} we have identified elements that serve as generators of the $K$-groups at the corners of this exact sequence:
\begin{align*}
\bigl[ \chi_{[\alpha,\beta,\alpha^\infty \beta^\infty]} + \chi_{[v_1,v_1,F_i^\infty \setminus \{\alpha^\infty \beta^\infty \}]} \bigr]_1 &\in K_1(C^*(G_i|_{F_i^\infty})) \\
\bigl[ \chi_{[v_1,v_1,\alpha^\infty \beta^\infty ]} \bigr]_0 &\in K_0(C^*(G_i|_{F_i^\infty})) \\
\bigl[ \chi_{[v_1,v_1,\alpha^\infty]} \bigr]_0, \bigl[ \chi_{[v_1,v_1,\beta^\infty]} \bigr]_0 &\in K_0(C^*(G_i|_{F_i^0})).
\end{align*}
We will compute the index map $\partial_1$.

\begin{Lemma}

The map $\partial_1$ in the above long exact sequence is given by $\partial_1(1) = (1,-1)$.

\end{Lemma}

\begin{proof}
The unitary
\[
u = \chi_{[\alpha,\beta,\alpha^\infty \beta^\infty ]} + \chi_{[v_1,v_1,F_i^\infty \setminus \{\alpha^\infty \beta^\infty \}]} \in C^*(G_i|_{F_i^\infty})
\]
lifts to the following partial isometry in $C^*(G_i|_{F_i})$:
\[
w = \chi_{[\alpha,\beta, \Lambda_1^\infty \cap X_2]} + \chi_{[v_1,v_1,F_i \setminus v_1 \Lambda_1^\infty]}
\in C^*(G_i|_{F_i}).
\]
Then
\begin{align*}
w^*w &= \chi_{[\beta,\beta, \Lambda_1^\infty \cap X_2]} + \chi_{[v_1,v_1, F_i \setminus v_1 \Lambda_1^\infty]} \\
&= \chi_{[v_1,v_1,\beta (\Lambda_1^\infty \cap X_2)]} + \chi_{[v_1,v_1, F_i \setminus v_1 \Lambda_1^\infty]} \\
&= \chi_{[v_1,v_1, \Lambda_1^\infty \setminus \{\alpha^\infty\}]} + \chi_{[v_1,v_1,F_i \setminus v_1 \Lambda_1^\infty]} \\
&= \chi_{[v_1,v_1,F_i \setminus \{\alpha^\infty\}]}, \\
\noalign{and similarly,}
ww^* &= \chi_{[v_1,v_1, F_i \setminus \{\beta^\infty\}]}.
\end{align*}
Therefore
\begin{align*}
\partial_1([u]_1) &= [1 - w^*w]_0 - [1 - ww^*]_0 \\
&= \bigl[ \chi_{[v_1,v_1,\alpha^\infty]} \bigr]_0 - \bigl[ \chi_{[v_1,v_1,\beta^\infty]} \bigr]_0 \\
&\equiv (1,-1) \in \IZ^2. \qedhere
\end{align*} 
\end{proof}

\begin{Proposition} \label{prop K-theory of G i F i}

$K_1(C^*(G_i|_{F_i})) = 0$ and $K_0(C^*(G_i|_{F_i})) \cong \IZ^2$, with generators $\bigl[ \chi_{[v_1,v_1,\alpha^\infty]} \bigr]_0$ and $\bigl[ \chi_{[v_1,v_1,Z(\beta^{i+1}) \cap F_i]} \bigr]_0$.

\end{Proposition}

\begin{proof}
Since $\partial_1$ is injective we know that $\pi_{*,1} = 0$, and hence that $K_1(C^*(G_i|_{F_i})) = 0$. We may choose $(1,0)$ and $(1,-1)$ as basis for $\IZ^2$. Then since $\partial_1(\IZ) = \IZ (1,-1)$ we obtain a short exact sequence
\[
0 \to \IZ(1,0) \to K_0(C^*(G_i|_{F_i})) \to \IZ \to 0.
\]
Therefore $K_0(C^*(G_i|_{F_i})) \cong \IZ^2$. As generators we use $\iota_{*,0}(1,0) = \bigl[ \chi_{[v_1,v_1,\alpha^\infty]} \bigr]_0$, and a lift of $\bigl[ \chi_{[v_1,v_1,\alpha^\infty \beta^\infty]} \bigr]_0$. In choosing the lift, we may choose any compact open subset of $F_i$ whose intersection with $F_i^\infty$ equals $\{\alpha^\infty \beta^\infty\}$. The set $Z(\beta^{i+1}) \cap F_i$ satisfies this requirement, so our second generator can be chosen to be $\bigl[ \chi_{[v_1,v_1,Z(\beta^{i+1}) \cap F_i]} \bigr]_0$.
\end{proof}

Now we describe the $K$-theory of $C^*(G_i)$.

\begin{Theorem} \label{K-theory of G i}

$K_1(C^*(G_i)) = 0$, and $K_0(C^*(G_i)) \cong \bigoplus_{\ell > i} C(X_{\ell+1},\IZ)^{k_\ell} \oplus \IZ^2$. The generators of the first term are as in Proposition \ref{prop K-theory of G i U i}, and the generators of the second term may be chosen to be $\bigl[ \chi_{[v_1,v_1,Z(\alpha^i) \setminus Z(\alpha^i\beta)]} \bigr]_0$ and $\bigl[ \chi_{[v_1,v_1,Z(\beta^{i+1})]} \bigr]_0$.

\end{Theorem}

\begin{proof}
The long exact sequence in $K$-theory associated to the sequence \eqref{equation first exact sequence} reduces to
\[
0 \to K_0(C^*(G_i|_{U_i})) \to K_0(C^*(G_i)) \to K_0(C^*(G_i|_{F_i})) \to 0.
\]
Since $K_0(C^*(G_i|_{F_i})) \cong \IZ^2$ is free abelian, the central group is isomorphic to the direct sum of the other two groups. For generators we may use generators of the subgroup, together with lifts of the generators of the quotient group. Thus we must choose lifts for the generators given in Proposition \ref{prop K-theory of G i F i}. We note that $(Z(v_1 \alpha^i) \setminus Z(v_1 \alpha^i\beta)) \cap F_i = \{ v_1 \alpha^\infty \}$. Therefore $\bigl[ \chi_{[v_1,v_1,Z(\alpha^i) \setminus Z(\alpha^i\beta)]} \bigr]_0$ is a lift of $\bigl[ \chi_{[v_1,v_1,\alpha^\infty]} \bigr]_0$. It is clear that $\bigl[ \chi_{[v_1,v_1,Z(\beta^{i+1})]} \bigr]_0$ is a lift of $\bigl[ \chi_{[v_1,v_1,Z(\beta^{i+1}) \cap F_i]} \bigr]_0$.
\end{proof}

\begin{Corollary} \label{cor typical element of K-theory}

A typical element of $K_0(C^*(G_i))$ has the form
\begin{equation} \label{eqn typical positive element}
a = 
\sum_{\ell > i} \sum_{r = 1}^{k_\ell} \sum_{j = 1}^{h_\ell} c_{\ell,r,j} [\chi_{\alpha^{\ell-1} \gamma_\ell^{(r)} F_{\ell,j}}]_0 + m [\chi_{Z(\alpha^i) \setminus Z(\alpha^i\beta)}]_0 + n [\chi_{Z(\beta^{i+1})}]_0,
\end{equation}
where for each $\ell$, $\{F_{\ell,j} : 1 \le j \le h_\ell \}$ is a partition of $X_{\ell+1}$ by compact open subsets (and the coefficients $c_{\ell,r,j}$, $m$, $n$ are integers).

\end{Corollary}

Now we characterize the positive elements of $K_0(C^*(G_i))$.

\begin{Theorem} \label{thm positive elements}

Let $a \in K_0(C^*(G_i))$ be as in equation \eqref{eqn typical positive element}. Then $a \ge 0$ if and only if for all $\ell > i$, $r$ and $j$, $c_{\ell,r,j} + m + n(\ell - i - 1) \ge 0$.

\end{Theorem}

%\begin{Remark} \label{rem index shift}

%The theorem has an equivalent statement, where the inequality is replaced by $c_{\ell,r,j} + m + n(\ell - i)$, which may appear to be a more natural version (as will be apparent during the proof). However the statement we chose will turn out to be more useful later on.

%\end{Remark}

We give some lemmas before the proof.

\begin{Lemma} \label{lem positivity one}

For $p \ge \ell \ge 1$,

\begin{enumerate}[(i)]

\item \label{lem positivity one (i)} $Z(\alpha^ \ell) \setminus Z(\alpha^\ell\beta) = \bigl( \sqcup_{q = \ell +1}^p \sqcup_{r=1}^{k_q} Z(\alpha^{q-1} \gamma_q^{(r)} )\bigr) \sqcup \bigl( Z(\alpha^p) \setminus Z(\alpha^p\beta) \bigr)$

\item \label{lem positivity one (ii)} $Z(\alpha^ \ell)
= \bigl( \sqcup_{r=1}^{k_{\ell +1}} Z(\alpha^ \ell \gamma_{\ell +1}^{(r)}) \bigr) \sqcup Z(\alpha^ \ell\beta) \sqcup \bigl( Z(\alpha^{\ell +1}) \setminus Z(\alpha^{\ell +1}\beta) \bigr)$.

\end{enumerate}

\end{Lemma}

\begin{proof}
\eqref{lem positivity one (i)} $\subseteq$: Let $x \in Z(\alpha^ \ell) \setminus Z(\alpha^ \ell\beta)$. Write the normal form of $x$ as $\theta_1 \theta_2 \cdots$. Then $\theta_1 =\alpha^{q-1}$ for some $q > \ell$. If $q \le p$ then $x \in Z(\alpha^{q-1} \gamma_q^{(r)})$ for some $r$. If $q > p$ then $x \in Z(\alpha^p) \setminus Z(\alpha^p\beta)$. The reverse containment is clear.

\noindent
\eqref{lem positivity one (ii)} $\subseteq$: Let $x \in Z(\alpha^\ell)$. Write $x = \theta_1 \theta_2 \cdots$ in normal form. Then $\theta_1 \in\alpha^\ell \Lambda_1$. Write $\theta_1 =\alpha^{\ell + m} \beta^n$ with $m$, $n \ge 0$. If $m = n = 0$ then $\theta_2 \in\gamma_{\ell + 1}^{(r)} \Lambda_2$ for some $r$, and then $x \in Z(\alpha^\ell \gamma_{\ell + 1}^{(r)})$. If $n > 0$ then $\theta_1 \in Z(\alpha^\ell \beta)$. Finally, if $n = 0$ and $m > 0$ then $x \in Z(\alpha^\ell) \setminus Z(\alpha^\ell \beta)$. The reverse containment is clear.
\end{proof}

\begin{Lemma} \label{lem positivity two}

Let $p \ge i$. Then
\[
[\chi_{Z(\alpha^i)}]_0
= \sum_{\ell = i+1}^p (\ell - i) \sum_{r=1}^{k_\ell} [\chi_{Z(\alpha^{\ell-1}\gamma_\ell^{(r)})}]_0
+ (p - i) [\chi_{Z(\alpha^p) \setminus Z(\alpha^p\beta)}]_0
+ [\chi_{Z(\alpha^p)}]_0.
\]

\end{Lemma}

\begin{proof}
By Lemma \ref{lem equivalent projections} we have that $[\chi_{Z(\alpha^i\beta)}]_0 = [\chi_{Z(\alpha^{i+1})}]_0$. Then Lemma \ref{lem positivity one}\eqref{lem positivity one (ii)} implies that
\begin{align*}
[\chi_{Z(\alpha^i)}]_0
&= \sum_{r=1}^{k_{i+1}} [\chi_{Z(\alpha^i \gamma_{i+1}^{(r)})}]_0 + [\chi_{Z(\alpha^{i+1})}]_0 + [\chi_{Z(\alpha^{i+1}) \setminus Z(\alpha^{i+1}\beta)}]_0. \\
\noalign{Now we may appy this fact to the middle term, repeatedly, to obtain}
&= \sum_{\ell = i+1}^p \sum_{r=1}^{k_\ell} [\chi_{Z(\alpha^{\ell-1}\gamma_\ell^{(r)})}]_0
+ \sum_{\ell=i+1}^p [\chi_{Z(\alpha^\ell) \setminus Z(\alpha^\ell\beta)}]_0
+ [\chi_{Z(\alpha^p)}]_0. \displaybreak[0] \\
\noalign{ Now using Lemma \ref{lem positivity one}\eqref{lem positivity one (i)} we obtain}
&= \sum_{\ell = i+1}^p \sum_{r=1}^{k_\ell} [\chi_{Z(\alpha^{\ell-1}\gamma_\ell^{(r)})}]_0 \displaybreak[0] 
+ \sum_{\ell=i+1}^p \bigl( \sum_{q = \ell+1}^p \sum_{r = 1}^{k_q} [\chi_{Z(\alpha^{q-1} \gamma_q^{(r)})}]_0 + [\chi_{Z(\alpha^p) \setminus Z(\alpha^p\beta)}]_0 \bigr) \displaybreak[0] \\
&\hspace*{.5 in} + [\chi_{Z(\alpha^p)}]_0 \\
&= \sum_{\ell = i+1}^p \sum_{r=1}^{k_\ell} [\chi_{Z(\alpha^{\ell-1}\gamma_\ell^{(r)})}]_0
+ \sum_{q=i+2}^p \sum_{\ell = i+1}^{q-1} \sum_{r = 1}^{k_q} [\chi_{Z(\alpha^{q-1} \gamma_q^{(r)})}]_0 \displaybreak[0] \\
&\hspace*{.5 in} + (p - i) [\chi_{Z(\alpha^p) \setminus Z(\alpha^p\beta)}]_0 + [\chi_{Z(\alpha^p)}]_0.
\end{align*}
The second term can be rewritten:
\begin{align*}
\sum_{q=i+2}^p \sum_{\ell = i+1}^{q-1} \sum_{r = 1}^{k_q} [\chi_{Z(\alpha^{q-1} \gamma_q^{(r)})}]_0
&= \sum_{q=i+2}^p (q - i - 1) \sum_{r=1}^{k_q} [\chi_{Z(\alpha^{q-1}\gamma_q^{(r)})}]_0 \\
&= \sum_{q=i+1}^p (q - i - 1) \sum_{r=1}^{k_q} [\chi_{Z(\alpha^{q-1}\gamma_q^{(r)})}]_0 \\
&= \sum_{\ell=i+1}^p (\ell - i - 1) \sum_{r=1}^{k_\ell} [\chi_{Z(\alpha^{\ell-1}\gamma_\ell^{(r)})}]_0.
\end{align*}
Making this replacement in the previous calculation, we find that
\[
[\chi_{Z(\alpha^i)}]_0
= \sum_{\ell = i+1}^p (\ell - i) \sum_{r=1}^{k_\ell} [\chi_{Z(\alpha^{\ell-1}\gamma_\ell^{(r)})}]_0
+ (p - i) [\chi_{Z(\alpha^p) \setminus Z(\alpha^p\beta)}]_0 + [\chi_{Z(\alpha^p)}]_0. \qedhere
\]
\end{proof}

\begin{proof} \textit{(of Theorem \ref{thm positive elements}.)}
First suppose that $c_{r,\ell,j} + m + n(\ell - i - 1) \ge 0$ for all $\ell$, $r$ and $j$. Choose $p > i + 1$ such that $c_{\ell,r,j} = 0$ for $\ell > p$. Then for $\ell > p$ we have $m + n(\ell - i - 1) \ge 0$, hence $n \ge -\frac{m}{\ell-i-1}$. Letting $\ell \to \infty$ we see that $n \ge 0$. Now we use the fact that for each $\ell$, $\{F_{\ell,j} : 1 \le j \le h_\ell\}$ is a partition of $X_{\ell+1}$, and also Lemmas \ref{lem positivity one}\eqref{lem positivity one (i)} and \ref{lem positivity two}, to obtain
\begin{align*}
a &= \sum_{\ell > i} \sum_{r = 1}^{k_\ell} \sum_{j = 1}^{h_\ell} c_{\ell,r,j} [\chi_{\alpha^{\ell-1} \gamma_\ell^{(r)} F_{\ell,j}}]_0 + m [\chi_{Z(\alpha^i \setminus \alpha^i\beta)}]_0 + n [\chi_{Z(\beta^{i+1})}]_0 \displaybreak[0] \\
&= \sum_{\ell = i+1}^p \sum_{r = 1}^{k_\ell} \sum_{j = 1}^{h_\ell} c_{\ell,r,j} [\chi_{\alpha^{\ell-1} \gamma_\ell^{(r)} F_{\ell,j}}]_0 \displaybreak[0] \\
 &\hspace*{.5 in} + m \bigl( \sum_{\ell = i+1}^p \sum_{r=1}^{k_\ell} [\chi_{Z(\alpha^{\ell - 1} \gamma_\ell^{(r)})}]_0 + [\chi_{Z(\alpha^p) \setminus Z(\alpha^p\beta)}]_0 \bigr) \displaybreak[0] \\
 &\hspace*{1 in} + n \bigl( \sum_{\ell = i+2}^p (\ell - i - 1) \sum_{r=1}^{k_{\ell}} [\chi_{Z(\alpha^{\ell-1}\gamma_\ell^{(r)})}]_0
+ (p - i - 1) [\chi_{Z(\alpha^p) \setminus Z(\alpha^p\beta)}]_0
+ [\chi_{Z(\alpha^p)}]_0 \displaybreak[0] \bigr) \\
&= \sum_{\ell = i+1}^p \sum_{r = 1}^{k_\ell} \sum_{j=1}^{h_\ell} (c_{\ell,r,j} + m + n(\ell - i - 1)) [\chi_{\alpha^{\ell-1} \gamma_\ell^{(r)} F_{\ell,j}}]_0 \displaybreak[0] \\
 &\hspace*{.5 in} + (m + n(p - i - 1)) [\chi_{Z(\alpha^p) \setminus Z(\alpha^p\beta)}]_0
 + n [\chi_{Z(\alpha^p)}]_0 \\
&\ge 0.
\end{align*}

Now suppose that $a \ge 0$. Fix $x \in U_i$. There is $\ell(x) > i$ such that $x \in E_{\ell(x)}$. Then there are $\mu(x) \in \Omega_{\ell(x)}$, $1 \le t(x) \le k_{\ell(x)}$, and $x' \in X_{\ell(x) + 1}$ such that $x = \mu(x) \gamma_{\ell(x)}^{(t(x))} x'$. Let $1 \le j(x) \le h_{\ell(x)}$ with $x' \in F_{\ell(x),j(x)}$. Recalling Notation \ref{notation C*(Gi,Ui)}, define $\pi_x : C^*(G_i|_{U_i}) \to M_{\Omega_{\ell(x)} \times \Omega_{\ell(x)}}$ by $\pi_x(f) = f_{\ell(x),t(x)} (x')$. Since $C^*(G_i|_{U_i})$ is an ideal in $C^*(G_i)$, $\pi_x$ extends uniquely to a $*$-homomorphism $\widetilde{\pi}_x : C^*(G_i) \to M_{\Omega_{\ell(x)} \times \Omega_{\ell(x)}}$. For $\nu \in \Omega_{\ell(x)}$ let $p_\nu = \chi_{\nu \gamma_{\ell(x)}^{(t(x))} X_{\ell(x) + 1}}\in C(X)$. Then $\widetilde{\pi}_x(p_\nu) = \pi_x(p_\nu) = e_{\nu,\nu} \in M_{\Omega_{\ell(x)},\Omega_{\ell(x)}}$. We calculate $\widetilde{\pi}_{x\,*}(a)$. First, for any $\ell > i$, $r$, and $j$,
\[
\widetilde{\pi}_x ( \chi_{\alpha^{\ell-1} \gamma_\ell^{(r)} F_{\ell,j}})
= \begin{cases}
 e_{\alpha^{\ell(x) - 1},\alpha^{\ell(x) - 1}}, &\text{ if } \ell=\ell(x), r=t(x), j=j(x) \\
0, &\text{ otherwise.}
 \end{cases}
\]
Next we let $\nu \in \Omega_{\ell(x)}$ and calculate (in $C^*(G_i|_{U_i})$).
\begin{align*}
p_\nu \cdot \chi_{Z(\alpha^i) \setminus Z(\alpha^i\beta)}
&= \chi_{\nu \gamma_{\ell(x)}^{(t(x))} X_{\ell(x) + 1}} \cdot \chi_{Z(\alpha^i) \setminus Z(\alpha^i\beta)} \\
&= \chi_{\nu \gamma_{\ell(x)}^{(t(x))} X_{\ell(x) + 1} \cap (Z(\alpha^i) \setminus Z(\alpha^i\beta))} \\
&= \delta_{\nu,\alpha^{\ell(x)-1}} \chi_{\alpha^{\ell(x)-1} \gamma_{\ell(x)}^{(t(x))} X_{\ell(x) + 1}} \\
&= \delta_{\nu,\alpha^{\ell(x)-1}} p_{\alpha^{\ell(x)-1}}.
\end{align*}
Now we have
\begin{align*}
\widetilde{\pi}_x(\chi_{Z(\alpha^i) \setminus Z(\alpha^i\beta)})
&= \sum_{\nu \in \Omega_{\ell(x)}} \pi_x (p_\nu \cdot \chi_{Z(\alpha^i) \setminus Z(\alpha^i\beta)}) \\
&= \pi_x(p_{\alpha^{\ell(x)-1}}) \\
&= e_{\alpha^{\ell(x) - 1},\alpha^{\ell(x) - 1}}.
\end{align*}
Next note that $p_\nu \cdot \chi_{Z(\beta^{i+1})} = \chi_{\nu \gamma_{\ell(x)}^{(t(x))} X_{\ell(x) + 1} \cap Z(\beta^{i+1})}$. The intersection in the last expression is nonempty only if $\nu = \beta^{i+1} \alpha^c \beta^d$, where $c + d = \ell(x) - i - 2$. There are $\ell(x) - i - 1$ choices where this happens. Therefore $\widetilde{\pi}_x(\chi_{Z(\beta^i)})$ has rank $\ell(x) - i - 1$ as a projection in $M_{\Omega_{\ell(x)},\Omega_{\ell(x)}}$.

Putting these three calculations together we obtain
\[
0 \le
\widetilde{\pi}_{x\,*}(a)
= \bigl( c_{\ell(x),t(x),j(x)} + m + n(\ell(x) - i - 1) \bigr) [e_{\alpha^{\ell(x) - 1},\alpha^{\ell(x) - 1}}]_0.
\]
Since for any choice of $\ell$, $r$, $j$ there exists $x \in U_i$ with $\ell=\ell(x)$, $r=t(x)$ and $j = j(x)$, it follows that $c_{\ell,r,j} + m + n(\ell - i - 1) \ge 0$ for all $\ell$, $r$ and $j$.
\end{proof}

Theorem \ref{thm positive elements} and its proof have a corollary which we will need later (namely, in the proof of Theorem \ref{thm K-theory of G}).

\begin{Corollary} \label{cor positive elements}
	
Let $a \in K_0(C^*(G_i))$ with
\begin{equation} \label{eqn simple positive elements}
a = m [\chi_{Z(\alpha^i) \setminus Z(\alpha^i\beta)}]_0 + n [\chi_{Z(\beta^{i+1})}]_0.
\end{equation}
Then $a \ge 0$ if and only if $m, n \ge 0$.
	
\end{Corollary}

\begin{proof}
If $m, n, \ge 0$, then it is clear from Theorem \ref{thm positive elements} that $a \ge 0 $ (noting that $c_{\ell, r, j} = 0$ for all $\ell, r, j$). Conversely, suppose that $a \ge 0$. In the first lines of the proof of Theorem \ref{thm positive elements} it is shown that $n \ge 0$. Letting $\ell = i+1$, the theorem implies that $m \ge 0$.
\end{proof}

\section{$K$-theory of $C^*(G)$}
\label{section K-theory of G}

To analyze $C^*(G)$, we begin with the following observation:

\begin{Theorem} \label{thm inductive sequence}
	
$C^*(G)$ is the limit of the inductive sequence 
\begin{equation} \label{eqn inductive sequence of G i}
C^*(G_1) \rightarrow C^*(G_2) \rightarrow \cdots
\end{equation}
where the connecting maps are induced from the inclusion maps $C_c(G_i) \hookrightarrow C_c(G_{i+1})$.

\end{Theorem}

\begin{proof}
For $\mu$, $\nu \in \Lambda$ we say that $\mu$ and $\nu$ \textit{have a common suffix} if $\mu = \mu'\theta$ and $\nu = \nu' \theta$, where $\theta \not\in \Lambda^0$. Let
\[
S = \{(\mu,\nu) \in v_1 \Lambda \times v_1 \Lambda : s(\mu) = s(\nu) \text{ and } \mu, \nu \text{ do not have a common suffix} \}.
\]
Thus $G = \{ [\mu,\nu,x] : (\mu,\nu) \in S, x \in s(\mu)\partial \Lambda \}$. We define a subset $S_i \subseteq S$ as follows. Recall the set $\Phi_i$ from Definition \ref{def Phi_i}. Let
\[
S_i = \{(\mu\theta,\nu\phi) \in S : \mu,\nu \in \Phi_i,\ \theta,\phi \in \Lambda_1,\ d(\theta) \perp d(\phi) \}.
\]
(Note that since a pair in $S$ does not have a common suffix, if $\phi=\theta=s(\mu)$ then $|\mu| = |\nu|$ and $\mu_{|\mu|} \not= \nu_{|\nu|} \}$.) Then $G_i = \{[\xi,\eta,x] : (\xi,\eta) \in S_i \}$. Note that for each $(\mu,\nu) \in S$, the set $[\mu,\nu,s(\mu)\partial\Lambda]$ is a compact open subset of $G$. Therefore $G_i = \cup_{(\mu,\nu) \in S_i} [\mu,\nu, s(\mu) \partial \Lambda]$ is an open subset of $G_{i+1}$, and of $G$, and also $G_{i+1} \setminus G_i = \cup_{(\mu,\nu) \in S_{i+1} \setminus S_i} [\mu,\nu,s(\mu)\partial\Lambda]$ is an open subset of $G_{i+1},$ and similarly, $G \setminus G_i$ is an open subset of $G$. Therefore $G_i$ is a clopen subgroupoid of $G_{i+1}$ and of $G$. By (the proof of) \cite[Theorem 14.2]{spiel2}, the inclusion $C_c(G_i) \hookrightarrow C_c(G_{i+1})$ induces an injective $*$-homomorphism $C^*(G_i) \hookrightarrow C^*(G_{i+1})$. Since $C_c(G) = \cup_i C_c(G_i)$, the limit of the system $C^*(G_i) \to C^*(G_{i+1})$ equals $C^*(G)$.
\end{proof}

{
\begin{Remark} \label{remark G amenable}

Since $G_i$ is a clopen subgroupoid of $G$ (as seen in the proof of Theorem \ref{thm inductive sequence}) and is amenable by Proposition \ref{prop G i amenable}, it follows from \cite[Proposition 9.84]{will} that $G$ is amenable. In particular, all $C^*$-algebras we discuss are nuclear. Moreover since the groupoids we discuss are \'etale, it follows from \cite{tu} that the $C^*$-algebras satisfy the UCT.

\end{Remark}

We will see shortly that $K_0(C^*(G))\cong \mathbb{Z}^2$ (Theorem \ref{thm K-theory of G}), but to prove this we will need some preliminary results. (Recall that a partition \textit{refines} a set $S$ if $S$ equals the union of a subcollection of the partition.)

\begin{Lemma} \label{lem partition lemma}

For $m \ge 1$ and $n\geq 0$, let
$W_{m,n} = \bigcup_{p+q = n} Z(v_m \alpha^p \beta^q)$. Define collections of subsets of $W_{m,n}$ by
\begin{align*}
	P_{m,n} ^{(1)} & =  \{ Z(v_m \alpha^{n+1}\beta^j)\setminus Z(v_m \alpha^{n+1}\beta^{j+1}): j\leq n\}\\
	P_{m,n} ^{(2)} & =  \{ Z(v_m \alpha^{i}\beta^{n+1})\setminus Z(v_m \alpha^{i+1}\beta^{n+1}): i\leq n\}\\
	P_{m,n} ^{(3)} & =  \{ Z(v_m \alpha^{i}\beta^{j}\lambda): i,j \le n \le i + j, \lambda \in v_{m+i+j} \Lambda_2^1 \}\\
	P_{m,n} ^{(4)} & =  \{ Z(v_m \alpha^{n+1}\beta^{n+1})\}\\
	P_{m,n} & =  \cup_{r=1}^4 P_{m,n} ^{(r)}.
\end{align*}
Then $P_{m,n}$ is a partition of $W_{m,n}$ that refines $Z(\nu)$ and $Z(\nu) \setminus Z(\nu\lambda)$ for every $\nu \in v_m \Lambda_1$ with $|\nu| = n$ and every $\lambda \in s(\nu)\Lambda$ with $|\lambda| = 1$.

\end{Lemma}

\begin{proof}
First we show that $W_{m,n} = \cup P_{m,n}$. It is clear that $\cup P_{m,n} \subseteq W_{m,n}$. To show the other inclusion, fix $x \in W_{m,n}$ and let $p = \min\{h: x_h \in \Lambda_2 \}$. Then $n+1 \leq p \leq \infty$. We have several cases to consider.
\begin{enumerate}
\item Suppose $p\leq 2n+1$. Then $x\in Z(v_m \alpha^i\beta^j x_p)$ for some $i$ and $j$ with $n \le i + j \le 2n$.
\begin{enumerate}
\item If $i>n$, then $j<n$. Then $Z(v_m \alpha^i\beta^j x_p)\subseteq Z(v_m \alpha^{n+1}\beta^j)\setminus Z(v_m \alpha^{n+1}\beta^{j+1})$, which is a set in $P_{m,n}^{(1)}$.
\item If $j>n$, then $i<n$. Then $Z(v_m \alpha^i\beta^j x_p)\subseteq Z(v_m \alpha^i\beta^{n+1})\setminus Z(v_m \alpha^{i+1}\beta^{n+1})$, which is a set in $P_ {m,n}^{(2)}$.
\item If $i, j\leq n$, then $Z(v_m \alpha^i\beta^j x_p) \in P_{m,n}^{(3)}$.
\end{enumerate}
\item Suppose $p > 2n+1$. Note that $Z(v_m \alpha^{n+1}\beta^{n+1}) \in P_{m,n}^{(4)}$. If $x\notin Z(v_m \alpha^{n+1}\beta^{n+1})$ then there are two possibilities. We know that $x_1 \cdots x_{p-1} \in \Lambda_1$. Let $\ell = d(x_1 \cdots x_{p-1}) \in \IN^2$.
\begin{enumerate}
\item If $\ell_2 \leq n$, then $x\in Z(v_m \alpha^{n+1}\beta^{\ell_2})\setminus Z(v_m \alpha^{n+1}\beta^{\ell_2+1})$, a set in $P_{m,n}^{(1)}$.
\item If $\ell_1 \leq n$, then $x\in Z(v_m \alpha^{\ell_1}\beta^{n+1})\setminus Z(v_m \alpha^{\ell_1+1}\beta^{n+1})$, a set in $P_{m,n}^{(2)}$.
\end{enumerate}
\end{enumerate}
Thus in all cases we have $x\in \cup P_{m,n}$ so $W_{m,n} = \cup P_{m,n}$.
	
Now we show that the elements of $P_{m,n}$ are pairwise disjoint. Consider the normal form of an element $x \in W_{m,n}$: $x = \theta_1 \theta_2 \cdots$, where $\theta_1 \in \Lambda_1$ and $|\theta_1| \ge n$. Let $d(\theta_1) = \ell$. From the definition of the sets $P_{m,n}^{(r)}$ we see that

$x \in \cup P_{m,n}^{(1)}$ if and only if $\ell_1 \ge n+1$ and $\ell_2 \le n$

$x \in \cup P_{m,n}^{(2)}$ if and only if $\ell_1 \le n$ and $\ell_2 \ge n + 1$

$x \in \cup P_{m,n}^{(3)}$ if and only if $\ell_1, \ell_2 \le n \le \ell_1 + \ell_2$

$x \in \cup P_{m,n}^{(4)}$ if and only if  $\ell_1,\ \ell_2 \ge n+1$.

\noindent
It follows that for any $1 \le r_1 < r_2 \le 4$, each element of $P_{m,n}^{(r_1)}$ is disjoint from each element of $P_{m,n}^{(r_2)}$. The elements of $P_{m,n}^{(1)}$ with different values of $\ell_2$ are disjoint, and similarly the elements of $P_{m,n}^{(2)}$ with different values of $\ell_1$ are disjoint. The same is true for $P_{m,n}^{(3)}$, for the same reason. Therefore $P_{m,n}$ is a partition of $W_{m,n}$.

Next let $\nu\in v_m \Lambda_1$ with $|\nu| = n$. Let $\nu = v_m \alpha^p\beta^q$ (with $p + q = n$), so that $Z(\nu) = Z(v_m \alpha^p)\cap Z(v_m \beta^q)$. We will show that $P_{m,n}$ refines $Z(\nu)$. Let $S\in P_{m,n}$ be such that $S\cap Z(\nu)\neq \varnothing$. We will show that $S\subseteq Z(\nu)$.

First suppose $S\in P_{m,n}^{(1)}$, so $S = Z(v_m \alpha^{n+1}\beta^j)\setminus Z(v_m \alpha^{n+1}\beta^{j+1})$ for some $j\leq n$. Then $S = Z(v_m \alpha^{n+1})\cap (Z(v_m \beta^j) \setminus Z(v_m \beta^{j+1}))$. Since $Z(\nu)\cap S\neq \varnothing$, we know $Z(v_m\beta^q)\cap (Z(v_m \beta^j)\setminus Z(v_m \beta^{j+1}))\neq \varnothing$, hence $q\leq j$. Then $Z(v_m \beta^j)\setminus Z(v_m \beta^{j+1})\subseteq Z(v_m \beta^q)$, and since $p\leq n$, $Z(v_m \alpha^{n+1})\subseteq Z(v_m \alpha^p)$. Therefore $S = Z(v_m \alpha^{n+1}) \cap (Z(v_m \beta^j)\setminus Z(v_m \beta^{j+1}))\subseteq Z(v_m \alpha^p)\cap Z(v_m \beta^q) = Z(\nu)$. It follows that $Z(v_m \alpha^{n+1}\beta^j)\setminus Z(v_m \alpha^{n+1}\beta^{j+1})$ is contained in $Z(\nu)$ if $q \le j$ and is disjoint from $Z(\nu)$ if $j < q$. If $S\in P_{m,n}^{(2)}$, an analogous argument shows that $Z(v_m \alpha^i\beta^{n+1})\setminus Z(v_m \alpha^{i+1}\beta^{n+1})$ is contained in $Z(\nu)$ if $p \le i$ and is disjoint from $Z(\nu)$ if $i < p$.

Now suppose $S\in P_{m,n}^{(3)}$ so $S = Z(v_m \alpha^i\beta^j \lambda)$ with $i, j \le n \le i+j$, and $\lambda \in v_{m+i+j} \Lambda_2^1$. Then $S \subseteq (Z(v_m \alpha^i)\setminus Z(v_m \alpha^{i+1}))\cap (Z(v_m \beta^j)\setminus Z(v_m \beta^{j+1}))$. Since $Z(\nu)\cap S\neq \varnothing$, $p\leq i$ and $q\leq j$, and in this case, $S\subseteq Z(\nu)$. Thus $Z(v_m \alpha^i \beta^j \lambda)$ is contained in $Z(\nu)$ if $p \le i \le n$ and $q \le j \le n$, and is disjoint from $Z(\nu)$ if $i < p$ or if $j < q$.

Finally, since $p, q < n+1$, we have $Z(v_m \alpha^{n+1}\beta^{n+1})\subseteq Z(\nu)$. Therefore $P_{m,n}$ refines $Z(\nu)$.

Lastly, let $\nu = v_m \alpha^p \beta^q$ with $p + q = n$ as above, and let $\lambda \in s(\nu)\Lambda$ with $|\lambda| = 1$. We will show that $P_{m,n}$ refines $E = Z(\nu) \setminus Z(\nu\lambda)$. We first consider the case that $\lambda = \alpha_{m+n}$. Then $E = Z(v_m \alpha^p \beta^q) \setminus Z(v_m \alpha^{p+1} \beta^q)$. Since $p \le n$ we know that $E$ is disjoint from every element of $P_{m,n}^{(1)} \cup P_{m,n}^{(4)}$. Let $S \in P_{m,n}^{(2)}$ with $S \cap E \not= \varnothing$. Then it must be the case that $S = Z(v_m \alpha^p \beta^{n+1}) \setminus Z(v_m \alpha^{p+1} \beta^{n+1})$, and hence $S \subseteq E$. Now we consider elements of $P_{m,n}^{(3)}$. Since $p,q \le n = p + q$, $S = Z(v_m \alpha^p \beta^j \gamma_{m+p+j}^{(r)}) \in P_{m,n}^{(3)}$ for $q \le j \le n$ and $1 \le r \le k_{m+p+j}$. Moreover for such $S$, $S \subseteq E$. These are the only elements $S \in P_{m,n}^{(3)}$ with $S \cap E \not= \varnothing$. Therefore $P_{m,n}$ refines $E$. The case where $\lambda = \beta_{m+n}$ is analogous. Let $\lambda = \gamma_{m+n}^{(r)}$ for some $1 \le r \le k_{m+n}$. Then $E = Z(v_m \alpha^p \beta^q) \setminus Z(v_m \alpha^p \beta^q \gamma_{m+n}^{(r)})$. Note that $Z(v_m \alpha^p \beta^q \gamma_{m+n}^{(r)})$ is one of the sets from $P_{m,n}^{(3)}$ contained in $Z(\nu)$. Therefore $E$ equals the union of the other sets from $P_{m,n}$ that are contained in $Z(\nu)$.
%Let $S \in P_{m,n}^{(1)}$. Then there is $j \le n$ such that $S = Z(v_m \alpha^{n+1} \beta^j) \setminus Z(v_m \alpha^{n+1} \beta^{j+1})$. If $j < q$ then $S \cap E = \varnothing$, while if $j \ge q$ then $S \subseteq E$. Now let $S \in P_{m,n}^{(2)}$, so $S = Z(v_m \alpha^i \beta^{n+1}) \setminus Z(v_m \alpha^{i+1} \beta^{n+1})$ for some $i \le n$. If $i < p$ then $S \cap E = \varnothing$, while if $i \ge p$ then $S \subseteq E$. Next let $S \in P_{m,n}^{(3)}$, so $S = Z(v_m \alpha^i \beta^j \gamma_{m+i+j}^{(s)})$ for some $i,j \le n \le i+j$ and $1 \le s \le k_{m+i+j}$. We see that if $p \le i \le n$, $q \le j \le n$, and $s \not= r$ if $p = i$ and $q = j$, then $S \subseteq E$, while $S \cap E = \varnothing$ otherwise.  Finally, since $p$, $q \le n$ we have that $Z(v_m \alpha^{n+1} \beta^{n+1}) \cap E = \varnothing$.
Therefore in all cases, $P_{m,n}$ refines $Z(\nu) \setminus Z(\nu\lambda)$.
\end{proof}

\begin{Remark} \label{rmk refinements}

It will be convenient to list the explicit refinements given in the proof of Lemma \ref{lem partition lemma}. Let $p + q = n$.
\begin{align}
Z(v_m \alpha^p \beta^q) 
&= \bigsqcup_{j=q}^n ( Z(v_m \alpha^{n+1} \beta^j) \setminus Z(v_m \alpha^{n+1} \beta^{j+1}) ) \label{rmk refinements one} \\
&\quad \sqcup \bigsqcup_{i=p}^n ( Z(v_m \alpha^i \beta^{n+1}) \setminus Z(v_m \alpha^{i+1} \beta^{n+1}) ) \notag \\
&\qquad \sqcup \bigsqcup_{\substack{p \le i \le n \\ q \le j \le n }}
\bigsqcup_{r=1}^{k_{m+i+j}} Z(v_m \alpha^i \beta^j \gamma_{m+i+j}^{(r)}) \notag \\
&\qquad \sqcup Z(v_m \alpha^{n+1} \beta^{n+1}) \notag %\\ \notag
\end{align}
\begin{align}
Z(v_m \alpha^p \beta^q) \setminus Z(v_m \alpha^{p+1} \beta^q)
&= ( Z(v_m \alpha^p \beta^{n+1}) \setminus Z(v_m \alpha^{p+1} \beta^{n+1}) ) \label{rmk refinements two} \\
&\quad \sqcup \bigsqcup_{j=q}^n \bigsqcup_{r=1}^{k_{m+p+j}} Z(v_m \alpha^p \beta^j \gamma_{m+p+j}^{(r)}) \notag %\\ \notag
\end{align}
\begin{align}
Z(v_m \alpha^p \beta^q) \setminus Z(v_m \alpha^p \beta^{q+1})
&= ( Z(v_m \alpha^{n+1} \beta^q) \setminus Z(v_m \alpha^{n+1} \beta^{q+1}) ) \label{rmk refinements three} \\
&\quad \sqcup \bigsqcup_{i=p}^n \bigsqcup_{r=1}^{k_{m+i+q}} Z(v_m \alpha^i \beta^q \gamma_{m+i+q}^{(r)}) \notag %\\ \notag
\end{align}
\begin{align}
Z(v_m \alpha^p \beta^q) \setminus Z(v_m \alpha^p \beta^q \gamma_{m+n}^{(r)})
&= \bigsqcup_{j=q}^n ( Z(v_m \alpha^{n+1} \beta^j) \setminus Z(v_m \alpha^{n+1} \beta^{j+1}) ) \label{rmk refinements four} \displaybreak[0] \\
&\quad \sqcup \bigsqcup_{i=p}^n ( Z(v_m \alpha^i \beta^{n+1}) \setminus Z(v_m \alpha^{i+1} \beta^{n+1}) ) \notag \displaybreak[0] \\
&\qquad \sqcup \bigsqcup_{\substack{p \le i \le n \\ q \le j \le n \\ (i,j) \not= (p,q) }} \bigsqcup_{s=1}^{k_{m+i+j}} Z(v_m \alpha^i \beta^j \gamma_{m+i+j}^{(s)}) \notag \displaybreak[0] \\
&\quad \sqcup \bigsqcup_{\substack{s=1 \\ s \not= r}}^{k_{m+n}} Z(v_m \alpha^p \beta^q \gamma_{m+n}^{(s)}) \notag \\
&\qquad \sqcup Z(v_m \alpha^{n+1} \beta^{n+1}). \notag
\end{align}

\end{Remark}

\begin{Notation}

If $\lambda \in \Lambda$ and $E \subseteq s(\lambda) \partial \Lambda$ then we write $\lambda E := \{ \lambda x : x \in E \}$. If $P$ is a family of subsets of $s(\lambda) \partial \Lambda$ we write $\lambda P := \{ \lambda E : E \in P \}$.

\end{Notation}

\begin{Proposition} \label{prop partition proposition}

For $n \geq 1$ let
$Q_n = \bigcup_{\mu \in \Phi_n} \mu P_{|\mu|+1, n - |\mu|}$, where $P_{r,s}$ is as in Lemma \ref{lem partition lemma} and $\Phi_n$ is as in Definition \ref{def Phi_i}. Then $Q_n$ is a partition of $X$ that refines $Z(\nu)$ and $Z(\nu) \setminus Z(\nu\lambda)$ for all $\nu\in v_1\Lambda$ with $|\nu|\leq n$ and $\lambda \in s(\nu)\Lambda$ with $|\lambda| = 1$.

\end{Proposition}

\begin{proof}
Let $n \ge 1$. We first show that the sets in $Q_n$ are pairwise disjoint. Let $\mu$, $\nu \in \Phi_n$ with $\mu \not= \nu$. Let $x \in \cup \mu P_{|\mu|+1, n - |\mu|}$ and $y \in \cup \nu P_{|\nu|+1, n - |\nu|}$. Then $\mu \in x$ and $\nu \in y$, and $x_{|\mu|+1} \cdots x_n$, $y_{|\nu|+1} \cdots y_n \in \Lambda_1$. Without loss of generality suppose that $|\mu| \le |\nu|$. Then $y_{|\nu|} \in \Lambda_2$ and $y_{|\nu|} = \nu_{|\nu|}$. If $|\mu| < |\nu|$ then $x_{|\nu|} \in \Lambda_1$, and hence $x_{|\nu|} \not= \nu_{|\nu|}$, hence $x \not= y$. If $|\mu| = |\nu|$ and $\mu_{|\nu|} = \nu_{|\nu|}$ then since $\mu \not= \nu$, $x \not= y$ by Lemma \ref{lemma:common extension}. Thus if $S \in P_{|\mu|+1, n - |\mu|}$ and $T \in P_{|\nu|+1, n - |\nu|}$ then $\mu S \cap \nu T = \varnothing$. Since $P_{|\mu|+1, n - |\mu|}$ is a pairwise disjoint family, so is $\mu P_{|\mu|+1, n - |\mu|}$. Therefore $Q_n$ is pairwise disjoint.

Next we show that $\cup Q_n = X$. Let $x \in X$. Write $x = x_1 x_2 \cdots$ as an infinite word. Let $p = \max \{ j \le n : x_j \in \Lambda_2 \}$, and let $\mu = x_1 \cdots x_p$. Then $\mu \in \Phi_n$ and $x_{p+1} \cdots x_n \in \Lambda_1$, and hence $x \in \cup \mu P_{p+1,n-p} \subseteq \cup Q_n$.

Finally, let $\nu \in v_1 \Lambda_1$ with $|\nu| \le n$. Write $\nu = \mu \theta$ where $\mu \in \Phi_n$ and $\theta \in \Lambda_1$. By Lemma \ref{lem partition lemma}, $Z(v_{|\mu|+1}\theta)$ and $Z(v_{|\mu|+1} \theta) \setminus Z(v_{|\mu|+1} \theta \lambda)$ are refined by $P_{|\mu|+1,|\theta|}$. Then $Z(\nu) = Z(\mu\theta)$ and $Z(\nu) \setminus Z(\nu\lambda) = Z(\mu\theta) \setminus Z(\mu\theta\lambda)$ are refined by $\mu P_{|\mu|+1,|\theta|}$.
\end{proof}

Recall from Theorem \ref{K-theory of G i} and Corollary \ref{cor typical element of K-theory} that
$$K_0(C^*(G_i)) \cong \bigoplus_{\ell > i} C(X_{\ell+1},\IZ)^{k_\ell} \oplus \IZ^2$$
where a typical generator of the $\ell^{th}$ summand is $[\chi_{\alpha^{\ell-1}\gamma_\ell^{(r)} F}]_0$, where $F \subseteq X_{\ell+1}$ is a compact open subset, and the generators of the right summand are $[\chi_{Z(\alpha^i) \setminus Z(\alpha^i\beta)}]_0$ and $[\chi_{Z(\beta^{i+1})}]_0$.

\begin{Lemma} \label{lem connecting maps in K 0 two}
	
The induced map $K_0(C^*(G_i))\rightarrow K_0(C^*(G_{i+1}))$ carries the summand $\IZ^2$ of $K_0(C^*(G_i))$ into the summand $\mathbb{Z}^2$ of $K_0(C^*(G_{i+1}))$. Using the generators chosen in Theorem \ref{K-theory of G i} for the summand $\IZ^2$, this restriction is implemented by the matrix
$B_i = \begin{psmallmatrix}
k_{i+1} + 1 & 1\\
k_{i+1} & 1
\end{psmallmatrix}$.

\end{Lemma}

\begin{proof}
From Lemma \ref{lem positivity one}, letting $\ell = i$ and $p = i+1$, we have
\begin{align*}
[\chi_{Z(\alpha^i)\setminus Z(\alpha^i\beta)}]_0 & =  [\chi_{Z(\alpha^{i+1})\setminus Z(\alpha^{i+1}\beta)}]_0 + \sum_{r = 1}^{k_{i+1}} [\chi_{Z(\alpha^i\gamma_{i+1}^{(r)})}]_0\\
& = [\chi_{Z(\alpha^{i+1})\setminus Z(\alpha^{i+1}\beta)}]_0 + \sum_{r = 1}^{k_{i+1}} ([\chi_{Z(\alpha^i\gamma_{i+1}^{(r)})\setminus Z(\alpha^i\gamma_{i+1}^{(r)} \beta)}]_0
+ [\chi_{Z(\alpha^i\gamma_{i+1}^{(r)} \beta)}]_0).
\end{align*}
By Lemma \ref{lem equivalent projections} we have
\begin{align*}
[\chi_{Z(\alpha^i\gamma_{i+1}^{(r)})\setminus Z(\alpha^i\gamma_{i+1}^{(r)} \beta)}]_0 & = [\chi_{Z(\alpha^{i+1})\setminus Z(\alpha^{i+1}\beta)}]_0\\
\noalign{and }
[\chi_{Z(\alpha^i\gamma_{i+1}^{(r)} \beta)}]_0 & = [\chi_{Z(\beta^{i+2})}]_0
\end{align*}
in $K_0(C^*(G_{i+1}))$. Therefore
\[
[\chi_{Z(\alpha^{i})\setminus Z(\alpha^{i}\beta)}]_0\mapsto (k_{i+1}+1)[\chi_{Z(\alpha^{i+1})\setminus Z(\alpha^{i+1}\beta)}]_0 + k_{i+1}[\chi_{Z(\beta^{i+2})}]_0.
\]
Finally,
\[
[\chi_{Z(\beta^{i+1})}]_0 = [\chi_{Z(\beta^{i+1})\setminus Z(\beta^{i+2})}]_0 + [\chi_{Z(\beta^{i+2})}]_0
\]
so that
\[
[\chi_{Z(\beta^{i+1})}]_0\mapsto [\chi_{Z(\alpha^{i+1})\setminus Z(\alpha^{i+1}\beta)}]_0 + [\chi_{Z(\beta^{i+2})}]_0
\]
again, using equivalences in $C^*(G_{i+1})$. The matrix of the restriction follows from these formulas.
\end{proof}

\begin{Lemma} \label{lem connecting maps in K 0 one}
	
Let $1 \le i < \ell$, $1 \le r \le k_\ell$, and $F \subseteq X_{\ell + 1}$ a compact open subset. There is $i' > i$ such that the image of $[\chi_{\alpha^{\ell - 1} \gamma_\ell^{(r)} F}]_0$ under the induced map $K_0(C^*(G_i)) \to K_0(C^*(G_{i'}))$ lies in the $\IZ^2$ summand (of $K_0(C^*(G_{i'}))$).

\end{Lemma}

\begin{proof}
We may suppose that $F = Z(\nu) \setminus \cup_{j=1}^m Z(\nu_j)$, where $\nu \in v_{\ell+1} \Lambda$ and $\nu_1$, $\ldots$, $\nu_m$ extend $\nu$. Choose $n \ge \ell + 1 + \max \{|\nu_1|, \ldots, |\nu_m| \}$. By Proposition \ref{prop partition proposition}, we have
\begin{align*}
Z(\alpha^{\ell-1} \gamma_\ell^{(r)} \nu)
&= \sqcup \{ w \in Q_n : w \subseteq Z(\alpha^{\ell-1} \gamma_\ell^{(r)} \nu) \} \\
\noalign{and for each $j$,}
Z(\alpha^{\ell-1} \gamma_\ell^{(r)} \nu_j)
&= \sqcup \{ w \in Q_n : w \subseteq Z(\alpha^{\ell-1} \gamma_\ell^{(r)} \nu_j) \},\ 1 \le j \le m. \\
\noalign{Hence}
\alpha^{\ell-1} \gamma_\ell^{(r)} F
&= Z(\alpha^{\ell-1} \gamma_\ell^{(r)} \nu) \setminus \cup_{j=1}^m Z(\alpha^{\ell-1} \gamma_\ell^{(r)} \nu_j) \\
&= \sqcup \{w \in Q_n : w \subseteq \alpha^{\ell-1} \gamma_\ell^{(r)} F \}.
\end{align*}
Now fix $w \in Q_n$ with $w \subseteq \alpha^{\ell-1} \gamma_\ell^{(r)} F$. Then there is $\mu \in \Phi_n$ such that $w \in \mu P_{|\mu|+1,n-|\mu|}$. Therefore $w$ has one of the following forms:

\begin{enumerate}[(i)]

\item \label{lem K 0 one i} $Z(\mu \alpha^{n-|\mu|+1} \beta^t) \setminus Z(\mu \alpha^{n-|\mu|+1} \beta^{t+1}), \text{ where } t \le n - |\mu|$

\item \label{lem K 0 one ii} $Z(\mu \alpha^s \beta^{n-|\mu|+1}) \setminus Z(\mu \alpha^{s+1} \beta^{n-|\mu|+1}), \text{ where } s \le n - |\mu|$

\item \label{lem K 0 one iii} $Z(\mu \alpha^s \beta^t \lambda), \text{ where } s,t \le n - |\mu| \le s+t,\ \lambda \in v_{|\mu|+s+t+1} \Lambda_2^1$

\item \label{lem K 0 one iv} $Z(\mu \alpha^{n-|\mu|+1} \beta^{n-|\mu|+1})$.

\end{enumerate}

Using Lemma \ref{lem equivalent projections} we see that in case \eqref{lem K 0 one i}, since $|\mu| \le n$, we have that $\chi_w$ is equivalent in $C^*(G_{n+t+1})$ to $\chi_{Z(\alpha^{n+t+1}) \setminus Z(\alpha^{n+t+1}\beta)}$. Case \eqref{lem K 0 one ii} is nearly identical. In case \eqref{lem K 0 one iii} we have that $\chi_w$ is equivalent in $C^*(G_{|\mu|+s+t})$ to $\chi_{Z(\beta^{|\mu|+s+t+1})}$. In case \eqref{lem K 0 one iv} we have that $\chi_w$ is equivalent in $C^*(G_{2n-|\mu|})$ to $\chi_{Z(\beta^{2n-|\mu|+1})}$. Choose $i'$ so large that the equivalences described above occur in $C^*(G_{i'})$. By Lemma \ref{lem connecting maps in K 0 two}, for all $w \in Q_n$ with $w \subseteq \alpha^{\ell-1} \gamma_\ell^{(r)} F$, $[\chi_w]_0$ lies in the summand $\IZ^2$ in $K_0(C^*(G_{i'}))$. Since $\chi_{\alpha^{\ell-1}\gamma_\ell^{(r)} F} = \sum \{ \chi_w : w \in Q_n,\ w \subseteq \alpha^{\ell-1} \gamma_\ell^{(r)} F \}$, the lemma follows.
\end{proof}

We can now give the ordered $K$-theory for $C^*(G)$.

\begin{Theorem} \label{thm K-theory of G}
	
$K_1(C^*(G)) = 0$ and $K_0(C^*(G)) \cong \IZ^2$. Moreover, $K_0(C^*(G))$ is realized as the limit of the inductive sequence $\IZ^2 \to \IZ^2 \to \cdots$, where the $i^{th}$ connecting map is given by the matrix $B_i$ from Lemma \ref{lem connecting maps in K 0 two}, and each term in the sequence has the standard positive cone $\mathbb{N}^2$.

\end{Theorem}

\begin{proof}
By Theorem \ref{K-theory of G i} we know that $K_1(C^*(G_i)) = 0$ for all $i$, and hence that $K_1(C^*(G)) = 0$. Let $\eta_{i,i'} : C^*(G_i) \to C^*(G_{i'})$ and $\eta_i : C^*(G_i) \to C^*(G)$ denote the inclusions of Theorem \ref{thm inductive sequence}. Let $a \in K_0(C^*(G))$. Then there is $i \ge 1$ and $b \in K_0(C^*(G_i))$ such that $a = \eta_{i,*}(b)$. Let
\[
b = \sum_{\ell > i} \sum_{r = 1}^{k_\ell} \sum_{j = 1}^{h_\ell} c_{\ell,r,j} [\chi_{\alpha^{\ell-1} \gamma_\ell^{(r)} F_{\ell,j}}]_0 + m [\chi_{Z(\alpha^i) \setminus Z(\alpha^i\beta)}]_0 + n [\chi_{Z(\beta^{i+1})}]_0
\]
as in \eqref{eqn typical positive element}. Note that this is a finite sum. By Lemmas \ref{lem connecting maps in K 0 two} and \ref{lem connecting maps in K 0 one} there is $i' > i$ such that ${\eta_{i,i'}}_*(b)$ lies in the summand $\IZ^2$ of $K_0(C^*(G_{i'}))$. Let $H_i$ denote the $\IZ^2$ summand of $K_0(C^*(G_i))$. Then $K_0(C^*(G))$ is the inductive limit of the sequence $(H_i, {\eta_{i,i+1}}_*|_{H_i})$. By Lemma \ref{lem connecting maps in K 0 two}, the connecting maps in this sequence are given by the matrices $B_i$, and hence are invertible. Therefore the limit $K_0(C^*(G))$ is isomorphic to $\IZ^2$. The fact that each term in the sequence has standard positive cone follows from Corollary \ref{cor positive elements}.
\end{proof}

Before determining the positive cone $K_0(C^*(G))_+$, we establish notation and give a theorem, both of which are taken from \cite{effshen}.

\begin{Notation}
	
For $\sigma \in \mathbb{R}_+ \setminus \mathbb{Q}$, let
\[
P_\sigma : = \{ (m, n)\in \mathbb{Z}^2 : \sigma m + n \ge 0 \}
\]
and $(\mathbb{Z}^2, P_\sigma)$ be the Riesz group given by the (total) ordering of $\mathbb{Z}^2$ by $P_\sigma$. We write $\sigma = [c_0, c_1, c_2, \dots]$ for a continued fraction expansion of $\sigma$,
\[
c_0 +
 \cfrac{1}{c_1 +
  \cfrac{1}{c_2 +\cdots}}
\]
where $c_0 \in \IZ$ and $c_i \in \IN$ for $i > 0$. The continued fraction expansion is \textit{simple} if $c_i > 0$ for $i \ge 1$.

\end{Notation}

It is shown in \cite[Lemma 3.1]{effshen} that the continued fraction expansion makes sense if $c_1 > 0$ and if $c_i > 0$ for infinitely many even and infinitely many odd values of $i$. We will need this generality later. The following result is given in this context.

\begin{Theorem} \label{thm effros shen}
	
(\cite[Theorem 3.2]{effshen}) Suppose that $\sigma \in \mathbb{R}_+ \setminus \mathbb{Q}$ has continued fraction expansion $\sigma = [c_0, c_1, \dots ]$ with $c_1 > 0$ and $c_{2j} > 0 , c_{2k+1} > 0$ for infinitely many $j$ and $k$. Then
\[
(\mathbb{Z}^2, P_\sigma) = \lim (\mathbb{Z}^2, \phi_{m, n})
\]
where $\phi_{m, m+1} =
\begin{psmallmatrix}
c_m & 1\\
1 & 0
\end{psmallmatrix}$ and the terms in the sequence $(\mathbb{Z}^2, \phi_{m, n})$ have the standard positive cone.
%If $\sigma = [c_0, c_1, \dots ]$ and $\tau = [d_0, d_1, \dots ]$ are simple continued fraction expansions, then $(\mathbb{Z}^2, P_\sigma)$ and $(\mathbb{Z}^2, P_\tau)$ are order isomorphic if and only if $c_m = d_m$ for sufficiently high $m$.
	
\end{Theorem}

We can now conclude the following:

\begin{Theorem} \label{thm positive cone}

Let $\sigma \in \mathbb{R}^+ \setminus \mathbb{Q}$. Let $\sigma$ have simple continued fraction expansion $[c_0,c_1,\ldots]$. Define integers $k_i \ge 0$ as follows. Let $k_1$ be arbitrary, and for $p \ge 0$,
\[
k_i =
\begin{cases}
 0, &\text{ for } c_1 + c_3 + \cdots + c_{2p-1} + 2 < i < c_1 + c_3 + \cdots + c_{2p+1} + 2 \\
c_{2p}, &\text{ for } i = c_1 + c_3 + \cdots + c_{2p-1} + 2.
\end{cases}
\]
We may indicate this visually as 
\[
(k_i)_{i=1}^\infty
= (k_1, c_0, \underbrace{0,\ldots,0}_{c_1 - 1}, c_2, \underbrace{0,\ldots,0}_{c_3 - 1}, c_4, \underbrace{0,\ldots,0}_{c_5 - 1}, c_6 \ldots).
\]
Let $\Lambda$ be the category of paths as in Definition \ref{def Lambda}. Then
\[
(K_0(C^*(G)), K_0(C^*(G))_+, [1]_0) \cong (\mathbb{Z}^2, P_\sigma, \begin{psmallmatrix} k_1 + 2 \\ k_1 + 1 \end{psmallmatrix}).
\]

\end{Theorem}

\begin{proof}
Note that $\begin{psmallmatrix}
k + 1 & 1\\
k & 1
\end{psmallmatrix}
=
\begin{psmallmatrix}
1 & 1\\
1 & 0
\end{psmallmatrix}
\begin{psmallmatrix}
k & 1\\
1 & 0
\end{psmallmatrix}$. Note also that ``collapsing'' in the sequence $\IZ^2 \xrightarrow{B_1} \IZ^2 \xrightarrow{B_2} \cdots$ yields
\[
\IZ^2 \xrightarrow{B_{j_1} \cdots B_2 B_1} \IZ^2 \xrightarrow{B_{j_2} \cdots B_{j_1 + 2} B_{j_1 + 1}} \IZ^2 \to \cdots
\]
Let $T = \begin{psmallmatrix} 1 & 1 \\ 0 & 1 \end{psmallmatrix}$. When $k_i = 0$ we have $\begin{psmallmatrix} k_i + 1 & 1 \\ k_i & 1 \end{psmallmatrix} = T$, so the $p$th string of zeros in $(k_i)$ collapses to $\IZ^2 \xrightarrow{T^{c_{2p-1}-1}} \IZ^2$. Then the portion $\IZ^2 \xrightarrow{B_{c_{2p}}} \IZ^2 \xrightarrow{T^{c_{2p+1}-1}} \IZ^2$ collapses to the single map
\begin{align*}
T^{c_{2p+1}-1} B_{c_{2p}}
&= \begin{psmallmatrix} 1 & c_{2p+1} - 1 \\ 0 & 1 \end{psmallmatrix} \begin{psmallmatrix} c_{2p} + 1 & 1 \\ c_{2p} & 1 \end{psmallmatrix} \\
&= \begin{psmallmatrix} 1 & c_{2p+1} - 1 \\ 0 & 1 \end{psmallmatrix} \begin{psmallmatrix} 1 & 1 \\ 1 & 0 \end{psmallmatrix} \begin{psmallmatrix} c_{2p} & 1 \\ 1 & 0 \end{psmallmatrix} \\
&=  \begin{psmallmatrix} c_{2p+1} & 1 \\ 1 & 0 \end{psmallmatrix} \begin{psmallmatrix} c_{2p} & 1 \\ 1 & 0 \end{psmallmatrix}
\end{align*}
Thus the sequence $\IZ^2 \xrightarrow{B_1} \IZ^2 \xrightarrow{B_2} \IZ^2 \to \cdots$ collapses to
\[
\IZ^2 \xrightarrow{\begin{psmallmatrix} c_0 & 1 \\ 1 & 0 \end{psmallmatrix} }
\IZ^2 \xrightarrow{\begin{psmallmatrix} c_1 & 1 \\ 1 & 0 \end{psmallmatrix}}
\IZ^2 \to \cdots
\]
Now the fact that $(K_0(C^*(G)), K_0(C^*(G))_+) \cong (\mathbb{Z}^2, P_\sigma)$ follows from Theorems \ref{thm K-theory of G} and \ref{thm effros shen}.  To see that $[1]_0 = \binom{k_1+2}{k_1+1}$, first note that
\begin{align*}
Z(v_1) & = Z(\alpha) \cup Z(\beta) \cup \bigl( \cup_{j=1}^{k_1} Z(\gamma_1^{(j)}) \bigr)\\
& = (Z(\alpha)\setminus Z(\alpha\beta)) \sqcup Z(\beta) \sqcup \bigl( \sqcup_{j=1}^{k_1} Z(\gamma_1^{(j)}) \bigr)\\
& = (Z(\alpha)\setminus Z(\alpha\beta)) \sqcup (Z(\beta)\setminus Z(\beta^2)) \sqcup Z(\beta^2)\\
& \sqcup (\sqcup_{j=1}^{k_1} Z(\gamma_1^{(j)})\setminus Z(\gamma_1^{(j)} \beta)) \sqcup (\sqcup_{j=1}^{k_1} Z(\gamma_1^{(j)} \beta)).\\
\end{align*}
Now applying Lemma \ref{lem equivalent projections}, we have
\[
[1]_0 = [\chi_{Z(v_1)}]_0 = (k_1 + 2)[\chi_{Z(\alpha)\setminus Z(\alpha\beta)}]_0 + (k_1 + 1)[\chi_{Z(\beta^2)}]_0. \qedhere
\]
\end{proof}

\section{Invariant measures on $G^{(0)}$}
\label{section invariant measure}

\begin{Definition} \label{def invariant measure}

Suppose $\Gamma$ is an \'etale groupoid and $\mu$ is a measure on $\Gamma^{(0)}$. We say $\mu$ is \textit{invariant} if $\mu(s(E)) = \mu (r(E))$ for every open bisection $E \subseteq \Gamma$.

\end{Definition}

\begin{Remark}

The groupoids studied in this paper are ample as well as \'etale, in that the unit space is totally disconnected. It follows that there is a base for the topology consisting of compact-open bisections. Then every open bisection is a countable disjoint union of compact-open bisections. Therefore in the condition in Definition \ref{def invariant measure} it suffices to consider only compact-open bisections.

\end{Remark}

Let $\sigma$ and $G$ be as in Theorem \ref{thm positive cone}. We will show that there exists a unique invariant Borel probability measure on $G^{(0)}$ (Theorems \ref{thm uniqueness of invariant measure} and \ref{thm existence of invariant measure}). Both theorems will require several lemmas of preparation. We mention that the existence can be established abstractly. (Namely, the ordered $K$-theory of $C^*(G)$ has a state (the same one as for the corresponding continued fraction AF algebra). Since $C^*(G)$ is nuclear (Remark \ref{remark G amenable}), a quasitrace inducing this state must actually be a trace, which must arise from an invariant measure. See \cite[section 6.9]{blackadar}.) However the analysis used in the proof of uniqueness provides the basis for an explicit proof of existence, and we felt it worthwhile to give this proof.

It will be convenient to use the map from invertible integer matrices to fractional linear transformations, and the action of these on the extended real line, to extend the notation for finite (and infinite) continued fractions to include the case where some of the coefficients are 0, and the final coefficient is an extended real number. We identify $\bigl[ \begin{smallmatrix} a & b \\ c & d \end{smallmatrix} \bigr] \in PGL(2,\IZ)$ with the fractional linear transformation $z \mapsto \frac{az + b}{cz + d}$. We have the quotient maps $\pi : GL(2,\IZ) \to PGL(2,\IZ)$, and $\nu : \IC^2 \setminus \{(0,0)\} \to \IC P^1 \equiv \widetilde \IC$ by $\nu(z_1,z_2) = \frac{z_1}{z_2}$. Then for $T \in GL(2,\IZ)$ and $z \in \widetilde{\IC}$, $\pi(T)(z) = \nu\bigl( T \begin{psmallmatrix}z \\ 1\end{psmallmatrix} \bigr)$. Then we note that $\pi\begin{psmallmatrix} a & 1 \\ 1 & 0 \end{psmallmatrix} (z) = \frac{az + 1}{z} = a + \frac{1}{z}$. For $a_0$, $\ldots$, $a_n \in \IZ$ nonegative, and $\alpha \in \widetilde{\IR}$, we write $[a_0,a_1, \ldots, a_n,\alpha] := \pi\bigl( \begin{psmallmatrix} a_0 & 1 \\ 1 & 0 \end{psmallmatrix} \ldots \begin{psmallmatrix} a_n & 1 \\ 1 & 0 \end{psmallmatrix} \bigr) (\alpha)$. Thus
\[
[a_0,a_1, \ldots, a_n,\alpha]
= a_0 + \cfrac{1}{a_1 + \cfrac{1}{a_2 + \cdots + \cfrac{1}{a_n + \alpha.}}}
\]
The usual arithmetic in $\widetilde{\IR}$ manages the situations where some of the $a_i$ equal 0. Also we have the usual identity (even allowing zero for coefficients):
\begin{align*}
[a_0, \ldots, a_n]
&= \pi\bigl( \begin{psmallmatrix} a_0 & 1 \\ 1 & 0 \end{psmallmatrix} \cdots \begin{psmallmatrix} a_{n-1} & 1 \\ 1 & 0 \end{psmallmatrix} \bigr) (a_n) \\
&= \pi\bigl( \begin{psmallmatrix} a_0 & 1 \\ 1 & 0 \end{psmallmatrix} \cdots \begin{psmallmatrix} a_{i-1} & 1 \\ 1 & 0 \end{psmallmatrix} \bigr) \bigl( \pi\bigl( \begin{psmallmatrix} a_i & 1 \\ 1 & 0 \end{psmallmatrix} \cdots \begin{psmallmatrix} a_{n-1} & 1 \\ 1 & 0 \end{psmallmatrix} \bigr)(a_n) \bigr) \\
&= \pi\bigl( \begin{psmallmatrix} a_0 & 1 \\ 1 & 0 \end{psmallmatrix} \cdots \begin{psmallmatrix} a_{i-1} & 1 \\ 1 & 0 \end{psmallmatrix} \bigr) \bigl( [a_i, \ldots, a_n] \bigr) \\
&= \bigl[ a_0, \ldots, a_{i-1}, [a_i, \ldots, a_n] \bigr].
\end{align*}

Now we begin the preparation for the proof of uniqueness of invariant measures.

\begin{Lemma} \label{lem unique invariant measure one}

Let $c$, $d \ge 0$, $c$, $d \in \IR$, with $c + d = 1$, and let $k \in \IZ$ with $k \ge 0$. Set $B = \begin{psmallmatrix} k+1 & 1 \\ k & 1 \end{psmallmatrix}$. Suppose that $e$, $f \ge 0$, $e$, $f \in \IR$, are such that $\begin{psmallmatrix} c \\ d \end{psmallmatrix} = B^t \begin{psmallmatrix} e \\ f \end{psmallmatrix}$. Then $c < 1$, $e + f = 1 - c$, and $[0,1,k] \le c \le [0,1,k,1]$.

\end{Lemma}

\begin{proof}
We have
\[
\begin{psmallmatrix} e \\ f \end{psmallmatrix}
= \begin{psmallmatrix} k+1 & k \\ 1 & 1 \end{psmallmatrix}^{-1} \begin{psmallmatrix} c \\ d \end{psmallmatrix}
= \begin{psmallmatrix} 1 & -k \\ -1 & k+1 \end{psmallmatrix}
\begin{psmallmatrix} c \\ 1-c \end{psmallmatrix}
= \begin{psmallmatrix} (k+1)c - k \\ -(k+2)c + k+1 \end{psmallmatrix}.
\]
Then $e + f = \bigl( (k+1)c - k \bigr) + \bigl( -(k+2)c + k+1 \bigr)
= 1 - c$. Moreover,
\begin{align*}
0 &\le e = (k+1)c - k, \\
\noalign{hence}
c &\ge \frac{k}{k+1} = \cfrac{1}{1 + \cfrac{1}{k}} = [0,1,k]. \\
\noalign{Also}
0 &\le f = -(k+2)c + k+1, \\
\noalign{hence}
c &\le \frac{k+1}{k+2} = \cfrac{1}{1 + \cfrac{1}{k + \cfrac{1}{1}}} = [0,1,k,1].
\end{align*}
Finally, since $\frac{k+1}{k+2} < 1$ we have that $c < 1$.
\end{proof}

\begin{Lemma} \label{lem unique invariant measure two}

In the context of Lemma \ref{lem unique invariant measure one}, let $\alpha$, $\beta \in \IR$ with $0 \le \alpha \le \frac{e}{1-c} \le \beta \le 1$. Then $[0,1,k+\alpha] \le c \le [0,1,k+\beta]$.

\end{Lemma}

\begin{proof}
Using the equation $e = (k+1)c - k$ from the proof of Lemma \ref{lem unique invariant measure one}, the hypothesis gives
\begin{align*}
\alpha \le \frac{(k+1)c - k}{1-c} &= \frac{c}{1-c} - k \le \beta \\
k + \alpha &\le \frac{c}{1-c} \le k + \beta \\
\frac{1}{k + \alpha} & \ge \frac{1}{c} - 1 \ge \frac{1}{k + \beta} \\
1 + \frac{1}{k + \alpha} &\ge \frac{1}{c} \ge 1 + \frac{1}{k + \beta} \\
[0,1,k + \alpha] = \cfrac{1}{1 + \cfrac{1}{k+\alpha}} &\le c \le \cfrac{1}{1 + \cfrac{1}{k + \beta}} = [0,1,k+\beta]. \qedhere
\end{align*}
\end{proof}

\begin{Proposition} \label{prop approximation of measure}

Let $\sigma \in \IR^+ \setminus \IQ$ have simple continued fraction expansion $[c_0,c_1,\ldots]$ as in Theorem \ref{thm positive cone}, and let $(k_i)_{i=1}^\infty$ and $\Lambda$ be as in that theorem. Suppose that $\mu$ is an invariant Borel probability measure on $G^{(0)}$. Let $a_0 = \mu(Z(\beta)^c)$. Then for each $n \ge 1$,
\[
[0,1,k_1,1,k_2,\ldots,1,k_n] \le a_0 \le [0,1,k_1,1,k_2,\ldots,1,k_n,1].
\]

\end{Proposition}

\begin{proof}
Recall from the proof of Theorem \ref{thm K-theory of G} that $H_i$ denotes the summand $\IZ^2$ of $K_0(C^*(G_i))$, with generators $p_i = [\chi_{Z(\alpha^i) \setminus Z(\alpha^i \beta)}]$ and $q_i = [\chi_{Z(\beta^{i+1})}]$. Let $a_i = \mu_*(p_i) = \mu(Z(\alpha^i) \setminus Z(\alpha^i \beta))$ and $b_i = \mu_*(q_i) = \mu(Z(\beta^{i+1}))$. Then $\begin{psmallmatrix} a_i & b_i \end{psmallmatrix} \in \IR^2 = Hom(H_i,\IR)$ represents $\mu_*|_{H_i}$. It follows that $\begin{psmallmatrix} a_i & b_i \end{psmallmatrix} = \begin{psmallmatrix} a_{i+1} & b_{i+1} \end{psmallmatrix} B_i$. We extend this to the case $i = 0$ using the same formulas.

We apply Lemma \ref{lem unique invariant measure one} with $(c,d) = (a_0,b_0)$, $(e,f) = (a_1,b_1)$, $k = k_1$, and $B = B_0$. In particular we have $a_0 < 1$. Let $(a_n^{(1)},b_n^{(1)}) = (1 - a_0)^{-1} (a_n,b_n)$ for $n \ge 1$. Then
\[
\begin{psmallmatrix} a_n^{(1)} & b_n^{(1)} \end{psmallmatrix} = \begin{psmallmatrix} a_{n+1}^{(1)} & b_{n+1}^{(1)} \end{psmallmatrix} B_n,
\]
for $n \ge 1$. Then $a_1^{(1)} + b_1^{(1)} = 1$, and Lemma \ref{lem unique invariant measure one} implies that $a_1^{(1)} < 1$, $a_2^{(1)} + b_2^{(1)} = 1 - a_1^{(1)}$, and $[0,1,k_2] \le a_1^{(1)} \le [0,1,k_2,1]$. Let $(a_n^{(2)},b_n^{(2)}) = (1 - a_1^{(1)})^{-1} (a_n^{(1)},b_n^{(1)})$ for $n \ge 2$. Then
\[
\begin{psmallmatrix} a_n^{(2)} & b_n^{(2)} \end{psmallmatrix} = \begin{psmallmatrix} a_{n+1}^{(2)} & b_{n+1}^{(2)} \end{psmallmatrix} B_n,
\]
for $n \ge 2$, and Lemma \ref{lem unique invariant measure one} implies that $a_2^{(2)} < 1$, $a_3^{(2)} + b_3^{(2)} = 1 - a_2^{(2)}$, and $[0,1,k_3] \le a_2^{(2)} \le [0,1,k_3,1]$. Continuing this process we define 
$a_n^{(m)}$ and $b_n^{(m)}$ for $n \ge m \ge 1$ so that $a_n^{(m)}$, $b_n^{(m)} \ge 0$, $a_n^{(m)} + b_n^{(m)} = 1$, and
\[
\begin{psmallmatrix} a_n^{(m)} & b_n^{(m)} \end{psmallmatrix} = \begin{psmallmatrix} a_{n+1}^{(m)} & b_{n+1}^{(m)} \end{psmallmatrix} B_n,
\]
for $n \ge m$. For each $n$, Lemma \ref{lem unique invariant measure one} implies that $[0,1,k_n] \le a_{n-1}^{(n-1)} \le [0,1,k_n,1]$. Now we apply Lemma \ref{lem unique invariant measure two} repeatedly: since $a_{n-1}^{(n-1)} = \frac{a_{n-1}^{(n-2)}}{1 - a_{n-2}^{(n-2)}}$ plays the role of $\frac{e}{1 - c}$,
\begin{align*}
[0,1,k_{n-1} + [0,1,k_n]] &\le a_{n-2}^{(n-2)} \le [0,1,k_{n-1} + [0,1,k_n,1]], \text{ that is,} \\
[0,1,k_{n-1},1,k_n] &\le a_{n-2}^{(n-2)} \le [0,1,k_{n-1},1,k_n,1]; \\
\noalign{repeating, we get}
[0,1,k_{n-2} + [0,1,k_{n-1},1,k_n]] &\le a_{n-3}^{(n-3)} \le [0,1,k_{n-2} + [0,1,k_{n-1},1,k_n,1]], \text{ that is,} \\
[0,1,k_{n-2},1,k_{n-1},1,k_n] &\le a_{n-3}^{(n-3)} \le [0,1,k_{n-2},1,k_{n-1},1,k_n,1]; \\
&\qquad \boldsymbol{\cdots} \\
[0,1,k_1,1,k_2,\ldots,1,k_n] &\le a_0 \le [0,1,k_1,1,k_2,\ldots,1,k_n,1]. \qedhere
\end{align*}
\end{proof}

\begin{Theorem} \label{thm uniqueness of invariant measure}

Let $\sigma \in \IR^+ \setminus \IQ$ have simple continued fraction expansion $[c_0,c_1,\ldots]$ as in Theorem \ref{thm positive cone}, and let $(k_i)_{i=1}^\infty$ and $\Lambda$ be as in that theorem. If there exists an invariant Borel probability measure on $G^{(0)}$ it is unique.

\end{Theorem}

\begin{proof}
Write
\begin{align*}
[g_0,g_1,\ldots] &= [0,1,k_1,1,k_2, \ldots] \\
&= [0,1,k_1,1,c_0,1,(0,1)^{c_1-1}, c_2,1,(0,1)^{c_3-1}, c_4,1,\ldots].
\end{align*}
Then $g_1 > 0$, and the consecutive terms $c_{2p},1$ imply that $g_{2j} > 0$ and $g_{2k+1} > 0$ for infinitely many $j$ and $k$. By \cite[Lemma 3.1]{effshen}, $\theta = \lim_{i \to \infty} [g_0,\ldots,g_i]$ exists. It follows from Proposition \ref{prop approximation of measure} that $a_0 = \theta$. Since $b_0 = 1 - a_0$, and $( a_i \ b_i ) = ( a_{i+1} \ b_{i+1} ) B_i$, all $a_i$ and $b_i$ are determined by $a_0$. Then the measure is determined on all sets of the form $Z(\mu_1 \cdots \mu_n)$ and $Z(\mu_1 \cdots \mu_n) \setminus Z( \mu_1 \cdots \mu_{n+1})$. By Proposition \ref{prop partition proposition}, these sets form a base for the topology of $G^{(0)}$. Since Borel probability measures on a compact metrizable space are regular, the measure is completely determined.
\end{proof}

Now we begin the preparation for the proof of the existence of an invariant measure.

\begin{Definition} \label{def basic sets}

For $h \ge 0$ let $\CE_h = \{ Z(\nu), Z(\nu) \setminus Z(\nu \lambda) : \nu \in v_1 \Lambda, \lambda \in s(\nu)\Lambda, |\nu| = h, |\lambda| = 1 \}$, and let $\CE = \cup_h \CE_h$.

\end{Definition}

\begin{Lemma} \label{lem existence invariant measure one}

Let $a_i$, $b_i$ be as in the proof of Theorem \ref{thm uniqueness of invariant measure}. Define $\mu : \CE \to \IR$ by $\mu(Z(\lambda_1 \cdots \lambda_h)) = b_{h-1}$ and $\mu(Z(\lambda_1 \cdots \lambda_h) \setminus Z(\lambda_1 \cdots \lambda_{h+1})) = a_h$. Let $\nu \in v_1 \Lambda$ with $|\nu| = h$ and let $\lambda \in s(\nu) \Lambda$ with $|\lambda| = 1$. For $E = Z(\nu)$ or $Z(\nu) \setminus Z(\nu\lambda)$ we have that $\mu(E) = \sum \{ \mu(S) : S \in Q_h,\ S \subseteq E \}$. (Recall the definition of $Q_h$ from Proposition \ref{prop partition proposition}.)

\end{Lemma}

\begin{proof}
We may write $\nu = \eta \alpha^p \beta^q$ where $\eta \in \Phi_h$ and $p+q = h - |\eta|$. Let 
$n = h - |\eta|$ and let $m = |\eta| + 1$, so that $s(\eta) = v_m$. Then $Z(\nu) = \eta Z(v_m \alpha^p \beta^q)$. We will use equation \eqref{rmk refinements one} from Remark \ref{rmk refinements}. For the moment we will write $[c,d] = \{c,c+1,\ldots,d\}$ for $c \le d$ integers. Then we claim that
\begin{equation} \label{eqn lemma existence one}
[p,n] \times [q,n] = \bigsqcup_{\ell=0}^{\min\{p,q\}} \bigl( ([p+\ell,n] \times \{q + \ell\}) \sqcup (\{p + \ell\} \times [q+\ell+1,n]) \bigr).
\end{equation}
To prove this claim, we first check the disjointness. For fixed $\ell$, the elements of $[p+\ell,n] \times \{q+\ell\}$ and those of $\{p+\ell\} \times [q+\ell+1,n]$ differ in their second coordinate. Let $\ell_1 < \ell_2$, and suppose that $(i,j)$ is in the $\ell_1$-term and in the $\ell_2$-term. If $j = q + \ell_1$ then $j < q + \ell_2$, contradicting the fact that $(i,j)$ is in the $\ell_2$-term. Therefore $j \ge q + \ell_1 + 1$. But then $i = p + \ell_1 < p + \ell_2$, again contradicting the fact that $(i,j)$ is in the $\ell_2$-term. Next we verify the equality. It is clear that the righthand side is contained in the left. (We note that if $p \le q$ then the second part of the last term will be empty.) For the reverse containment, let $(i,j) \in [p,n] \times [q,n]$. First suppose that $i - p < j - q$. Put $\ell = i - p$. Note that then $\ell \le n - p = q$, and $\ell  < j - q \le n - q = p$. Thus $\ell \le \min \{p,q\}$. Then $i = p + \ell$, and $j > q + (i - p) = q + \ell$. Therefore $(i,j) \in \{p+\ell\} \times [q + \ell + 1,n]$. Next suppose that $i - p \ge j - q$. Put $\ell = j - q$. Note that $\ell \le n - q = p$, and $\ell \le i - p \le n - p = q$. Thus again we have that $\ell \le \min \{p,q\}$. Then $j = q + \ell$, and $i \ge p + (j - q) = p + \ell$. Therefore $(i,j) \in [p+\ell,n] \times \{q+\ell\}$. This finishes the proof of the claim.

We use the above claim to rewrite \eqref{rmk refinements one}, as the third term is indexed by $[p,n] \times [q,n]$. For definiteness we assume that $p \le q$ (the other possibility is handled in a similar manner). Then we have
\begin{align*}
Z(v_m \alpha^p \beta^q)
&= \bigsqcup_{\ell=0}^p \Bigl( \bigsqcup_{i=p+\ell}^n \bigsqcup_{r=1}^{k_{m+i+q+\ell}} Z(v_m \alpha^i \beta^{q+\ell} \gamma_{m+i+q+\ell}^{(r)}) \\
&\hspace*{1 in} \sqcup \bigl( Z(v_m \alpha^{n+1} \beta^{q+\ell}) \setminus Z(v_m \alpha^{n+1} \beta^{q+\ell+1}) \bigr) \\
&\quad \sqcup \bigsqcup_{j=q+\ell+1}^n \bigsqcup_{r=1}^{k_{m+p+\ell+j}} Z(v_m \alpha^{p+\ell} \beta^j \gamma_{m+p+\ell+j}^{(r)}) \\
&\hspace*{1 in} \sqcup \bigl( Z(v_m \alpha^{p+\ell} \beta^{n+1}) \setminus Z(v_m \alpha^{p+\ell+1} \beta^{n+1}) \bigr) \Bigr) \\
&\sqcup \bigsqcup_{i=2p+1}^n \bigl( Z(v_m \alpha^i \beta^{n+1}) \setminus Z(v_m \alpha^{i+1} \beta^{n+1}) \bigr) \\
&\sqcup Z(v_m \alpha^{n+1} \beta^{n+1}).
\end{align*}
Now we prepend $\eta$ to all of the sets in the above, and apply $\mu$ to each. The lefthand side becomes $\mu(Z(\eta\alpha^p\beta^q)) = b_{h-1}$. To compute the righthand side we first note that
\begin{align*}
m+i+q+\ell &= |\eta| + 1 + i + q + \ell \\
&= h - p - q + 1 + i + q + \ell \\
&= h + i - p + \ell + 1; \\
|\eta \alpha^i \beta^{q+\ell} \gamma_{m+i+q+\ell}^{(r)}| &= m + i + q + \ell \\
&= h + i - p + \ell + 1; \\
|\eta \alpha^{n+1} \beta^{q+\ell}| &= |\eta| + n+1 + q + \ell \\
&= h + q + \ell + 1; \\
m+p+\ell+j &= |\eta| + 1 + p + \ell + j \\
&= h - p - q + 1 + p + \ell + j \\
&= h + \ell + j - q + 1; \\
|\eta \alpha^{p+\ell} \beta^j \gamma_{m+p+\ell+j}^{(r)}| &= m + p + \ell + j \\
&= h + \ell + j - q + 1 \\
|\eta \alpha^{p+\ell} \beta^{n+1}| &= |\eta| + p + \ell + n+1 \\
&= h + p + \ell + 1; \\
|\eta \alpha^i \beta^{n+1}| &= |\eta| + i + n + 1 \\
&= h + i + 1; \\
|\eta \alpha^{n+1} \beta^{n+1}| &= |\eta| + 2n + 2 \\
&= h + n + 2.
\end{align*}

Then the righthand side becomes
\begin{align}
\sum_{\ell=0}^p \Bigl( \sum_{i=p+\ell}^n & k_{h+i-p+\ell+1} b_{h+i-p+\ell} + a_{h+q+\ell+1} \label{eqn lemma existence two} \\
&+ \sum_{j=q+\ell+1}^n k_{h+\ell+j-q+1} b_{h+\ell+j-q} + a_{h+p+\ell+1} \Bigr) \notag \\
&\quad + \sum_{i=2p+1}^n a_{h+i+1} + b_{h+n+1}. \notag
\end{align}
Recall that
\begin{align*}
\begin{psmallmatrix} a_i & b_i \end{psmallmatrix}
&= \begin{psmallmatrix} a_{i+1} & b_{i+1} \end{psmallmatrix} B_i
= \begin{psmallmatrix} a_{i+1} & b_{i+1} \end{psmallmatrix} \begin{psmallmatrix} k_{i+1} + 1 & 1 \\ k_{i+1} & 1 \end{psmallmatrix}
= \begin{psmallmatrix} (k_{i+1}+1)a_{i+1} + k_{i+1} b_{i+1} & a_{i+1} + b_{i+1} \end{psmallmatrix} \\
\begin{psmallmatrix} a_{i+1} & b_{i+1} \end{psmallmatrix}
&= \begin{psmallmatrix} a_i & b_i \end{psmallmatrix} B_i^{-1} 
= \begin{psmallmatrix} a_i & b_i \end{psmallmatrix} \begin{psmallmatrix} 1 & -1 \\ -k_{i+1} & k_{i+1} + 1\end{psmallmatrix}
= \begin{psmallmatrix} a_i - k_{i+1} b_i & -a_i + (k_{i+1} + 1) b_i \end{psmallmatrix}.
\end{align*}
From these equations we will need
\begin{align}
b_i &= a_{i+1} + b_{i+1} \label{eqn lemma existence three} \\
k_{i+1}b_i &= a_i - a_{i+1}. \label{eqn lemma existence four}
\end{align}
Using \eqref{eqn lemma existence four} we have
\begin{align*}
\sum_{i=p+\ell}^n k_{h+i-p+\ell+1} b_{h+i-p+\ell} + a_{h+q+\ell+1}
&= \sum_{i=p+\ell}^n (a_{h+i-p+\ell} - a_{h+i-p+\ell+1}) + a_{h+q+\ell+1} \\
&= (a_{h+2\ell} - a_{h+n-p+\ell+1}) +a_{h+q+\ell+1} \\
&= a_{h+2\ell}, \text{ since } n = p + q, \\
\sum_{j=q+\ell+1}^n k_{h+\ell+j-q+1} b_{h+\ell+j-q} + a_{h+p+\ell+1}
&= \sum_{j=q+\ell+1}^n (a_{h+\ell+j-q} - a_{h+\ell+j-q+1}) + a_{h+p+\ell+1} \\
&= (a_{h+2\ell+1} - a_{h+\ell+n-q+1}) + a_{h+p+\ell+1} \\
&= a_{h+2\ell+1}. \\
\noalign{Using \eqref{eqn lemma existence three} we have that the sum in \eqref{eqn lemma existence two} equals}
\sum_{\ell=0}^p (a_{h+2\ell} + a_{h+2\ell+1}) + \sum_{i=2p+1}^n a_{h+i+1} + b_{h+n+1}
&= \sum_{i=0}^{n+1} a_{h+i} + b_{h+n+1} \\
&= \sum_{i=0}^n a_{h+i} + b_{h+n} \\
&= \cdots \\
&= a_h + b_h \\
&= b_{h-1}.
\end{align*}
This proves that $\mu(Z(\nu)) = \sum\{\mu(S) : S \in Q_h, S \subseteq Z(\nu)\}$.

Now we consider sets of the form $Z(\nu) \setminus Z(\nu\lambda)$. First we let $\lambda = \alpha$. Prepending $\eta$ to the sets on the righthand side of equation \eqref{rmk refinements two} and applying $\mu$ gives
\begin{align*}
a_{h+p+1} + \sum_{j=q}^n k_{h+j-q+1} b_{h+j-q}
&= a_{h+p+1} + \sum_{j=q}^n (a_{h+j-q} - a_{h+j-q+1}) \\
&= a_{h+p+1} + (a_h - a_{h+p+1}) \\
&= a_h \\
&= \mu(Z(\eta \alpha^p \beta^q) \setminus Z(\eta \alpha^{p+1} \beta^q)).
\end{align*}
The case where $\lambda = \beta$ is similar. Finally, if $\lambda = \gamma_{h+1}^{(r)}$ for some $1 \le r \le k_{h+1}$, the result of prepending $\eta$, then applying $\mu$, to the sets on the righthand side of equation \eqref{rmk refinements four} yields the corresponding result for equation \eqref{rmk refinements one} less the one term $\mu(Z(\eta \alpha^p \beta^q \gamma_{h+1}^{(r)}))$. Using equation \eqref{eqn lemma existence three}, this gives $b_{h-1} - b_h = a_h = \mu(Z(\nu) \setminus Z(\nu\lambda))$, as required.
\end{proof}

\begin{Lemma} \label{lem existence invariant measure two}

Let $0 \le g \le h$. Every set in $\CE_g$ is refined by $Q_h$. For $E \in \CE_g$, $\mu(E) = \sum\{\mu(S) : S \in Q_h, S \subseteq E\}$.

\end{Lemma}

\begin{proof}
When $g = h$ this follows from Proposition \ref{prop partition proposition} and Lemma \ref{lem existence invariant measure one}. We consider the case $h = g+1$. Let $\nu \in v_1 \Lambda$ with $|\nu| = g$. We have
\begin{align*}
Z(\nu)
&= (Z(\nu\alpha) \setminus Z(\nu\alpha\beta)) \sqcup Z(\nu\beta) \sqcup \bigsqcup_{r=1}^{k_{g+1}} Z(\nu \gamma_{g+1}^{(r)}). \\
\noalign{Then clearly we have}
Z(\nu) \setminus Z(\nu\beta)
&= (Z(\nu\alpha) \setminus Z(\nu\alpha\beta)) \sqcup \bigsqcup_{r=1}^{k_{g+1}} Z(\nu \gamma_{g+1}^{(r)}), \\
Z(\nu) \setminus Z(\nu \gamma_{g+1}^{(r)})
&= (Z(\nu\alpha) \setminus Z(\nu\alpha\beta)) \sqcup Z(\nu\beta) \sqcup \bigsqcup_{\substack{s=1 \\ s \not= r}}^{k_{g+1}} Z(\nu \gamma_{g+1}^{(s)}), \\
\noalign{and by symmetry,}
Z(\nu) \setminus Z(\nu\alpha)
&= (Z(\nu\beta) \setminus Z(\nu\alpha\beta)) \sqcup \bigsqcup_{r=1}^{k_{g+1}} Z(\nu \gamma_{g+1}^{(r)}).
\end{align*}
Now we compute the sum of $\mu$-values of the sets in the above decompositions. On the righthand side of the first equation above we get
\begin{align*}
a_{g+1} + b_g + k_{g+1}b_g
&= a_{g+1} + b_g + (a_g - a_{g+1}), \text{ by \eqref{eqn lemma existence four}}, \\
&= a_g + b_g \\
&= b_{g-1}, \text{ by \eqref{eqn lemma existence three}}, \\
&= \mu(Z(\nu)).
\end{align*}
The remaining three equations all lead to the same equation involving $\mu$-values. For the first of these three equations, say, the righthand side gives
\begin{align*}
a_{g+1} + k_{g+1} b_g
&= a_{g+1} + (a_g - a_{g+1}) \\
&= a_g \\
&= \mu(Z(\nu) \setminus Z(\nu\beta)).
\end{align*}
This proves the lemma when $h = g+1$. For the general case, repeating the result just proved shows that each set in $\CE_g$ is refined by sets from $\CE_h$, and the sum of $\mu$-values of the refining sets equals the $\mu$-value of the original set. By Lemma \ref{lem existence invariant measure one} we know that each set in $\CE_h$ is refined by $Q_h$, and the sum of the $\mu$-values of the refining sets equals the $\mu$-value of the original set. Combining these facts finishes the proof of the lemma.
\end{proof}

\begin{Theorem} \label{thm existence of invariant measure}

Let $\sigma \in \IR^+ \setminus \IQ$ have continued fraction expansion $[c_0,c_1,\ldots]$ as in Theorem \ref{thm positive cone}, and let $(k_i)_{i=1}^\infty$ and $\Lambda$ be as in that theorem. There exists an invariant Borel probability measure on $G^{(0)}$.

\end{Theorem}

\begin{proof}
Let $\CA$ be the algebra of subsets of $G^{(0)}$ generated by $\CE$. (Thus $\CA$ is the algebra of all compact-open subsets of $G^{(0)}$.) Let $\mu : \CE \to \IR$ be defined as in Lemma \ref{lem existence invariant measure one}. We claim first that if $F = \sqcup_{i=1}^m E_i$ with $F$, $E_i \in \CE$, then $\mu(F) = \sum_{i=1}^m \mu(E_i)$. To prove this, choose $h$ such that $Q_h$ refines $F$ and all $E_i$; $h$ exists by Lemma \ref{lem existence invariant measure two}. Then $\{S \in Q_h : S \subseteq F \} = \bigsqcup_{i=1}^m \{S \in Q_h : S \subseteq E_i \}$. Then Lemma \ref{lem existence invariant measure two} implies that
\[
\mu(F)
= \sum_{\substack{S \in Q_h \\ S \subseteq F}} \mu(S) 
= \sum_{i=1}^m \sum_{\substack{S \in Q_h \\ S \subseteq E_i}} \mu(S) 
= \sum_{i=1}^m \mu(E_i).
\]
Next, observe that every set in $\CA$ can be written as a (finite) disjoint union of sets from $\CE$. To see this, consider $A \in \CA$. There are $\nu_1$, $\ldots$, $\nu_n \in v_1 \Lambda$ such that $A$ can be constructed from the $Z(\nu_i)$ by intersection, union and difference. By Lemma \ref{lem existence invariant measure two} there is $h$ such that $Z(\nu_1)$, $\ldots$, $Z(\nu_n)$ are all refined by $Q_h$. Then any combination of the $Z(\nu_i)$ using intersection, union and difference, e.g. $A$, will equal the union of a subcollection of $Q_h$. It follows that $A = \bigsqcup \{ S \in Q_h : S \subseteq A \}$. We claim that $\sum\{ \mu(S) : S \in Q_h, S \subseteq A \}$ is independent of $h$, for $h$ large enough that $Q_h$ refines $A$. To see this, let $h_1$ and $h_2$ both be large, and let $h_3 \ge \max \{h_1,h_2\}$. By Lemma \ref{lem existence invariant measure two}, each set in $Q_{h_1} \cup Q_{h_2}$ is refined by $Q_{h_3}$, and the $\mu$-value of the set equals the sum of the $\mu$-values of the refining sets. We have
\[
\sum_{\substack{S \in Q_{h_1} \\ S \subseteq A}} \mu(S)
= \sum_{\substack{S \in Q_{h_1} \\ S \subseteq A}} \sum_{\substack{T \in Q_{h_3} \\ T \subseteq S}} \mu(T) 
= \sum_{\substack{T \in Q_{h_3} \\ T \subseteq A}} \mu(T) 
= \cdots 
= \sum_{\substack{S \in Q_{h_2} \\ S \subseteq A}} \mu(S).
\]
Now for $A \in \CA$ we define $\mu(A) = \sum\{ \mu(S) : S \in Q_h, S \subseteq A \}$ for any $h$ large enough so that $A$ is refined by $Q_h$. The previous claim shows that this is independent of $h$. Now we show that $\mu$ is finitely additive on $\CA$. Let $A_1$, $\ldots$, $A_m \in \CA$ be pairwise disjoint. Choose $h$ large enough that $A_1$, $\ldots$, $A_m$ are all refined by $Q_h$. Then
\[
\mu(\sqcup_{i=1}^m A_i)
= \sum_{\substack{S \in Q_h \\ S \subseteq \sqcup_{i=1}^m A_i}} \mu(S) 
= \sum_{i=1}^m \sum_{\substack{S \in Q_h \\ S \subseteq A_i}} \mu(S) 
= \sum_{i=1}^m \mu(A_i).
\]
Now, since all sets in $\CA$ are compact, $\mu$ is actually countably additive on $\CA$, i.e. $\mu$ is a premeasure on $\CA$. By Caratheodory's theorem, $\mu$ extends (uniquely) to a measure on the Borel sets of $G^{(0)}$.

Finally, we prove that $\mu$ is invariant. Let $\Delta \subseteq G$ be a compact-open bisection. For $g \in \Delta$ there are $\nu_1$, $\nu_2 \in v_1 \Lambda$ with $|\nu_1| = |\nu_2|$, and a compact-open set $E \subseteq s(\nu_1)\Lambda$, such that $\nu_2 E \subseteq s(\Delta)$ and $g \in [\nu_1,\nu_2,E] \subseteq \Delta$. Choose $h$ so that $\nu_2 E$ is refined by $Q_h$. If $S \in Q_h$ with $S \subseteq \nu_2 E$ then $S$ is of one of the forms $Z(\nu_2 \eta)$ or $Z(\nu_2 \eta) \setminus Z(\nu_2 \eta \lambda)$ with $|\lambda| = 1$. Then $\Delta S \Delta^{-1}$ has one of the forms $Z(\nu_1 \eta)$, respectively $Z(\nu_1 \eta) \setminus Z(\nu_1 \eta \lambda)$. In either case, $\mu(\Delta S \Delta^{-1}) = \mu(S)$. Since $\Delta$ is compact we may write $\Delta = \bigsqcup_{i=1}^m [\nu_1^{(i)}, \nu_2^{(i)}, E_i]$ with each term as above. Then we may choose $h$ so that all $\nu_2^{(i)} E_i$ are refined by $Q_h$. It now follows from the above calculation that $\mu(r(\Delta)) = \mu(s(\Delta))$. Therefore $\mu$ is invariant.
\end{proof}

\section{Identifying $C^*(G)$}
\label{section identifying C*(G)}

In their fundamental paper \cite{effshen} introducing the continued fraction AF algebras, Effros and Shen do not actually identify a specific $C^*$-algebra corresponding to an irrational number. In \cite[Section VI.3]{dav}, Davidson does make such an identification.

\begin{Definition} \label{def effros shen algebra} (\cite[Section 10]{hardywright}, \cite[Section 3]{effshen}, \cite[Section VI.3]{dav})

Let $\theta \in \IR \setminus \IQ$ have simple continued fraction expansion $\theta = [g_0,g_1,g_2, \ldots]$. Let the $n$th convergent be $\frac{p_n}{q_n}$. (Thus $p_0 = g_0$, $p_1 = g_0 g_1 + 1$, $q_0 = 1$, $q_1 = g_1$, and for $n \ge 2$, $p_n = g_n p_{n-1} + p_{n-2}$ and $q_n = g_n q_{n-1} + q_{n-2}$.) For $n \ge 0$ let $A_n = M_{q_n} \oplus M_{q_{n-1}}$ (where we let $A_0 = M_{q_0} = \IC \oplus 0$), and for $n \ge 1$ we include $A_{n-1} \hookrightarrow A_{n}$ with partial multiplicities given by $\begin{psmallmatrix} g_n & 1 \\ 1 & 0 \end{psmallmatrix}$. The \textit{Effros Shen algebra} $A_\theta$ is the AF algebra $\overline{\cup_{n \ge 1} A_n}$.

\end{Definition}

It is proved in \cite[VI.3]{dav} that $(K_0(A_\theta), K_0(A_\theta)_+,[1]_0) \cong (\IZ \theta + \IZ, P_\theta, \begin{psmallmatrix} 0 \\ 1 \end{psmallmatrix})$. Moreover, if $\tau$ is the unique tracial state on $A_\theta$, then $\tau_*(K_0(A_\theta)_+) = (\IZ \theta + \IZ)_+ = \begin{psmallmatrix} \theta & 1 \end{psmallmatrix} P_\theta$.

We remark that the above definition of $A_\theta$ does not involve the first quotient $g_0$ of $\theta$. Thus we should think of $A_\theta$ as depending only on the fractional part of $\theta$.

In this section we will prove the following theorem.

\begin{Theorem} \label{main theorem}

Let $\sigma$ be a positive irrational number, with simple continued fraction expansion $[c_0,c_1,\ldots]$. Let $(k_i)_{i \ge 1}$ be the sequence of nonnegative integers defined in Theorem \ref{thm positive cone}. Let $\Lambda$ be the category of paths associated to $(k_i)$ as in Definition \ref{def Lambda}, and let $G$ be the groupoid constructed from $\Lambda$ as in Definition \ref{def the groupoid G}. Let $\theta$ have the (non-simple) continued fraction expansion $[0,1,k_1,1,k_2,1,\ldots]$. Then $C^*(G)$ is isomorphic to the Effros Shen algebra $A_\theta$.

\end{Theorem}

\begin{Remark} \label{rmk sigma theta}

The number $\theta$ is the same as in the proof of Theorem \ref{thm uniqueness of invariant measure}. As observed at the beginning of the proof of that theorem,
\[
[0,1,k_1,1,k_2,\ldots] = [0,1,k_1,1,c_0,1,(0,1)^{c_1 - 1}, c_2, 1, (0,1)^{c_3 - 1}, \ldots].
\]
Using the notation introduced before Lemma \ref{lem unique invariant measure one} we have
\begin{align*}
[1,(0,1)^{a-1}, x]
&= \pi\bigl( \begin{psmallmatrix} 1 & 1 \\ 1 & 0 \end{psmallmatrix} \bigl( \begin{psmallmatrix} 0 & 1 \\ 1 & 0 \end{psmallmatrix} \begin{psmallmatrix} 1 & 1 \\ 1 & 0 \end{psmallmatrix} \bigr)^{a-1} \bigr)(x) \\
&= \pi\bigl( \begin{psmallmatrix}1 & 1 \\ 1 & 0 \end{psmallmatrix} \begin{psmallmatrix}1 & 0 \\ 1 & 1\end{psmallmatrix}^{a-1} \bigr)(x) \\
&= \pi\bigl( \begin{psmallmatrix}1 & 1 \\ 1 & 0\end{psmallmatrix} \begin{psmallmatrix}1 & 0 \\ a-1 & 1\end{psmallmatrix} \bigr)(x) \\
&= \pi\bigl( \begin{psmallmatrix}a & 1 \\ 1 & 0 \end{psmallmatrix} \bigr)(x) \\
&= [a,x].
\end{align*}
Thus $\theta = [0,1,k_1,1,c_0,c_1,c_2,\ldots] = [0,1,k_1,1,\sigma]$. In particular, $C^*(G)$ is Morita equivalent to $A_\sigma$.The reason is that the numbers of the continued fraction expansion determine the multiplicities of the edges in the Bratteli diagram of the corresponding Effros Shen algebra. Since the continued fraction expansion of $\sigma$ is a tail of that of $\theta$, the Bratteli diagram for $A_\sigma$ is a tail of that for $A_\theta$, except for the dimensions of the full matrix algebras at the vertices. It follows that tensoring with the algebra of compact operators renders the diagram for $A_\sigma$ identical to the tail of the diagram for $A_\theta$.

\end{Remark}

Broadly speaking, the proof of Theorem \ref{main theorem} is given in two steps: first showing that $C^*(G)$ is classified by its Elliott invariant, $Ell(C^*(G))$, and then by completing the computation of $Ell(C^*(G))$. To show that it is classifiable we must show that it is infinite dimensional, separable, unital, simple, nuclear, satisfies the UCT, and has finite nuclear dimension.

\begin{Lemma} \label{lem nuclear dimension}

$C^*(G)$ has nuclear dimension at most three.

\end{Lemma}

\begin{proof}

Recall the short exact sequences \ref{equation first exact sequence} and \ref{equation second exact sequence}:
$$ 0 \longrightarrow C^*(G_i|_{U_i}) \longrightarrow C^*(G_i) \longrightarrow C^*(G_i|_{F_i}) \longrightarrow 0 $$
$$ 0 \longrightarrow C^*(G_i|_{F_i^0}) \overset{\iota}{\longrightarrow} C^*(G_i|_{F_i}) \overset{\pi}{\longrightarrow} C^*(G_i|_{F_i^\infty}) \longrightarrow 0.$$
By Corollary \ref{cor K-theory of F i 0}, $C^*(G_i|_{F_i^0})$ is AF and therefore
$$\dim_{nuc}C^*(G_i|_{F_i^0}) = 0$$
by \cite[Remark 2.2(iii)]{winzac}. Corollary \ref{cor K-theory of F i infty} shows $C^*(G_i|_{F_i^\infty}) \cong  M_{\Phi_i \times \Phi_i} \otimes C(\IT)$ so that
$$\dim_{nuc}C^*(G_i|_{F_i^\infty}) = 1$$
by \cite[Proposition 2.4 and Corollary 2.8(i)]{winzac}. Then by \cite[Proposition 2.9]{winzac},
$$\dim_{nuc}C^*(G_i|_{F_i}) \le 2.$$
As in the proof of Proposition \ref{prop K-theory of G i U i}, $C^*(G_i|_{U_i})$ is AF and hence
$$\dim_{nuc}C^*(G_i|_{U_i}) = 0,$$
again by \cite[Remark 2.2(iii)]{winzac}, so that
$$\dim_{nuc} C^*(G_i) \le 3$$
by another application of \cite[Proposition 2.9]{winzac}. Since the above holds for all $i > 0$, \cite[Proposition 2.3(iii)]{winzac} implies
$$\dim_{nuc}C^*(G) \le 3. \qedhere $$
\end{proof}

To show that $C^*(G)$ is simple, we will show that $G$ has the following two properties:

\begin{Definition}

Let $\Gamma$ be a groupoid. We say that $\Gamma$ is \textit{topologically free} if the set $\{x \in \Gamma^{(0)} : x \Gamma x = \{x\}\}$ is dense in $\Gamma^{(0)}$. We say that $\Gamma$ is \textit{minimal} if for every $x \in \Gamma^{(0)}$, the orbit of $x$ is dense in $\Gamma^{(0)}$. We note that if $\Gamma$ is minimal, and if there exists a point in $\Gamma^{(0)}$ with trivial isotropy, then $\Gamma$ is topologically free.

\end{Definition}

\begin{Lemma} \label{lem G topologically free and minimal}

$G$ is topologically free and minimal.

\end{Lemma}

\begin{proof}
We first show that $G$ is minimal. Let $x \in X$ and let $U$ be a nonempty open subset of $X$. We will find $g \in G$ such that $s(g) = x$ and $r(g) \in U$. By Proposition \ref{prop partition proposition} we may assume that $U$ has one of the forms $Z(\mu)$, $Z(\mu) \setminus Z(\mu\alpha)$, $Z(\mu) \setminus Z(\mu\beta)$. We give the proof for the third of these; the proofs of the other two situations are similar. So we assume that $U = Z(\mu) \setminus Z(\mu\beta)$. Let $i \ge |\mu|$ be such that $k_i > 0$. Let $g = [\mu \alpha^{i-|\mu|} \gamma_i^{(1)}, x_1 x_2 \cdots x_i, x_{i+1} x_{i+2} \cdots ]$. Then $s(g) = x$ and $r(g) = \mu \alpha^{i-|\mu|} \gamma_i^{(1)} x_{i+1} x_{i+2} \cdots \in Z(\mu) \setminus Z(\mu\beta)$.

Since $G$ is minimal, to prove that $G$ is topologically free it suffices to exhibit a point of $x$ with trivial isotropy. Note that if $x$ has nontrivial isotropy in $G$ then for some $i$, $x$ has nontrivial isotropy in $G_i$. It then follows from Lemmas \ref{lemma structure of E ell}, \ref{lem F i infty} and \ref{lem F i 0} that only the points of $\cup_i F_i^\infty = \{\mu \alpha^\infty \beta^\infty : \mu \in \Lambda \}$ have nontrivial isotropy. Therefore all other points of $X$ have trivial isotropy.
\end{proof}

\begin{Corollary} \label{cor G simple}

$C^*(G)$ is a simple $C^*$-algebra.

\end{Corollary}

\begin{proof}
This follows from Lemma \ref{lem G topologically free and minimal} and \cite[Corollary 4.6]{ren2}.
\end{proof}

\begin{Lemma} \label{lem G separable and unital}

$C^*(G)$ is separable and unital.

\end{Lemma}

\begin{proof}
Separability follows from the fact that $\Lambda$ is countable. $C^*(G)$ is unital since $G^{(0)}$ is compact.
\end{proof}

\begin{Corollary} \label{cor classifiable}

$C^*(G)$ is classifiable.

\end{Corollary}

\begin{proof}
By Lemmas \ref{lem nuclear dimension} and \ref{lem G separable and unital}, Corollary \ref{cor G simple} and Remark \ref{remark G amenable}, we know that $C^*(G)$ is separable, unital, simple, nuclear, satisfies the UCT, and has finite nuclear dimension. It is clearly infinite dimensional. By \cite[Corollary D]{tww}, $C^*(G)$ is classified by its Elliott invariant.
\end{proof}

\begin{Theorem} \label{thm unique trace}

There exists a unique trace on $C^*(G)$.

\end{Theorem}

\begin{proof}

By Theorems \ref{thm uniqueness of invariant measure} and \ref{thm existence of invariant measure} there exists a unique invariant Borel probability measure, $\mu$, on $G^{(0)}$. Let $E:C^*(G) \rightarrow C_0(G^{(0)})$ be the canonical conditional expectation. Then $\mu$ gives rise to a state $\phi$ on $C^*(G)$ by taking
\[
\phi(f) = \int_{G^{(0)}} f \circ E \, d\mu
\]
for $f\in C_c(G)$. Since $\mu$ is invariant, $\phi$ is tracial.

To see that $\phi$ is the unique tracial state on $C^*(G)$, we show that the points of $G^{(0)}$ with non-trivial isotropy form a set of measure zero, and apply \cite[Corollary 1.2]{nesh}. Suppose $x, y \in G^{(0)}$ have non-trivial isotropy. As in the proof of Lemma \ref{lem G topologically free and minimal}, $x = \sigma \alpha^\infty \beta^\infty$ and $y = \tau \alpha^\infty \beta^\infty$ for some $\sigma, \tau \in \Lambda$. Supposing, say, $|\tau| \ge |\sigma|$, the element $[\sigma \alpha^{|\tau|-|\sigma|}, \tau, \alpha^\infty \beta^\infty]$ has range $x$ and source $y$. Thus the points of $G^{(0)}$ having nontrivial isotropy form a single countably infinite orbit. By invariance $\mu(\{x\}) = \mu(\{y\})$, and since $\mu$ is a finite measure it must be that $\mu(\{x\}) = 0$. \qedhere

\end{proof}

\begin{proof} \textit{(of Theorem \ref{main theorem})}
By Corollary \ref{cor classifiable} we know that $C^*(G)$ is in the class of $C^*$-algebras classified by \cite[Corollary D]{tww}. By Theorem \ref{thm positive cone}, $C^*(G)$ has ordered $K$-theory with position of the unit $(K_0(C^*(G)), K_0(C^*(G))_+,[1]_0) \cong (\IZ^2, P_\sigma, \begin{psmallmatrix} k_1 + 2 \\ k_1 + 1 \end{psmallmatrix})$.

Recall from the proof of Theorem \ref{thm K-theory of G} the maps $\eta_{i,i'} : C^*(G_i) \to C^*(G_{i'})$ and $\eta_i : C^*(G_i) \to C^*(G)$. We thus have the diagram
\[
\begin{tikzpicture}[xscale=2.5,yscale=1.5]

\node (m10) at (-1,0) [rectangle] {$\IZ^2$};
\node (00) at (0,0) [rectangle] {$\IZ^2$};
\node (10) at (1,0) [rectangle] {$\IZ^2$};
\node (20) at (2,0) [rectangle] {$\IZ^2$};
\node (30) at (3,0) [rectangle] {$\cdots$};
\node (40) at (3.75,0) [rectangle] {$\IZ^2$,};
\node (11) at (1,1) [rectangle] {$K_0(G_1)$};
\node (21) at (2,1) [rectangle] {$K_0(G_2)$};
\node (31) at (3,1) [rectangle] {$\cdots$};
\node (41) at (3.75,1) [rectangle] {$K_0(G)$};

\draw[-latex,thick] (m10) -- (00) node[pos=0.5,inner sep=0.5pt,above=1pt] {$T$};
\draw[-latex,thick] (00) -- (10) node[pos=0.5,inner sep=0.5pt,above=1pt] {$B_0$};
\draw[-latex,thick] (10) -- (20) node[pos=0.5,inner sep=0.5pt, above=1pt] {$B_1$};
\draw[-latex,thick] (20) -- (30) node[pos=0.5,inner sep=0.5pt, above=1pt] {$B_2$};
\draw[-latex,thick] (30) -- (40);
\draw[-latex,thick] (11) -- (21) node[pos=0.5,inner sep=0.5pt, below=1pt] {$\eta_{1,2}$};
\draw[-latex,thick] (21) -- (31) node[pos=0.5,inner sep=0.5pt, below=1pt] {$\eta_{2,3}$};
\draw[-latex,thick] (31) -- (41);
\draw[right hook-latex,thick] (10) -- (11);
\draw[right hook-latex,thick] (20) -- (21);
\draw[-latex,thick] (40) -- (41) node[pos=0.5,inner sep=0.5pt, left=1pt] {$\cong$};

\end{tikzpicture}
\]
where $T$ is the composition $\IZ^2 \xrightarrow{\begin{psmallmatrix} 0&1 \\ 1 & 0 \end{psmallmatrix} } \IZ^2 \xrightarrow{ \begin{psmallmatrix} 1&1 \\ 1 & 0 \end{psmallmatrix} } \IZ^2$ (and we have written $K_0(G_i)$ in place of $K_0(C^*(G_i))$). In the diagram above we have extended the sequence back two steps. Thus the lower sequence in the diagram above can be expressed as
\[
\IZ^2 \xrightarrow{\begin{psmallmatrix} 0 1 \\ 1 0 \end{psmallmatrix}}
\IZ^2 \xrightarrow{\begin{psmallmatrix} 1 1 \\ 1 0 \end{psmallmatrix}}
\IZ^2 \xrightarrow{\begin{psmallmatrix} k_1 1 \\ 1 0 \end{psmallmatrix}}
\IZ^2 \xrightarrow{\begin{psmallmatrix} 1 1 \\ 1 0 \end{psmallmatrix}}
\IZ^2 \xrightarrow{\begin{psmallmatrix} k_2 1 \\ 1 0 \end{psmallmatrix}}
\cdots.
\]
We claim that $B_0T$ is an isomorphism from $(\IZ^2, P_\theta, \begin{psmallmatrix} 0 \\ 1 \end{psmallmatrix})$ to $(\IZ^2, P_\sigma, \begin{psmallmatrix} k_1 + 2 \\ k_1 + 1 \end{psmallmatrix})$. Since $T, B_0 \in GL(2,\IZ)$, $B_0T$ is an isomorphism on $\IZ^2$, and an easy calculation shows that $B_0T \begin{psmallmatrix} 0 \\ 1 \end{psmallmatrix} = \begin{psmallmatrix} k_1 + 2 \\ k_1 + 1 \end{psmallmatrix}$. To see that $B_0 T(P_\theta) = P_\sigma$, we first note the following fact. Let $s = [d_0,d_1,t]$, where $d_0 \in \IZ$, $d_1 \in \IN$, and $t \in (0,\infty) \setminus \IQ$. Let $s' = [d_1,t]$. Then $\begin{psmallmatrix} d_0 & 1 \\ 1 & 0 \end{psmallmatrix}(P_s) =P_{s'}$. For the proof, we have
\begin{align*}
\begin{psmallmatrix} d_0 & 1 \\ 1 & 0 \end{psmallmatrix} (P_s)
&= \left\{ \begin{psmallmatrix} d_0 & 1 \\ 1 & 0 \end{psmallmatrix} \begin{psmallmatrix} m \\ n \end{psmallmatrix} : \begin{psmallmatrix} s & 1 \end{psmallmatrix}\begin{psmallmatrix} m \\ n \end{psmallmatrix} \ge 0 \right\} \\
&= \left\{ \begin{psmallmatrix} m \\ n \end{psmallmatrix} : \begin{psmallmatrix} s & 1 \end{psmallmatrix} \begin{psmallmatrix} 0 & 1 \\ 1 & -d_0 \end{psmallmatrix} \begin{psmallmatrix} m \\ n \end{psmallmatrix} \ge 0 \right\} \\
&= \left\{ \begin{psmallmatrix} m \\ n \end{psmallmatrix} : \begin{psmallmatrix} 1 & s - d_0 \end{psmallmatrix} \begin{psmallmatrix} m \\ n \end{psmallmatrix} \ge 0 \right\} \\
&= \left\{ \begin{psmallmatrix} m \\ n \end{psmallmatrix} : \begin{psmallmatrix} \frac{1}{s - d_0} & 1 \end{psmallmatrix} \begin{psmallmatrix} m \\ n \end{psmallmatrix} \ge 0 \right\}, \text{ since $s - d_0 > 0$,} \\
&= P_{\frac{1}{s - d_0}}.
\end{align*}
Since $s - d_0 = [0,d_1,t]$, we have $\frac{1}{s - d_0} = [d_1,t] = s'$. With this fact in hand, we see that
\[
B_0 T(P_\theta)
= \begin{psmallmatrix} 1 & 1 \\ 1 & 0 \end{psmallmatrix}
\begin{psmallmatrix} k_1 & 1 \\ 1 & 0 \end{psmallmatrix}
\begin{psmallmatrix} 1 & 1 \\ 1 & 0 \end{psmallmatrix}
\begin{psmallmatrix} 0 & 1 \\ 1 & 0 \end{psmallmatrix}
(P_{[0,1,k_1,1,\sigma]})
= P_\sigma.
\]
Moreover, $(\IZ^2,P_\theta,\begin{psmallmatrix} 0 \\ 1 \end{psmallmatrix})$ is a simple dimension group, with unique state $\rho$ given by $\rho(\begin{psmallmatrix} 1 \\ 0 \end{psmallmatrix}) = \theta$ and $\rho(\begin{psmallmatrix} 0 \\ 1 \end{psmallmatrix}) = 1$. Since $\tau_*$ does the same, $\tau_* = \rho$. Thus $C^*(G)$ and $A_\theta$ have the same Elliott invariant. By \cite[Theorem D]{tww} it follows that $C^*(G) \cong A_\theta$.
\end{proof}

\begin{Remark}

We have that $\theta = [0,1,k_1,1,k_2,\ldots] = [0,1,k_1,1,\sigma]$ in the notation introduced before Lemma \ref{lem unique invariant measure one}. Since $k_1 \in \IN$, and $\sigma \in (0,\infty) \setminus \IQ$ are arbitrary, a short calculation with a finite continued fraction shows that $\frac{k_1}{k_1 + 1} < \theta < \frac{k_1 + 1}{k_1 + 2}$. Thus $A_\theta$ can be realized in the form $C^*(G)$ for any $\theta \in (0,1) \setminus \IQ$.

\end{Remark}

\begin{Corollary} \label{cor cartan subalgebras}

For each irrational number $\theta \in (0,1)$, the Effros Shen algebra $A_\theta$ contains a Cartan subalgebra that is not conjugate to the standard diagonal subalgebra $D_\theta$ (as in \cite{strvoi}).

\end{Corollary}

\begin{proof}
As shown in the proof of Lemma \ref{lem G topologically free and minimal}, $G$ is not principal but is topologically free (or \textit{topologically principal}). By \cite[Proposition 5.11]{ren3}, $(C^*(G),C(G^{(0)}))$ is a Cartan pair not having the unique extension property, whereas $(A_\theta,D_\theta)$ is a Cartan pair that does have that property. Therefore $C(G^{(0)})$ and $D_\theta$ are nonconjugate Cartan subalgebras in $A_\theta \cong C^*(G)$.
\end{proof}

\begin{Remark}

It is shown in \cite{liren} that every classifiable $C^*$-algebra contains infinitely many pairwise nonisomorphic Cartan subalgebras. The proof shows that all but one of these have spectra with positive dimension, and in fact, the dimension tends to infinity. The examples we give are of a second nonconjugate Cartan subalgebra having zero dimensional spectrum. Both are isomorphic to the continuous functions on a Cantor space, and hence are isomorphic. Thus we prove that the Effros Shen algebras contain two isomorphic but nonconjugate Cartan subalgebras.

\end{Remark}

\section{A locally finite subalgebra of $C^*(G)$}
\label{section subalgebra}

In this section we show that $C^*(G)$ has an isomorphic subalgebra presented as an AF algebra in the usual way. Recall the set $\Phi_i$ from Definition \ref{def Phi_i}. We note that the sets $E_i$ and $F_i$ in the next definition (and the rest of this section) are unrelated to the sets of Definition \ref{def subsets of X}.

\begin{Definition}

For $i \ge 0$ let
\begin{align*}
E_i &= \{\mu \in v_1 \Lambda : |\mu| = i \} \\
F_i &= \{ \eta \beta^q : \eta \in \Phi_i, |\eta| + q = i + 1 \} \\
A_i &= \{ Z(\mu) \setminus Z(\mu\beta) : \mu \in E_i \} \\
B_i &= \{ Z(\mu) : \mu \in F_i \}.
\end{align*}

\end{Definition}

\begin{Lemma} \label{lem partitions for subalgebras}

For each $i \ge 0$,

\begin{enumerate}[1.]

\item $A_i \cup B_i$ is a partition of $X$,

\item $A_{i+1} \cup B_{i+1}$ refines $A_i \cup B_i$. (Here the word \emph{refine} is used in the usual way for a pair of partitions of a set.)

\end{enumerate}

\end{Lemma}

\begin{proof}
Fix $i \ge 0$. We first show that $A_i \cup B_i$ is pairwise disjoint. Let $\mu$, $\nu \in E_i$ be distinct. Write $\mu = \xi \alpha^p \beta^{i - |\xi| - p}$ and $\nu = \eta \alpha^q \beta^{i - |\eta| - q}$, where $\xi$, $\eta \in \Phi_i$. If $|\xi| = i$ or $|\eta| = i$ then it is clear that $Z(\mu) \cap Z(\nu) = \varnothing$. So we may suppose that $|\xi|,|\eta| < i$. Suppose there exists $x \in (Z(\mu) \setminus Z(\mu\beta)) \cap (Z(\nu) \setminus Z(\nu\beta))$. Letting $\xi = \theta_1 \cdots \theta_k$ and $\eta = \phi_1 \cdots \phi_\ell$ in normal form, we have that the normal form of $x$ is
\begin{align*}
x &= \theta_1 \cdots \theta_k \alpha^{p'} \beta^{i-|\xi|-p} \theta_{k+2} \cdots \\
&= \phi_1 \cdots \phi_\ell \alpha^{q'} \beta^{i-|\eta|-q} \phi_{\ell+2} \cdots,
\end{align*}
where $p \le p' \le \infty$ and $q \le q' \le \infty$. Then the uniqueness of the normal form implies that $\xi \alpha^{p'} \beta^{i-|\xi|-p} = \eta \alpha^{q'} \beta^{i-|\eta|-q}$. Then $\xi = \eta$ and $i - |\xi| - p = i - |\eta| - q$, hence $p = q$, hence $\mu = \nu$. Therefore $A_i$ is a pairwise disjoint family.  Now let $\mu$, $\nu \in F_i$. Write $\mu = \theta_1 \cdots \theta_m \beta^p$ and $\nu = \phi_1 \cdots \phi_n \beta^q$ in normal form, with $p$, $q \ge 1$. Suppose that $Z(\mu) \cap Z(\nu) \not= \varnothing$, i.e. that $\mu \Cap \nu$. Without loss of generality let $m \le n$. By Lemma  \ref{lemma:common extension}, $\phi_j = \theta_j$ for $j \le m$, and $\beta^p \Cap \phi_{m+1} \cdots \phi_n \beta^q$. If $m < n$ then we must have $\beta^p \in [\phi_{m+1}]$, and this contradicts the fact that $|\mu| = i = |\nu|$. Thus $m = n$. Then again since $|\mu| = |\nu|$ it follows that $p = q$, hence that $\mu = \nu$. Thus $B_i$ is pairwise disjoint. Finally, let $\mu \in E_i$ and $\nu \in F_i$. Suppose that there exists $x \in (Z(\mu) \setminus Z(\mu\beta)) \cap Z(\nu)$. Write $\nu = \eta \beta^q$ with $\eta \in \Phi_i$ and $q > 0$. Since $\mu \Cap \nu$, and $|\mu| < |\nu|$, we must have $\mu = \eta \alpha^p \beta ^{i - |\eta| - p}$, with $p \le i - |\eta|$. Since $x \not\in Z(\mu\beta)$ we must have $x = \eta \alpha^{p'}\beta^{i - |\eta| - p} y$, where $y \in \partial \Lambda$ either equals $\alpha^\infty$ or begins in $\Lambda_2$, and $p \le p'$. But $i - |\eta| - p \le i - |\eta| = q - 1$, so $x \not\in Z(\nu)$, a contradiction. Thus no such $x$ can exist, and it follows that $A_i \cup B_i$ is a pairwise disjoint family.

Next we show that $X$ equals the union of the family of sets $A_i \cup B_i$. Let $x \in X$. First we consider the case that the normal form of $x$ does not end with an element of $\Lambda_1^\infty$. Write $x = \eta \alpha^p \beta^q y$, where $\eta \in \Phi_i$, $|\eta| + p + q \ge i$, and $y$ begins in $\Lambda_2$. If $|\eta| + q \le i$, then $\xi := \eta \alpha^{i - |\eta| - q} \beta^q \in E_i$, and $x \in Z(\xi) \setminus Z(\xi\beta)$, an element of $A_i$.  Suppose that $q > i - |\eta|$. Then letting $\xi := \eta \beta^{i - |\eta| + 1} \in x$ we have $\xi \in F_i$. Thus $x \in Z(\xi)$, an element of $B_i$. Secondly we consider the case that $x = \eta \alpha^p \beta^q$, where $\eta \in \Phi_i$ and $p + q = \infty$.  If $q \le i - |\eta|$, let $\xi = \eta \alpha^{i - |\eta| - q} \beta^q \in E_i$. Then $x \in Z(\xi) \setminus Z(\xi\beta)$, an element of $A_i$. If $q > i - |\eta|$, let $\xi = \eta \beta^{i - |\eta| + 1} \in F_i$. Then $x \in Z(\xi)$, an element of $B_i$.

Finally we show that $A_{i+1} \cup B_{i+1}$ refines $A_i \cup B_i$. First let $\mu \in E_i$. For $1 \le r \le k_{i+1}$ we have that $\mu \gamma_{i+1}^{(r)} \in E_{i+1}$ and $\mu \gamma_{i+1}^{(r)} \beta \in F_{i+1}$. Then $Z(\mu \gamma_{i+1}^{(r)}) \setminus Z(\mu \gamma_{i+1}^{(r)} \beta) \in A_{i+1}$, $Z(\mu \gamma_{i+1}^{(r)} \beta) \in B_{i+1}$, and
\begin{align*}
Z(\mu) \setminus Z(\mu\beta)
&= (Z(\mu\alpha) \setminus Z(\mu\alpha\beta)) \\
&\quad \cup \bigl( \bigcup_{r=1}^{k_{i+1}} (Z(\mu \gamma_{i+1}^{(r)}) \setminus Z(\mu \gamma_{i+1}^{(r)} \beta)) \bigr) \\
&\qquad \cup \bigl( \bigcup_{r=1}^{k_{i+1}} Z(\mu \gamma_{i+1}^{(r)} \beta) \bigr).
\end{align*}
Thus $A_{i+1} \cup B_{i+1}$ refines all elements of $A_i$. Next let $\mu \in F_i$. Then $\mu\beta \in F_{i+1}$, so $Z(\mu\beta) \in B_{i+1}$. Moreover $\mu \in E_{i+1}$, so $Z(\mu) \setminus Z(\mu\beta) \in A_{i+1}$. Since $Z(\mu) = (Z(\mu) \setminus Z(\mu\beta)) \cup Z(\mu\beta)$, we have that $A_{i+1} \cup B_{i+1}$ refines all elements of $B_i$.
\end{proof}

\begin{Lemma} \label{lem inductive E and F}

For $i \ge 0$,
\begin{align*}
E_{i+1} &= F_i \sqcup \bigl( \bigsqcup_{r=1}^{k_{i+1}} E_i \gamma_{i+1}^{(r)} \bigr) \sqcup E_i \alpha \\
F_{i+1} &= F_i \beta \sqcup \bigl( \bigsqcup_{r=1}^{k_{i+1}} E_i \gamma_{i+1}^{(r)} \beta \bigr).
\end{align*}

\end{Lemma}

\begin{proof}
Both containments of the right side in the left are clear. For the reverse containments, first let $\mu \in E_{i+1}$. Write $\mu = \eta \alpha^p \beta^{i+1 - |\eta| - p}$, where $\eta \in \Phi_{i+1}$ and $p \le i + 1 - |\eta|$. If $|\eta| = i+1$ then the last edge in $\mu$ is in $\Lambda_2$, i.e. $\mu \in \bigsqcup_{r=1}^{k_{i+1}} E_i \gamma_{i+1}^{(r)}$. If $|\eta| \le i$ and $p > 0$ then $\mu \in E_i \alpha$. If $|\eta| \le i$ and $p = 0$ then $\mu \in F_i \beta$, proving the first equality. For the second, let $\mu \in F_{i+1}$. Then $\mu = \eta \beta^q$ with $\eta \in \Phi_{i+1}$ and $q = i+2 - |\eta|$. If $q \ge 2$, then $\mu \in F_i \beta$. If $q = 1$, then $|\eta| = i+1$. Therefore $\eta \in \bigsqcup_{r=1}^{k_{i+1}} E_i \gamma_{i+1}^{(r)}$, and hence $\mu \in \bigsqcup_{r=1}^{k_{i+1}} E_i \gamma_{i+1}^{(r)} \beta$.
\end{proof}

\begin{Theorem} \label{thm locally finite subalgebra}

There is a sequence $C_0 \subseteq C_1 \subseteq \cdots$ of finite dimensional $C^*$-subalgebras of $C^*(G)$ with $C_i \cong M_{|E_i|} \oplus M_{|F_i|}$. $C := \overline{\cup_{i=0}^\infty C_i}$ is isomorphic to $C^*(G)$, and the inclusion $C \hookrightarrow C^*(G)$ induces an isomorphism of Elliott invariants.

\end{Theorem}

\begin{proof}
Fix $i \ge 0$. Lemma \ref{lem equivalent projections} implies that for $\mu$, $\nu \in E_i$, $\chi_{[\nu,\mu,Z(\beta_{i+1})^c]}$ is a partial isometry in $C^*(G)$ with initial and final projections equal to $\chi_{Z(\mu) \setminus Z(\mu\beta)}$ and $\chi_{Z(\nu) \setminus Z(\nu\beta)}$, and that for $\xi$, $\eta \in F_i$, $\chi_{[\eta,\xi,Z(v_{i+2})]}$ is a partial isometry in $C^*(G)$ with initial and final projections $\chi_{Z(\xi)}$ and $\chi_{Z(\eta)}$. It is clear that the span of these partial isometries is a finite dimensional $C^*$-algebra $C_i \subseteq C^*(G)$, and that these partial isometries are matrix units defining an isomorphism of $C_i$ with $M_{|E_i|} \oplus M_{|F_i|}$. By Lemma \ref{lem partitions for subalgebras} it follows that $C_i \subseteq C_{i+1}$. Lemma \ref{lem inductive E and F} implies that for $\xi \in E_i$, $Z(\xi) \setminus Z(\xi\beta)$ is the disjoint union of $k_{i+1} + 1$ sets of the form $Z(\mu) \setminus Z(\mu\beta)$ with $\mu \in E_{i+1}$ and one set of the form $Z(\nu)$ with $\nu \in F_{i+1}$. Similarly, for $\eta \in F_i$, $Z(\eta)$ is the disjoint union of $k_{i+1}$ sets of the form $Z(\mu) \setminus Z(\mu\beta)$ with $\mu \in E_{i+1}$ and one set of the form $Z(\nu)$ with $\nu \in F_{i+1}$. Thus the matrix of multiplicities of the embeddings of the summands of $C_i$ into the summands of $C_{i+1}$ are given by the matrix $B_i$ from Lemma \ref{lem connecting maps in K 0 two}. Therefore $C$ is isomorphic to the Effros Shen algebra $A_\theta$ from Theorem \ref{main theorem}, and thus also to $C^*(G)$. Now it follows from Theorem \ref{K-theory of G i} that the inclusion induces an isomorphism of Elliott invariants.
\end{proof}

\begin{Remark} \label{remark nonrefining}

The partitions defined in Lemma \ref{lem partitions for subalgebras} are simpler than those defined in Lemma \ref{lem partition lemma} and Proposition \ref{prop partition proposition}. However the more complicated partitions are necessary to prove the results of the earlier sections. The reason is that the partitions of Lemma \ref{lem partitions for subalgebras} do not refine all cylinder sets. To see this we consider $Z(\alpha)$. For each $i$, $Z(\beta^{i+1}) \in B_i$. But $Z(\beta^{i+1}) \cap Z(\alpha) = Z(\alpha\beta^{i+1}) \not= \varnothing$, and also $Z(\beta^{i+1}) \setminus Z(\alpha)$ is nonempty, since for example $\beta^\infty$ is in it. Thus none of the partitions $A_i \cup B_i$ refines $Z(\alpha)$. One consequence of this is the following.

\end{Remark}

\begin{Theorem} \label{thm proper subalgebra}
The AF subalgebra $C$ in Theorem \ref{thm locally finite subalgebra} is proper.
\end{Theorem}

\begin{proof}
We will show that $\chi_{Z(\alpha)} \not\in C$. Let diag$\,C_i$ denote the diagonal subalgebra of $C_i$, i.e. the span of $\{\chi_D : D \in A_i \cup B_i \}$. Then the diagonal of $C$, diag$\,C$, is the closure of $\bigcup_i \text{diag}\,C_i$, and the closure occurs in $C(G^{(0)})$, with the uniform norm. Since diag$\,C$ is a masa in $C$, and commutes with $\chi_{Z(\alpha)}$, it suffices to show that $\chi_{Z(\alpha)} \not\in \text{diag}\,C$.

Suppose to the contrary that $\chi_{Z(\alpha)} \in \text{diag}\,C$. Then there are $i$, $D_1,\ldots,D_n \in A_i \cup B_i$, and $c_1,\ldots,c_n \in \IC$ such that $f = \sum_{j=1}^n c_j \chi_{D_j} \in \text{diag}\,C_i$ satisfies $\|\chi_{Z(\alpha)} - f \|_u < \frac{1}{2}$. If $j$ is such that $Z(\alpha) \cap D_j \not= \varnothing$, let $x \in Z(\alpha) \cap D_j$. Since $D_1,\ldots,D_n$ are pairwise disjoint we have $x \not\in D_\ell$ for $\ell \not= j$. Then
\[
\tfrac{1}{2} > |\chi_{Z(\alpha)}(x) - c_j \chi_{D_j}(x)| = |1 - c_j|.
\]
If $j$ is such that $D_j \setminus Z(\alpha) \not= \varnothing$, let $x \in D_j \setminus Z(\alpha)$. Again, $x \not\in D_\ell$ for $\ell \not= j$, and we have
\[
\tfrac{1}{2} > |\chi_{Z(\alpha)}(x) - c_j \chi_{D_j}(x)| = |0 - c_j| = |c_j|.
\]
Thus the same value of $j$ cannot satisfy both properties. It follows that for each $j$ either $D_j \subseteq Z(\alpha)$ or $D_j \cap Z(\alpha) = \varnothing$. Since $A_i \cup B_i$ does not refine $Z(\alpha)$, Remark \ref{remark nonrefining} implies that there exists $x \in Z(\alpha) \setminus \cup \{D_j : D_j \subseteq Z(\alpha)\}$. Then $\chi_{D_j}(x) = 0$ for all $j$, and hence
\[
\| \chi_{Z(\alpha)} - f \|_u
\ge |\chi_{Z(\alpha)}(x) - \sum_{j=1}^n c_j \chi_{D_j}(x)|
= |1 - 0|
= 1,
\]
a contradiction.
\end{proof}

\section{Stability of $C^*(G(\Lambda))$} \label{section stability}

In this section we fix the sequence $(k_i)_{i=1}^\infty$. We noted before Definition \ref{def the groupoid G} that $C^*(G)$ is Morita equivalent to $C^*(\Lambda)$. In this section we investigate $C^*(\Lambda)$. Since $C^*(\Lambda)$ is an AF algebra, we identify it by means of the scale it determines in $K_0(C^*(\Lambda))_+ \linebreak[2] = (\IZ \theta + \IZ)_+$. (The scale is the image in $K_0$ of the projections actually in the algebra \cite[Section 6.1]{blackadar}. We know, for example, that the scale determined by $C^*(G)$ is $(\IZ \theta + \IZ) \cap [0,1]$.) In $C^*(\Lambda)$ the projections $\chi_{[v_n,v_n,Z(v_n)]}$, $n \ge 1$, are pairwise orthogonal, and their partial sums form an approximate identity. Therefore the scale equals $\{t \in (\IZ \theta + \IZ)_+ : t < \sum_{i=1}^n [\chi_{[v_i,v_i,Z(v_i)]}]_0, \text{ for some }n \}$.

\begin{Lemma} \label{lem v_n equivalent to beta^{n-1}}

$\chi_{[v_n,v_n,Z(v_n)]}$ and $\chi_{[v_1,v_1,Z(\beta^{n-1})]}$ are equivalent projections in $C^*(\Lambda)$.

\end{Lemma}

\begin{proof}
The partial isometry $\chi_{[v_n,\beta^{n-1},Z(v_n)]}$ implements the claimed equivalence.
\end{proof}

Recall from the proof of Proposition \ref{prop approximation of measure} that we let $b_i = \mu(Z(\beta^{i+1})) \in \IZ \theta + \IZ$. Thus Lemma \ref{lem v_n equivalent to beta^{n-1}} implies that $[\chi_{[v_n,v_n,Z(v_n)]}]_0 = b_{n-2}$ for $n \ge 2$, and $[\chi_{[v_1,v_1, Z(v_1)]}]_0 = 1$.

\begin{Theorem} \label{thm scale}

For $i \ge 0$ let $\theta_i = [0,1,k_{i+1},1,k_{i+2},1,\ldots]$. The scale of $K_0(C^*(\Lambda))$ equals the set $\{t \in \IZ \theta + \IZ : 0 \le t < 1 + \sum_{j=1}^\infty \prod_{i=0}^j (1 - \theta_i)\}$.

\end{Theorem}

\begin{proof}
From the above remarks, and Lemma \ref{lem v_n equivalent to beta^{n-1}}, the theorem will follow if we show that $b_n = \prod_{i=0}^n (1 - \theta_i)$. Recalling the proof of Proposition \ref{prop approximation of measure}, we have $a_0 = \theta_0$, $b_0 = 1 - \theta_0$, and also $\begin{psmallmatrix} a_i & b_i \end{psmallmatrix} = \begin{psmallmatrix} a_{i+1} & b_{i+1} \end{psmallmatrix} B_i$, $i \ge 0$, where $B_i = \begin{psmallmatrix} k_{i+1} + 1 & 1 \\ k_{i+1} & 1 \end{psmallmatrix}$. We defined $a_i^{(1)}$ and $b_i^{(1)}$, $i \ge 1$, by
\begin{align*}
\begin{psmallmatrix} a_1 & b_1 \end{psmallmatrix}
&= (1 - a_0) \begin{psmallmatrix} a_1^{(1)} & b_1^{(1)} \end{psmallmatrix} \\
\begin{psmallmatrix} a_i^{(1)} & b_i^{(1)} \end{psmallmatrix} &= \begin{psmallmatrix} a_{i+1}^{(1)} & b_{i+1}^{(1)} \end{psmallmatrix} B_i.
\end{align*}
It follows from the proof of Proposition \ref{prop approximation of measure} that $a_1^{(1)} = \theta_1$ and $b_1^{(1)} = 1 - \theta_1$. Therefore $b_1 = (1 - \theta_0)(1 - \theta_1)$. Inductively we have that $b_n = \prod_{i=0}^n (1 - \theta_i)$.
\end{proof}

\begin{Example}

We consider the simplest example: $k_i = 1$ for $i \ge 1$. Then $\theta_i = [0,1,1,1,\ldots] = \tau - 1 = \tau^{-1}$ for all $i \ge 0$, where $\tau = \frac{1}{2} (\sqrt{5} + 1)$ is the golden ratio. Then $1 - \theta_i = 1 - \tau^{-1} = 2 - \tau = \tau^{-2}$. Then the scale of $K_0(C^*(\Lambda))$ is defined (as in Theorem \ref{thm scale}) by
\begin{align*}
1 + \sum_{j=1}^\infty (\tau^{-2})^{j+1}
&= 1 + \sum_{j=2}^\infty \tau^{-2j} \\
&= 1 + \tau^{-4} \frac{1}{1 - \tau^{-2}} \\
&= 1 + \frac{1}{\tau^2} \cdot \frac{1}{\tau^2 - 1} \\
&= 1 + (2 - \tau)(\tau - 1) \\
&= 2\tau - 2 \\
&= \frac{2}{\tau}.
\end{align*}

\end{Example}

We now wish to give reasonable bounds on $1 + \sum_{j=1}^\infty \prod_{i=0}^j (1 - \theta_i)$. Let $(k_{p_i})_{i=1}^\infty$ be the subsequence of $(k_j)$ consisting of those terms that are nonzero. Thus $(k_j)_{j=1}^\infty = (0^{p_1 - 1}, k_{p_1}, 0^{p_2 - p_1 - 1}, k_{p_2}, \ldots)$. We let $q_i = p_i - p_{i-1}$ for $i \ge 1$ (with $p_0 := 0$).

\begin{Lemma} \label{lem first estimates}

For $\ell \ge 0$,
\begin{gather*}
\tfrac{1}{4} (q_{\ell+1} - 1)
< \sum_{j=0}^{q_{\ell+1}-1} \prod_{i=0}^j (1 - \theta_{p_\ell+i})
< \tfrac{2}{3} q_{\ell+1} \\
\frac{e^{-1}}{q_{\ell+1}} e^{-\frac{1}{2} (k_{p_{\ell+1}} + 3)}
< \prod_{i=0}^{q_{\ell+1}-1} (1 - \theta_{p_\ell+i})
< \frac{2}{q_{\ell+1} + 2}.
\end{gather*}

\end{Lemma}

\begin{proof}
Recalling the remarks before Lemma \ref{lem unique invariant measure one}, we have for $p > 1$ and $y > 0$,
\begin{align*}
[(0,1)^p,y]
&= \pi \Bigl( \bigl( \begin{psmallmatrix} 0 & 1 \\ 1 & 0 \end{psmallmatrix} \begin{psmallmatrix} 1 & 1 \\ 1 & 0 \end{psmallmatrix} \bigr)^p \Bigr) (y) 
= \pi \Bigl( \begin{psmallmatrix} 1 & 0 \\ 1 & 1 \end{psmallmatrix}^p \Bigr) (y) \\
&= \pi \Bigl( \begin{psmallmatrix} 1 & 0 \\ p & 1 \end{psmallmatrix} \Bigr) (y) 
= \frac{y}{py+1} 
= \frac{1}{p + \frac{1}{y}} 
= [0,p,y].
\end{align*}
It follows that for $i > 1$, $k \ge 1$, and $x > 0$, we have that
\[
\frac{1}{i+1} < [(0,1)^i,k,1,x] < \frac{1}{i}.
\]
Note also that
\[
[0,1,k,1,x] < [0,1,k,1] = \frac{k+1}{k+2},
\]
and therefore that $\frac{1}{2} < [0,1,k,1,x] < \frac{k+1}{k+2}$. Recalling the definition of $\theta_i$ from Theorem \ref{thm scale} we note that $\theta_0 = [0,1,(0,1)^{q_1-1}, k_{p_1}, 1, (0,1)^{q_2-1}, k_{p_2}, 1, \ldots]$.
Then we have
\begin{align*}
\theta_{p_\ell + i}
&= [0,1,(0,1)^{q_{\ell+1} - i - 1}, k_{p_{\ell+1}}, 1, \ldots] \\
&= [(0,1)^{q_{\ell+1}-i}, k_{p_{\ell+1}}, 1, \ldots]
\in \Bigl( \frac{1}{q_{\ell+1} - i + 1}, \frac{1}{q_{\ell+1} - i} \Bigr), \text{ for } 0 \le i < q_{\ell+1} - 1, \\
\theta_{p_{\ell+1} - 1}
&= [0,1,k_{p_{\ell+1}},1,\ldots]
\in \Bigl( \frac{1}{2}, \frac{k_{p_{\ell+1}} + 1}{k_{p_{\ell+1}} + 2} \Bigr).
\end{align*}
We summarize this as
\begin{equation} \label{eqn theta_i estimates}
\frac{1}{q_{\ell+1} - i + 1}
< \theta_{p_\ell + i}
< \begin{cases}\displaystyle
\frac{1}{q_{\ell+1} - i}, &\text{ if } 0 \le i < q_{\ell+1} - 1 \\ \\ \displaystyle
\frac{k_{p_{\ell+1}} + 1}{k_{p_{\ell+1}} + 2}, &\text{ if } i = q_{\ell+1} - 1.
\end{cases}
\end{equation}
We use the estimates
\begin{equation} \label{eqn sum t^n/n}
t
< \sum_{n=1}^\infty \tfrac{1}{n} t^n
< t + \tfrac{1}{2} \sum_{n=2}^\infty t^n
= t + \tfrac{1}{2} \frac{t^2}{1-t}, \text{ for } 0 < t < 1,
\end{equation}
and for $1 \le m \le n$,
\begin{equation*}
\int_m^{n+1} \frac{1}{u} \, du
< \sum_{i=m}^n \frac{1}{i}
\le \begin{cases}
\int_{m-1}^n \frac{1}{u} \, du, &\text{ if } m \ge 2 \\ \\
1 + \int_1^n \frac{1}{u} \, du, &\text{ if } m = 1,
\end{cases}
\end{equation*}
or equivalently,
\begin{equation} \label{eqn partial harmonic series}
\log \frac{n+1}{m}
< \sum_{i=m}^n \frac{1}{i}
\le \begin{cases}
\log \frac{n}{m-1}, &\text{ if } m \ge 2 \\ \\
1 + \log n, &\text{ if } m = 1.
\end{cases}
\end{equation}
Since $\log(1-t) = -\sum_{n=1}^\infty \frac{1}{n} t^n$ for $|t| < 1$, equation \eqref{eqn sum t^n/n} implies
\begin{equation} \label{eqn log(1-t)}
-t - \tfrac{1}{2} \frac{t^2}{1 - t} < \log(1 - t) < -t.
\end{equation}
Now we establish the upper bounds. For $0 \le j \le q_{\ell + 1} - 1$ we have
\begin{align*}
\sum_{i=0}^j \log(1-\theta_{p_\ell+i})
&< -\sum_{i=0}^j \theta_{p_\ell+i}, \text{ by \eqref{eqn log(1-t)},} \\
&< -\sum_{i=0}^j \frac{1}{q_{\ell + 1} - i + 1}, \text{ by \eqref{eqn theta_i estimates},} \\
&= -\sum_{i = q_{\ell + 1} + 1 - j}^{q_{\ell + 1} + 1} \frac{1}{i} \\
&< -\log \frac{q_{\ell + 1} + 2}{q_{\ell + 1} + 1 - j}, \text{ by \eqref{eqn partial harmonic series},} \\
&= \log \bigl( 1 - \frac{j+1}{q_{\ell + 1} + 2} \bigr), \\
\prod_{i=0}^j (1 - \theta_{p_\ell+i})
&< 1 - \frac{j+1}{q_{\ell + 1} + 2}. \\
\noalign{Then}
\prod_{i=0}^{q_{\ell+1}-1} (1 - \theta_{p_\ell+i})
&< \frac{2}{q_{\ell+1}+2}, \\
\noalign{and}
\sum_{j=0}^{q_{\ell + 1}-1} \prod_{i=0}^j (1 - \theta_{p_\ell+i})
&< q_{\ell+1} - \frac{q_{\ell+1} (q_{\ell+1} + 1)}{2(q_{\ell+1} + 2)} \\
& = q_{\ell+1} \frac{q_{\ell+1} + 3}{2(q_{\ell+1} + 2)} \\
& \le \tfrac{2}{3} q_{\ell+1},
\end{align*}
since $\frac{t+3}{t+2}$ decreases for $t \ge 1$.

For the lower bounds we first note that since $\displaystyle t + \tfrac{1}{2} \frac{t^2}{1 - t} = t + \tfrac{1}{2} \sum_{i=2}^\infty t^i$ is increasing for $t \in (0,1)$, then for $0 \le i < q_{\ell+1} - 1$, \eqref{eqn theta_i estimates} gives
\begin{align}
\theta_{p_\ell+i} + \tfrac{1}{2} \frac{\theta_{p_\ell+i}^2}{1 - \theta_{p_\ell+i}}
&< \frac{1}{q_{\ell+1} - i} + \tfrac{1}{2} \tfrac{1}{(q_{\ell+1} - i)^2 \bigl(1 - \frac{1}{q_{\ell+1} - i} \bigr)} \notag \\
&= \frac{1}{q_{\ell+1} - i} + \tfrac{1}{2} \bigl( \frac{1}{q_{\ell+1} - i -1} - \frac{1}{q_{\ell+1} - i} \bigr) \notag \\
&= \tfrac{1}{2} \bigl( \frac{1}{q_{\ell+1} - i} + \frac{1}{q_{\ell+1} - i - 1} \bigr) \notag \\
&< \frac{1}{q_{\ell+1} - i - 1}, \label{eqn special one}
\end{align}
while
\begin{align}
\theta_{p_{\ell+1} - 1} + \tfrac{1}{2} \frac{\theta_{p_{\ell+1} - 1}^2}{1 - \theta_{p_{\ell+1} - 1}}
&< \frac{k_{p_{\ell+1}} + 1}{k_{p_{\ell+1}} + 2} + \tfrac{1}{2} \frac{(k_{p_{\ell+1}} + 1)^2}{(k_{p_{\ell+1}} + 2)^2 \bigl( 1 - \frac{k_{p_{\ell+1}} + 1}{k_{p_{\ell+1}} + 2} \bigr)} \notag \\
&< 1 + \tfrac{1}{2} \frac{(k_{p_{\ell+1}} + 1)^2}{k_{p_{\ell+1}} + 2} \notag \\
&< 1 + \tfrac{1}{2} (k_{p_{\ell+1}} + 1) \notag \\
&= \tfrac{1}{2} (k_{p_{\ell+1}} + 3). \label{eqn special two}
\end{align}
Now we have for $0 \le j < q_{\ell+1} - 2$ (in case $q_{\ell+1} \ge 2$),
\begin{align}
\sum_{i=0}^j \log(1 - \theta_{p_\ell+i})
&> -\sum_{i=0}^j \bigl( \theta_{p_\ell+i} + \tfrac{1}{2} \frac{\theta_{p_\ell+i}^2}{1 - \theta_{p_\ell+i}} \bigr), \text{ by \eqref{eqn log(1-t)},} \notag \\
&> - \sum_{i=0}^j \frac{1}{q_{\ell+1} - i - 1}, \text{ by \eqref{eqn special one},} \notag \\
&= -\sum_{i = q_{\ell+1} - (j + 1)}^{q_{\ell+1} - 1} \frac{1}{i} \label{eqn special three} \\
&> - \log \frac{q_{\ell+1} - 1}{q_{\ell+1} - 1 - (j + 1)}, \text{ by \eqref{eqn partial harmonic series},} \notag \\
&= \log \bigl(1 - \frac{j+1}{q_{\ell+1} - 1} \bigr). \label{eqn special 4}
\end{align}
For $j = q_{\ell+1} - 2$ (still in case $q_{\ell+1} \ge 2$), we have
\begin{align}
\sum_{i=0}^{q_{\ell+1} - 2} \log(1 - \theta_{p_\ell+i})
&> -\sum_{i = 1}^{q_{\ell+1} - 1} \frac{1}{i}, \text{ by \eqref{eqn special three}} \notag \\
&\ge -1 - \log(q_{\ell+1} - 1), \text{ by \eqref{eqn partial harmonic series}} \label{eqn lower bound 1.5} \\
&> -1 - \log q_{\ell+1}. \label{eqn lower bound two}
\end{align}
Now suppose $j = q_{\ell+1} - 1$. If $q_{\ell+1} \ge 2$, by \eqref{eqn lower bound 1.5}, \eqref{eqn lower bound two} and \eqref{eqn special two} we have
\begin{align}
\sum_{i=0}^{q_{\ell+1} - 1} \log(1 - \theta_{p_\ell+i})
&> -1 - \log (q_{\ell+1} - 1) - \tfrac{1}{2} (k_{p_{\ell+1}} + 3) \label{eqn special 5} \\
&> -1 - \log q_{\ell+1} - \tfrac{1}{2} (k_{p_{\ell+1}} + 3) \label{eqn special 6}.
\end{align}
If $q_{\ell+1} = 1$ then $p_\ell = p_{\ell+1} - 1$, and we have
\begin{align*}
\sum_{i=0}^{q_{\ell+1} - 1} \log(1 - \theta_{p_\ell+i})
&= \log(1 - \theta_{p_{\ell+1} - 1}) \\
&> -(\theta_{p_{\ell+1} - 1} + \tfrac{1}{2} \frac{\theta_{p_{\ell+1} - 1}^2}{1 - \theta_{p_{\ell+1} - 1}}) \\
&> -\tfrac{1}{2}(k_{p_{\ell+1}} + 3), \text{ by \eqref{eqn special two}}.
\end{align*}
Now we exponentiate: for $j < q_{\ell+1} - 2$ we have
\begin{equation*}
\prod_{i=0}^j (1 - \theta_{p_\ell + i})
> 1 - \frac{j+1}{q_{\ell+1} - 1}, \text{ by \eqref{eqn special 4}};
\end{equation*}
for $j = q_{\ell+1} - 2$ we have, by \eqref{eqn lower bound 1.5} and \eqref{eqn lower bound two},
\begin{align*}
\prod_{i=0}^{q_{\ell+1} - 2} (1 - \theta_{p_\ell + i})
&> \frac{e^{-1}}{q_{\ell+1} - 1} \\
&> \frac{e^{-1}}{q_{\ell+1}};
\end{align*}
and finally,
\begin{align*}
\prod_{i=0}^{q_{\ell+1}-1} (1 - \theta_{p_\ell+i})
&> \begin{cases}
\frac{e^{-1}}{(q_{\ell+1} - 1)} e^{-\tfrac{1}{2} (k_{p_{\ell+1}} + 3)}, &\text{ if } q_{\ell+1} > 1, \text{ by \eqref{eqn special 5}}, \\
e^{-\tfrac{1}{2} (k_{p_{\ell+1}} + 3)}, &\text{ if } q_{\ell+1} = 1.
\end{cases} \\
&> \frac{e^{-1}}{q_{\ell+1}} e^{-\tfrac{1}{2} (k_{p_{\ell+1}} + 3)}, \text{ by \eqref{eqn special 6}.}
\end{align*}
Now summing gives
\begin{align*}
\sum_{j=0}^{q_{\ell+1} - 1} \prod_{i=0}^j (1 - \theta_{p_{\ell}+i})
&> \sum_{j=0}^{q_{\ell+1} - 3} \bigl( 1 - \frac{j+1}{q_{\ell+1} - 1} \bigr) + \frac{e^{-1}}{q_{\ell+1} - 1} + \frac{e^{-1}}{q_{\ell+1} - 1} e^{-\frac{1}{2}(k_{p_{\ell+1}} + 3)} \\
&= q_{\ell+1} - 2 - \frac{1}{q_{\ell+1} - 1} \sum_{b=1}^{q_{\ell+1} - 2} b + \frac{e^{-1}}{q_{\ell+1} - 1} + \frac{e^{-1}}{q_{\ell+1} - 1} e^{-\frac{1}{2}(k_{p_{\ell+1}} + 3)} \\
&= \tfrac{1}{2}(q_{\ell+1} - 2) + \frac{e^{-1}}{q_{\ell+1} - 1} (1 + e^{-\frac{1}{2} (k_{p_{\ell+1}} + 3)}) \\
&> \tfrac{1}{2}(q_{\ell+1} - 2) + \frac{e^{-1}}{q_{\ell+1} - 1},
\end{align*}
in case $q_{\ell+1} \ge 2$, and in case $q_{\ell+1} = 1$,
\[
\sum_{j=0}^{q_{\ell+1} - 1} \prod_{i=0}^j (1 - \theta_{p_{\ell}+i})
= 1 - \theta_{p_\ell} > e^{-\frac{1}{2}(k_{p_{\ell+1}} + 3)}.
\]
Note that for $q \ge 2$ we have $\tfrac{1}{2}(q - 2) + \frac{e^{-1}}{q - 1} > \frac{1}{4} (q - 1)$, since this inequality reduces to $e^{-1} > \frac{1}{4}$ when $q = 2$, and since $\frac{1}{2}(q-2) \ge \frac{1}{4}(q-1)$ when $q > 2$. Thus we have
\[
\sum_{j=0}^{q_{\ell+1} - 1} \prod_{i=0}^j (1 - \theta_{p_{\ell}+i})
\ge
\tfrac{1}{4}(q_{\ell+1} - 1). \qedhere
\]
\end{proof}

\begin{Lemma} \label{lem upper and lower bounds}

Define $U$ and $L$ by
\begin{align*}
U &= \sum_{n=0}^\infty \tfrac{2}{3} q_{n+1} \prod_{\ell=0}^{n-1} \frac{2}{q_{\ell+1} + 2} \\
L &= \sum_{n=0}^\infty \tfrac{1}{4} (q_{n+1} - 1) e^{-n} \prod_{\ell=0}^{n-1} \frac{1}{q_{\ell+1}} e^{-\frac{1}{2}(k_{p_{\ell + 1}} + 3)}.
\end{align*}
Then $L \le 1 + \sum_{j=1}^\infty \prod_{i=0}^j (1 - \theta_i) \le U + \theta_0$.

\end{Lemma}

\begin{proof}
Let $a_i = 1 - \theta_i$. We have
\begin{align*}
\sum_{j=0}^\infty \prod_{i=0}^j a_i
&= \sum_{j=0}^{p_1 - 1} \prod_{i=0}^j a_i
  + \sum_{j=p_1}^{p_2 - 1} \prod_{i=0}^j a_i
    + \sum_{j=p_2}^{p_3 - 1} \prod_{i=0}^j a_i + \cdots \\
&= \sum_{j=0}^{q_1 - 1} \prod_{i=0}^j a_i
  + \sum_{j=0}^{q_2 - 1} \prod_{i=0}^{p_1 + j} a_i
    + \sum_{j=0}^{q_3 - 1} \prod_{i=0}^{p_2 + j} a_i + \cdots \\
&= \sum_{n=0}^\infty \sum_{j=0}^{q_{n+1} - 1} \prod_{i=0}^{p_n + j} a_i \\
&= \sum_{n=0}^\infty \Bigl( \prod_{i=0}^{p_n - 1} a_i \Bigr) \left( \sum_{j=0}^{q_{n+1} - 1} \prod_{i=0}^j a_{p_n + i} \right) \\
&= \sum_{n=0}^\infty \left( \prod_{\ell = 0}^{n - 1} \Bigl( \prod_{i=0}^{q_{\ell + 1} - 1} a_{p_\ell + i} \Bigr) \right) \left( \sum_{j=0}^{q_{n+1} - 1} \prod_{i=0}^j a_{p_n + i} \right).
\end{align*}
Thus we have
\[
\sum_{j=0}^\infty \prod_{i=0}^j (1 - \theta_i)
= \sum_{n=0}^\infty \left( \prod_{\ell = 0}^{n - 1} \Bigl( \prod_{i=0}^{q_{\ell + 1} - 1} (1 - \theta_{p_\ell + i}) \Bigr) \right) \left( \sum_{j=0}^{q_{n+1} - 1} \prod_{i=0}^j (1 - \theta_{p_n + i}) \right).
\]
Since $1 + \sum_{j=1}^\infty \prod_{i=0}^j (1 - \theta_i) = \theta_0 + \sum_{j=0}^\infty \prod_{i=0}^j (1 - \theta_i)$, we may use the estimates in Lemma \ref{lem first estimates} to finish the proof.
\end{proof}

Recall that the sequence $(q_j -1)_{j=1}^\infty$ gives the lengths of the strings of zeros between consecutive nonzero terms of the sequence $(k_i)_{i=1}^\infty$. We use Lemma \ref{lem upper and lower bounds} to give some general situations when $C^*(\Lambda)$ is, and is not, stable.

\begin{Theorem} \label{thm stability}

Let $\Lambda$ be defined by the sequence $(k_i)_{i=1}^\infty$ as in section \ref{section:one}. Let $(p_n)_{n=0}^\infty$ and $(q_n)_{n=1}^\infty$ be as above.

\begin{enumerate}

\item \label{thm stability one} If $(q_n)$ is a bounded sequence then $\sum_{j=0}^\infty \prod_{i=0}^j (1 - \theta_i) < \infty$, and hence $C^*(\Lambda)$ is not stable.

\item \label{thm stability two} If $q_n > 1 + e^{n-1} \prod_{i < n} q_i e^{\frac{1}{2}(k_{p_i} + 3)}$ for all $n$, then $\sum_{j=0}^\infty \prod_{i=0}^j (1 - \theta_i) = \infty$, and hence $C^*(\Lambda)$ is stable.

\end{enumerate}

\end{Theorem}

\begin{proof}
We first suppose that $(q_n)$ is a bounded sequence, say that $q_n \le C$ for all $n$. Then
\[
\sum_{n=0}^\infty \tfrac{2}{3} q_{n+1} \prod_{\ell=0}^{n-1} \frac{2}{q_{\ell+1} + 2}
\le \tfrac{2}{3} C \sum_{n=0}^\infty (\tfrac{2}{3})^n < \infty.
\]
By Lemma \ref{lem upper and lower bounds} it follows that $\sum_{j=0}^\infty \prod_{i=0}^j (1 - \theta_i) < \infty$.

Next we suppose that the condition in \eqref{thm stability two} holds. Then
\[
\sum_{n=0}^\infty \tfrac{1}{4} (q_{n+1} - 1) e^{-n} \prod_{\ell=0}^{n-1} \frac{1}{q_{\ell+1}} e^{-\frac{1}{2}(k_{p_{\ell + 1}} + 3)}
> \sum_{n=0}^\infty \tfrac{1}{4} = \infty.
\]
By Lemma \ref{lem upper and lower bounds} it follows that $\sum_{j=0}^\infty \prod_{i=0}^j (1 - \theta_i) = \infty$.
\end{proof}

A loose interpretation of Theorem \ref{thm stability} is that $C^*(\Lambda)$ is stable if and only if the sequence $(k_i)$ is rarely nonzero.


\begin{thebibliography}{xyx}

\bibitem{blackadar} B. Blackadar, \textit{$K$-Theory for Operator Algebras}, 2nd ed., Cambridge University Press, Cambridge, 1998.

\bibitem{bratkish} O. Bratteli and A. Kishimoto, \textit{Non-commutative spheres III}, Commun. Math. Phys. 147 (1992), 605-624.

\bibitem{dav} K.R. Davidson, \textit{C*-Algebras by Example}, Fields Institute Monographs, vol. 6, American Mathematical Society, Providence, RI, 1996.

\bibitem{drin} D. Drinen, \textit{Viewing AF-algebras as graph algebras}, Proc. American Mathematical Society 128 (2000) p. 1991-2000.

\bibitem{elpw} S. Echterhoff, W. L\"uck, N.C. Phillips and S. Walters, \textit{The structure of crossed products of irrational rotation algebras by finite subgroups of SL$_2(\IZ)$}, J. reine angew. Math. 639 (2010), 173-221.

\bibitem{effshen} E.G. Effros and C.L. Shen, \textit{Approximately finite C*-Algebras and continued fractions}, Indiana Univ. Math. J. 20 (1980) no. 2, p. 191 - 204.

\bibitem{evans} D.G. Evans, \textit{On higher rank graphs C*-algebras}, PhD thesis, Cardiff University, 2002.

\bibitem{evans2} D.G. Evans, \textit{On the K-theory of higher-rank graph C*-algebras}, New York J. Math, 14 (2008) p. 1-31.

\bibitem{evanssims} D.G. Evans and A. Sims, \textit{When is the Cuntz-Krieger algebra of a higher-rank graph approximately finite-dimensional?}, J. Funct. Anal. 263 (2012), no. 1, p. 183-215

\bibitem{hardywright} G.H. Hardy and E.M. Wright, \textit{An Introduction to the Theory of Numbers, 4th ed.}, Oxford University Press, London, 1975.

\bibitem{kumjian} A. Kumjian, \textit{An involutive automorphism of the Bunce-Deddens algebra}, C.R. Math. Rep. Acad. Sci. Canada, X (1988), no. 5, 217-218.

\bibitem{kumpas} A. Kumjian and D. Pask, \textit{Higher rank graph C*-algebras}, New York J. Math. 6 (2000), p. 1-20.

\bibitem{kumpasrae} A. Kumjian, D. Pask, and I. Raeburn, \textit{Cuntz-Krieger algebras of directed graphs}, Pacific J. Math.
184 (1998), 161-174.

\bibitem{laclarnes} M. Laca, N.S. Larsen, S. Neshveyev, \textit{On Bost-Connes type systems for number fields}, J. Number Theory 129 (2009), 325-338.

\bibitem{liren} X. Li and J. Renault, \textit{Cartan subalgebras in $C^*$-algebras. Existence and uniqueness}, Transactions AMS, Vol. 372(3) (2019), p.1985-2010.

\bibitem{mitscherthesis} I. Mitscher, \textit{Representing certain continued fraction AF algebras as $C^*$-algebras of categories of paths and non-AF groupoids}, PhD thesis, Arizona State University, 2020.

\bibitem{muhrenwil} P.S. Muhly, J.N. Renault, D. P. Williams, \textit{Equivalence and isomorphism for groupoid $C^*$-algebras}, J. Operator Theory 17 (1987), no. 1, p. 3-22.

\bibitem{nesh} S. Neshveyev, \textit{KMS States on the $C^*$-algebras of non-principal groupoids}, J. Operator Theory 70 (2013), no. 2, p. 513-530.

\bibitem{raeb} I. Raeburn, \textit{Graph Algebras}, CBMS Regional Conference Series in Mathematics, vol. 103, Published for the Conference Board of the Mathematical Science, Washington, DC; by the Amerian Mathematical Society, Providence, RI, 2005.

\bibitem{raesimyee} I. Raeburn, A. Sims and T. Yeend, \textit{The $C^*$-algebras of finitely aligned higher rank graphs}, J. Functional Analysis 213 (2004), 206-240.

\bibitem{raeszy} I. Raeburn and W. Szyma\'nski, \textit{Cuntz-Krieger algebras of infinite graphs and matrices}, Transactions AMS 356, no. 1 (2003), 39-59.

\bibitem{ren} J. Renault, \textit{A groupoid approach to C*-algebras}, Lecture Notes in Mathematics vol. 793, Springer, Berlin, 1980.

\bibitem{ren2} J. Renault, \textit{The ideal structure of groupoid crossed product $C^*$-algebras}, J. Operator Theory 25 (1991), 3-36.

\bibitem{ren3} J. Renault, \textit{Cartan subalgebras in $C^*$-algebras}, Irish Math. Soc. Bulletin 61 (2008), 29-63.

\bibitem{rord} M. R\o rdam, F. Larsen, and N.J. Lausten, \textit{An Introduction to K-Theory for $C^*$-Algebras}, London Mathematical Society Student Texts 49, Cambridge University Press, Cambridge, 2000.

\bibitem{spiel1} J. Spielberg, \textit{Groupoids and $C^*$-algebras for categories of paths}, Transactions of the AMS, 366, Number 11 (2014), p. 5771-5819.

\bibitem{spiel2} J. Spielberg, \textit{Groupoids and $C^*$-algebras for left cancellative small categories}, Indiana Univ. Math. J., 69 no. 5 (2020), 1579-1626.

\bibitem{strvoi} S. Str\u atil\u a and D. Voiculescu, \textit{Representations of AF-algebras and of the group $U(\infty)$}, Lecture Notes in Mathematics Vol. 486, Springer-Verlag, Berlin, Heidelberg, New York, 1975.

\bibitem{tww} A. Tikuisis, S. White and W. Winter, \textit{Quasidiagonality of nuclear $C^*$-algebras}, Annals of Math. (2) 185 (2017), no. 1, 229-284.

\bibitem{tu} J-L. Tu, \textit{La conjecture de Baum-Connes pour les feuilletages moyennables}, K-Theory {\bf 17} (1999), 215--264.

\bibitem{tyler} J. Tyler, \textit{Every AF-algebra is Morita equivalent to a graph algebra}, Bull. Austral. Math. Soc. 69 (2004), 237-240.

\bibitem{wall} H.S. Wall, \textit{Analytic Theory of Continued Fractions}, Dover, Mineola, New York, 2018.

\bibitem{will} D.P. Williams, \textit{A Toolkit for Groupoid $C^*$-algebras}, Mathematical Surveys and Monographs 241, American Mathematical Society, Providence, RI, 2019.

\bibitem{winzac} W. Winter and J. Zacharias, \textit{The nuclear dimension of C*-algebras}, Advances in Mathematics 224 (2010), no. 2, p. 461-498.

\end{thebibliography}
\end{document}